\pdfoutput=1

\documentclass[newpagesize]{crscat}
\svnInfo $Id: crscat-I.tex 353 2009-03-04 15:07:18Z vtl $
\renewcommand{\draftname}{PREPRINT}

\title{Coarse categories I: foundations}
\authortrunggwdg

\begin{document}

\begin{abstract}
Following Roe and others (see, e.g., \cite{MR1451755}), we (re)develop coarse 
geometry from the foundations, taking a categorical point of view. In this 
paper, we concentrate on the discrete case in which topology plays no role. Our 
theory is particularly suited to the development of the \emph{Roe 
($C^*$-)algebras} $C^*(X)$ and their $K$-theory on the analytic side; we also 
hope that it will be of use in the more strictly geometric/algebraic setting of 
controlled topology and algebra. However, we leave these topics to future 
papers.

Crucial to our approach are nonunital coarse spaces, and what we call 
\emph{locally proper} maps; these are actually not new, being implicit in 
\cite{MR1988817}. Our \emph{coarse category} $\CATCrs$ is a generalization of 
the usual one: its objects are (possibly nonunital) coarse spaces and its 
morphisms are (locally proper) coarse maps modulo \emph{closeness}. $\CATCrs$ 
is considerably richer than the usual coarse category of unital coarse spaces 
and proper coarse maps. As such, it has all nonzero limits and all colimits 
(all of which are easily constructed). We examine various other categorical 
issues. For example, $\CATCrs$ does not have a terminal object, so we 
substitute a \emph{termination functor}. This functor will be important in the 
development of exponential objects (i.e., ``function spaces'') 
\cite{crscat-II}, and also leads to a notion of \emph{quotient coarse spaces}. 
To connect our methods with the standard methods, we also examine the 
relationship between $\CATCrs$ and the usual coarse category.

Finally we briefly discuss some basic examples and applications. Topics include 
\emph{metric coarse spaces}, \emph{continuous control} \cite{MR1277522},
metric and continuously controlled \emph{coarse simplices},
%FIXME
%\emph{coarse homotopy} \cite{MR1243611},
\emph{$\sigma$-coarse spaces} \cite{MR2225040}, and the relation between 
quotient coarse spaces and the $K$-theory of Roe algebras (which is of 
particular interest for continuously controlled coarse spaces).
\end{abstract}

\maketitle

%%%%%%%%%%%%%%%%%%%%%%%%%%%%%%%%%%%%%%%%%%%%%%%%%%%%%%%%%%%%%%%%%%%%%%%%%%%%%%%%

\section*{Introduction}

Coarse, or large-scale, geometry has long been studied in various guises, but 
most notably in the context of metric spaces. Most generically, a \emph{coarse 
space} is a space together with some kind of large-scale structure (e.g., a 
metric modulo \emph{quasi-isometry}; see Remark~\ref{rmk:lsLip-qisom}). A 
\emph{coarse map} between coarse spaces is then a map which respects this 
structure (e.g., a large-scale Lipschitz map). Since the small-scale (i.e., the 
topology) is ignored, one can typically take coarse spaces to be 
\emph{discrete}, replacing any nondiscrete space by some ``coarsely dense'' 
subset.

In recent decades, coarse ideas have played an important role in the study of 
infinite discrete groups using the methods of geometric group theory, 
especially in the work of Gromov and his followers (see, e.g., 
\cite{MR1253544}). The most basic example here is that if $\Gamma$ is a 
finitely generated group, then the word length metric on $\Gamma$ is modulo 
quasi-isometry independent of the finite set of generators used in defining it.

Coarse ideas have also arisen in geometric topology, and more specifically 
controlled topology which primarily concerns itself with problems on the 
structure of manifolds. (We refer the reader to \cite{MR1308714}*{Ch.~9} for a 
survey of the topic and for references.) In this setting, one is interested in 
``operations'' (e.g., homotopies, surgeries) on spaces which respect some 
large-scale structure, i.e., are \emph{controlled}. As before, one may take the 
large-scale structure to be given by a metric (i.e., \emph{bounded control}). 
However, it is often more convenient to work with a coarser large-scale 
structure which is defined in purely topological terms (i.e., \emph{continuous 
control}; see \S\ref{subsect:cts-ctl}).

Controlled topology parallels the more classical theory for compact manifolds 
which relies on the use of algebraic invariants (e.g., algebraic $K$-theory). 
In controlled topology, one gets controlled versions of those invariants (in, 
e.g., \emph{bounded} and \emph{continuously controlled $K$-theory} 
\cites{MR1000388, MR802790, MR1277522}; see also \cite{MR1208729}). By 
considering the fundamental group of a space, a key object of study in the 
study of homotopy invariants (e.g., the Novikov Conjecture on higher 
signatures), many of the problems of geometric topology are related back to 
geometric group theory.

On the analytic side, to any coarse space $X$, Roe has associated a 
$C^*$-algebra $C^*(X)$ (the \emph{Roe algebra} of $X$), as well as various 
``(co)homology'' groups, e.g., \emph{coarse $K$-homology} $KX_\grstar(X)$. (For 
a good overview of this and the following, see \cite{MR1399087}.) On the other 
hand, one can also take the $K$-theory of $C^*(X)$; the Coarse Baum--Connes 
Conjecture is that a certain assembly map $KX_\grstar(X) \to K_\grstar(C^*(X))$ 
is an isomorphism, at least for suitably nice $X$.

The $K$-theory of Roe algebras arises in the index theory of elliptic operators 
on noncompact manifolds (on compact manifolds, the Roe algebra is just the 
compact operators and the results specialize to classical index theory). 
Indeed, historically it was the study of index theory on noncompact manifolds 
which led Roe to coarse geometry (see \cites{MR918459, MR1399087}), and not the 
other way around. In this way, the analytic approaches to the Novikov 
Conjecture (starting with the work of Lusztig \cite{MR0322889}) are again 
related to coarse geometry. (See \cite{MR1388295} for a nice survey of the 
different approaches to the Novikov Conjecture.)

\subsection*{Roe's coarse geometry}

After originally developing coarse geometry in the metric context 
\cite{MR1147350}, Roe (and his collaborators) realized that one can define an 
abstract notion of \emph{coarse space}, just as in small-scale geometry one has 
abstract \emph{topological spaces}. Just as the passage from metric space to 
topological space forgets large-scale (metric) information, the passage from 
metric space to coarse space should forget small-scale information. But an 
abstract coarse space retains enough structure to perform the large-scale 
constructions which were previously done in the metric context (e.g., construct 
the Roe algebras, coarse $K$-homology, etc.).

A \emph{coarse space} is a set $X$ together with a \emph{coarse structure}, 
which is a collection $\calE_X$ of subsets of $X^{\cross 2} \defeq X \cross X$ 
(called the \emph{entourages} of $X$) satisfying various axioms. When $X$ is a 
(proper) metric space, $\calE_X$ consists of the subsets $E \subseteq X^{\cross 
2}$ such that
\[
    \sup \set{d_X(x,x') \suchthat (x,x') \in E} < \infty.
\]
A subset $K \subseteq X$ is \emph{bounded} if and only if $K^{\cross 2}$ is an 
entourage of $X$; when $X$ is a metric space, $K$ is bounded if and only if it 
is metrically bounded. If $X$ is a discrete set, one typically axiomatically 
insists that the bounded subsets of $X$ be finite (we call this the 
\emph{properness axiom}; see Definition~\ref{def:prop-ax}); more generally, if 
$X$ is a topological space, the bounded subsets are required to be compact.

A set map $f \from Y \to X$ is a \emph{coarse map} if $f$ is \emph{proper} in 
the sense that the inverse image of any bounded subset of $X$ is a bounded 
subset of $Y$ and if $f$ \emph{preserves entourages} in the sense that 
$f^{\cross 2}(F) \defeq (f \cross f)(F)$ is an entourage of $X$. In the metric 
case, $f$ is a coarse map if it is metrically proper and ``nonexpansive''.

There is an obvious notion of closeness for maps into a metric space: maps $f, 
f' \from Y \to X$ are \emph{close} if
\[
    \sup \set{d_X(f(y),f'(y)) \suchthat y \in Y} < \infty.
\]
This generalizes to the case when $X$ is a general coarse space: $f$, $f'$ are 
close if $(f \cross f')(1_Y)$ is an entourage of $X$, where $1_Y$ is the 
diagonal set $\set{(y,y) \suchthat y \in Y}$.

Roe's \emph{coarse category} has coarse spaces as objects, and closeness 
classes of coarse maps as morphisms. (A coarse map is a \emph{coarse 
equivalence} if it represents an isomorphism in the coarse category.) Coarse 
invariants are defined on this category, either as functions on the isomorphism 
classes of the coarse category (e.g., \emph{asymptotic dimension}) or as 
functors from the coarse category to some other category (e.g., \emph{coarse 
$K$-homology}). Though coarse invariants are the primary object of study in 
coarse geometry, the coarse category is rarely analyzed directly, and there is 
some confusion in the literature about what the coarse category is (some 
authors take its arrows to be actual coarse maps; we will call this the 
\emph{precoarse category}).

There is an obvious ``product coarse structure'' on the cartesian (set) product 
$X \cross Y$. The entourages are exactly the subsets of $(X \cross Y)^{\cross 
2}$ which project to entourages of $X$ and $Y$ in the obvious way. However, 
this is not (usually) a product in the coarse category: the projection maps are 
not proper, unless both $X$ and $Y$ are finite/compact. This problem already 
arises in the category of proper metric spaces and proper maps (modulo 
closeness).

\begin{UNremark}
The above does \emph{not} prove that $X$ and $Y$ do not have a product in the 
coarse category. Certain products (of infinite/noncompact coarse spaces) 
\emph{do} exist in the coarse category; indeed, there is an infinite space $X$
(namely the continuously controlled ray, or equivalently a countable set 
equipped with the \emph{terminal}, i.e., ``indiscrete'', coarse structure) such 
that the product of $X$ with every countable coarse space exists 
(Remark~\ref{rmk:term-unital-prod}). The above does not even prove that the 
``product coarse space'' $X \cross Y$, as defined above, is not a product of 
$X$ and $Y$ if equipped with suitable maps $X \cross Y \to X$ and $X \cross Y 
\to Y$ (not the set projections).
\end{UNremark}

\subsection*{Nonunital coarse spaces and locally proper maps}

Metric spaces always yield \emph{unital} coarse spaces, i.e., coarse spaces $X$ 
such that $1_X \defeq \set{(x,x) \suchthat x \in X}$ is an entourage. Though 
Roe defines nonunital coarse spaces, unitality is usually a standing 
assumption, presumably since nonunital coarse spaces have no obvious use.

\emph{The} major innovation of this paper is the following: We relax the 
requirement that coarse maps be proper, to a requirement that we call 
\emph{locally properness}. Local properness is not new: it is actually included 
in Bartels's definition of ``coarse map'' \cite{MR1988817}*{Def.~3.3}. When the 
domain is a unital coarse space, local properness is equivalent to (``global'') 
properness (Corollary~\ref{cor:loc-prop-uni}). However, when the domain is 
nonunital, we get many more coarse maps. Consequently, using nonunital coarse 
spaces, it becomes extremely easy to construct (nonzero) categorical products 
in the coarse category. Indeed, we can do much more.

\begin{UNexample}
Suppose $X'$ is a (closed) subspace of a proper metric space $X$, so that $X'$ 
is itself a coarse space. There is an obvious \emph{ideal} $\lAngle 1_{X'} 
\rAngle_X$ of $\calE_X$ generated by $1_{X'}$ (see Definition~\ref{def:ideal}). 
The coarse space $|X|_{\lAngle 1_{X'} \rAngle_X}$ with underlying set $X$ and 
coarse structure $\lAngle 1_{X'} \rAngle_X$ is nonunital, unless $X'$ is 
``coarsely dense'' in $X$.

Define a (set) map $p \from X \to X'$ which sends each $x \in X$ to a point 
$p(x)$ in $X'$ closest to $x$. Then $p$ is usually not proper, hence is not 
coarse as a map $X \to X'$. However, it \emph{is} locally proper and coarse (in 
our generalized sense) as a map $|X|_{\lAngle 1_{X'} \rAngle_X} \to X'$, and is 
actually a coarse equivalence. (We leave it to the reader to verify this, after 
locating the required definitions.)
\end{UNexample}

For simplicity as well as for philosophical reasons, we only consider 
\emph{discrete} coarse spaces; hence for us a map is (globally) proper if and 
only if the inverse image of any point is a finite set. If a map $f \from Y \to 
X$ between coarse spaces is proper, then $f^{\cross 2}$ is a proper map, and 
hence the restriction of $f^{\cross 2}$ to any entourage $F \subseteq Y^{\cross 
2}$ of $Y$ is proper. We take the latter as the definition of local properness: 
A map $f \from Y \to X$ between coarse spaces (not necessarily unital) is 
\emph{locally proper} if, for all entourages $F$ of $Y$, the restriction 
$f^{\cross 2} |_F \from F \to X^{\cross 2}$ is a proper map. There are a number 
of equivalent ways of defining local properness, the most intuitive of which is 
the following. For a nonunital coarse space, there is an obvious notion of 
\emph{unital subspace}; a map is locally proper if and only if the restriction 
to every unital subspace of its domain is a proper map 
(Corollary~\ref{cor:loc-prop-uni}).

When $X$ is nonunital, we must modify the the definition of closeness, lest the 
identity map on $X$ not be close to itself. We modify it in a simple way, now 
requiring that the domain also be a coarse space: Coarse maps $f, f' \from Y 
\to X$ (between possibly nonunital coarse spaces) are \emph{close} if $(f 
\cross f')(F)$ is an entourage of $X$ for every entourage $F$ of $Y$. After 
checking the usual things, we get our nonunital \emph{coarse category}, whose 
objects are (possibly nonunital) coarse spaces and whose arrows are closeness 
classes of (locally proper) coarse maps.

\begin{UNremark}
Emerson--Meyer have defined a notion of \emph{$\sigma$-coarse spaces}, coarse 
maps between such spaces, and an appropriate notion of closeness 
\cite{MR2225040}. A $\sigma$-coarse space is just the colimit of an increasing 
sequence of unital coarse spaces. In fact, we show that the (pre)coarse 
category of discrete $\sigma$-coarse spaces is equivalent to a subcategory of 
our (pre)coarse category consisting of the \emph{$\sigma$-unital coarse spaces} 
(we do not examine the situation when one allows $\sigma$-coarse spaces to have 
nontrivial topologies).
\end{UNremark}

\subsection*{Products, limits, etc.}

Let us see how to construct the product of coarse spaces $X$ and $Y$ in this 
category. We do so by putting a \emph{nonunital} coarse structure on the set $X 
\cross Y$. The entourages of the \emph{coarse product} $X \cross Y$ are the $G 
\subseteq (X \cross Y)^{\cross 2}$ such that:
\begin{enumerate}
\item the restricted projections $\pi_1 |_G, \pi_2 |_G \from G \to X \cross Y$ 
    are proper maps (this is the aforementioned properness axiom);
\item $\pi_X |_G \from G \to X^{\cross 2}$ and $\pi_Y |_G \from G \to Y^{\cross 
    2}$ are proper maps; and
\item $(\pi_X)^{\cross 2}(G)$ is an entourage of $X$ and $(\pi_Y)^{\cross 
    2}(G)$ is an entourage of $Y$.
\end{enumerate}
One can then check that this is a product in our nonunital coarse category 
(indeed, it is a product in our nonunital \emph{precoarse category}). We must 
emphasize that the coarse structure on the set product is crucial: If $\ast$ is 
a one-point coarse space, then $X \cross \ast \cong X$ as a set, but unless $X$ 
is bounded the coarse product $X \cross \ast$ is \emph{not} coarsely equivalent 
to $X$.

The above construction generalizes to all nonzero products (by nonzero product, 
we mean a product of a nonempty collection of spaces), including infinite 
products (Theorem~\ref{thm:PCrs-lim} and Proposition~\ref{prop:Crs-prod}). We 
will then proceed to examine equalizers in the nonunital coarse category, and 
discover that it has all equalizers of pairs of maps 
(Proposition~\ref{prop:Crs-equal}). A standard categorical corollary is that 
the nonunital coarse category has all nonzero (projective) limits 
(Theorem~\ref{thm:Crs-lim}). One can similarly analyze coproducts (i.e., sums 
or ``disjoint unions'') and coequalizers, and get that the nonunital coarse 
category has all colimits, i.e., inductive limits 
(Theorem~\ref{thm:Crs-colim}).

\subsection*{Terminal objects and quotients}

For set theoretic reasons, the coarse category does not have a terminal object. 
(As we shall see in \S\ref{subsect:rest-Crs}, one way of obtaining a terminal 
object is to restrict the cardinality of coarse spaces, though there is a 
better way to proceed. For most purposes, no such restriction is needed.) 
However, there is a plethora of coarse spaces which behave like terminal 
objects. The \emph{terminal coarse structure} on a set $X$ consists of the 
subsets of $X^{\cross 2}$ which are both ``row- and column-finite''; denote the 
resulting coarse space by $|X|_1$. A rather underappreciated fact about such 
coarse spaces is that, for any coarse space $Y$, \emph{all} coarse maps $Y \to 
|X|_1$ are close to one another. An immediate categorical consequence of this 
is that, assuming that any such coarse map exists, the product of $|X|_1$ and 
$Y$ in the (unital or nonunital) coarse category is just $Y$ itself 
(Proposition~\ref{prop:term-id}).

In the unital coarse category, $X \mapsto |X|_1$ is a functor. In the nonunital 
coarse category, one must replace $|X|_1$ with a different coarse space, 
denoted $\Terminate(X)$ (with $\Terminate(X) = |X|_1$ for $X$ unital), to 
obtain a functor. In an abelian category, one can define a quotient $X/f(Y)$ 
(for $f \from Y \to X$) as push-out $X \copro_Y 0$. This generalizes to any 
category with zero objects and push-out squares. We will see that in fact we 
can generalize this to the coarse setting, defining $X/[f](Y)$ to be the 
push-out $X \copro_Y \Terminate(Y)$ in the (nonunital) coarse category. 
(Indeed, one can do the same in the category of topological spaces, noting that 
there are two cases: ``$\Terminate(X)$'' is a one-point space if $X \neq 
\emptyset$ and the empty space otherwise.)

\subsection*{Applications}

Our development of coarse geometry is a strict generalization of Roe's, despite 
our assumption of discreteness (see \S\ref{sect:top-crs}). Most of the standard 
constructions in Roe's coarse geometry (such as those mentioned above) 
generalize easily to our nonunital, locally proper version. (Note, however, 
that our theory does not encompass what one may call, following the language of 
\cite{MR1817560}*{Ch.~12}, the ``uniform category'' in which both the coarse 
structure and the topology are important. For example, Roe's $C^*$-algebras 
$D^*(X)$, which are functorial for uniform maps, require a notion of 
\emph{topological coarse space}. We defer this task to \cite{crscat-III}; see 
Remark~\ref{rmk:top-crs-sp}.) However, we will refrain from fully developing 
these applications in this paper. For the purposes of this paper, we briefly 
examine some things enabled by our generalization.

Having examined the coarse category from the categorical point of view, many 
standard constructions from topology transfer easily over to the coarse 
setting. For example, one obtains a notion of coarse simplicial complex. Of 
course, it is easy to deal with finite complexes explicitly in the unital 
coarse category. However, one result of having \emph{all} colimits, including 
infinite ones, is that we actually obtain infinite coarse simplicial complexes. 
This should enable one to apply simplicial methods in coarse geometry.

%FIXME
%%%
\begin{comment}
Another application is to coarse homotopy. There are various notions of 
homotopy used in coarse geometry, e.g., the coarse homotopy of Higson and Roe 
\cite{MR1243611}. However, the standard description of coarse homotopy is not 
``categorical'', for the obvious reason that the standard, unital coarse 
category does not seem to have products in general. We rectify this, and 
reformulate coarse homotopy in much more familiar categorical terms: We find a 
coarse space $I$ such that, for (at most) countable coarse spaces $X$ and $Y$, 
a coarse homotopy of maps $Y \to X$ is exactly given by a coarse map $(h_t) 
\defeq Y \cross I \to X$. This $I$ comes equipped with coarse maps $\delta_j 
\from P \defeq |\setN|_1 \to I$, $j = 0, 1$ (or, indeed, $j \in \ccitvl{0,1}$), 
which allows one to recover the coarse maps ``at the endpoints'' via the 
compositions
\[
    Y \isoto Y \cross P \nameto{\smash{\id \cross \delta_j}} Y \cross I \to X.
\]
As a historical note, we mention that this description is motivated by the 
continuously controlled case in which the connection to topology is much more 
obvious.
\end{comment}
%%%

\subsection*{Notes on history and references}

The framework and terminology we use are essentially due to Roe and his 
collaborators (see \cites{MR1147350, MR1451755}, in particular). Since our 
development differs in various details and in the crucial concept of local 
properness, and for the sake of completeness, we provide a complete exposition 
from basic principles; other, more standard, expositions include 
\cites{MR1451755, MR1399087, MR1817560, MR2007488}. In the basics, we do not 
claim much originality and most of the results will be known to those familiar 
with coarse geometry. However, in the context of locally proper maps, we have 
found certain methods of proof (in particular, the use of 
Proposition~\ref{prop:prop}) to be particularly effective, and have emphasized 
the use of these methods. Thus our proofs of standard results may differ from 
the usual proofs.

We have endeavoured to provide reasonably thorough references. However, it is 
often unwieldy to provide complete data for things which have been generalized 
and refined over the years. In such cases, rather than providing references to 
the original definition and all the subsequent generalizations, we simply 
reference a work (often expository in nature) which provides the current 
standard definition; often, such definitions can be found in a number of 
places, such as the aforementioned standard expositions.

\subsection*{Organization}

The rest of this paper is organized into five (very unequal) sections:
\begin{description}
\item[\S\ref{sect:crs-geom}] We define our basic framework of coarse 
    structures, coarse spaces, and coarse maps, together with important results 
    on local properness, and push-forward and pull-back coarse structures.
\item[\S\ref{sect:Crs}] We consider the precoarse categories (and $\CATPCrs$ in 
    particular) and their properties; the arrows in these categories are actual 
    coarse maps. We discuss limits and colimits in these categories, as well as 
    the relation between the general category $\CATPCrs$ and the subcategories 
    of unital and/or connected coarse spaces.
\item[\S\ref{sect:Crs}] We discuss the relation of closeness on coarse maps, 
    establish basic properties of closeness, and consider the quotient coarse 
    categories ($\CATCrs$ in particular). We show that $\CATCrs$ has all 
    nonzero products and all equalizers, hence all nonzero limits. Similarly, 
    it has all coproducts and all coequalizers, hence all colimits. We define 
    the termination functor $\Terminate$, and examine some of its properties; 
    in particular, it provides ``identities'' for the product. We characterize 
    the monic arrows of $\CATCrs$ and show that $\CATCrs$ has categorical 
    images, and dually we do the same for epi arrows and coimages. We apply 
    $\Terminate$, together with push-outs, to define quotient coarse spaces. 
    Finally, we discuss ways to ``restrict'' $\CATCrs$ to obtain subcategories 
    with terminal objects.
\item[\S\ref{sect:top-crs}] We examine Roe's formalization of coarse geometry, 
    which allows coarse spaces to carry topologies, and the relation 
    between the Roe coarse category and ours. In particular, we discuss how, 
    given a ``proper coarse space'' (in the sense of Roe), one can functorially 
    obtain a (discrete) coarse space (in our sense). We show that this gives a 
    fully faithful functor from the Roe coarse category to $\CATCrs$.
\item[\S\ref{sect:ex-appl}] We give the basic examples of coarse spaces: those 
    which come from proper metric spaces, and those which come from 
    compactifications (i.e., continuously controlled coarse spaces). We define 
    corresponding metric and continuously controlled coarse simplices, and 
    indicate how one might then develop coarse simplicial theory.
%FIXME
%    We recall the notion of coarse homotopy, and show that it has a satisfying 
%    description (which parallels the usually description of homotopy in 
%    topology) in $\CATCrs$.
    We compare Emerson--Meyer's $\sigma$-coarse spaces to our nonunital coarse 
    spaces (in the discrete case). Finally, we briefly examine the relation 
    between quotients of coarse spaces, the $K$-theory of Roe algebras, and 
    Kasparov $K$-homology.
\end{description}

\subsection*{Acknowledgements}

This work has been greatly influenced by many people, too many to enumerate. 
However, I would like to specifically thank my thesis advisor John Roe for his 
guidance over the years, as well as Heath Emerson and Nick Wright for helpful 
discussions. I would also like to thank Marcelo Laca, John Phillips, and Ian 
Putnam for their support at the University of Victoria.

%%%%%%%%%%%%%%%%%%%%%%%%%%%%%%%%%%%%%%%%%%%%%%%%%%%%%%%%%%%%%%%%%%%%%%%%%%%%%%%%

\section{Foundations of coarse geometry}\label{sect:crs-geom}

Throughout this section, $X$, $Y$, and $Z$ will be sets (sometimes with extra 
structure), and $f \from Y \to X$ and $g \from Z \to Y$ will be (set) maps. We 
denote the restriction of $f$ to $T \subseteq Y$ by $f |_T \from T \to X$. When 
$f(Y) \subseteq S \subseteq X$, we denote the range restriction of $f$ to $S$ 
by $f |^S \from Y \to S$. Thus if $T \subseteq Y$ and $f(T) \subseteq S 
\subseteq X$, we have a restriction $f |_T^S \from T \to S$.

\subsection{\pdfalt{\maybeboldmath Operations on subsets of $X \cross X$}%
        {Operations on subsets of X x X}}

Much of the following can be developed in the more abstract context of 
groupoids \cite{MR1451755}, but we will refrain from doing so. The basic object 
in question is the pair groupoid $X^{\cross 2} \defeq X \cross X$. Recall that 
$X^{\cross 2}$ has object set $X$ and set of arrows $X \cross X$. The map $X 
\injto X^{\cross 2}$ is $x \mapsto (x,x) \eqdef 1_x$ for $x \in X$. The target 
and source maps are the projections $\pi_1, \pi_2 \from X^{\cross 2} \to X$, 
respectively. For $x,x',x'' \in X$, composition is given by $(x,x') \circ 
(x',x'') \defeq (x,x'')$ and the inverse by $(x,x')^{-1} \defeq (x',x)$. Any 
set map $f \from Y \to X$ induces a groupoid morphism
\[
    f^{\cross 2} \defeq f \cross f \from Y^{\cross 2} \to X^{\cross 2}
\]
which in turn induces a map $\powerset(Y^{\cross 2}) \to \powerset(X^{\cross 
2})$, again denoted $f^{\cross 2}$, between power sets.

\begin{definition}[\maybeboldmath operations on $\powerset(X^{\cross 2})$]
For $E, E' \in \powerset(X^{\cross 2})$:
\begin{enumerate}
\item (\emph{addition}) $E + E' \defeq E \union E'$;
\item (\emph{multiplication}) $E \circ E' \defeq \set{e \circ e' \suchthat 
    \text{$e \in E$, $e' \in E'$, and $\pi_2(e) = \pi_1(e')$}}$; and
\item (\emph{transpose}) $E^\transpose \defeq \set{e^{-1} \suchthat e \in E}$.
\end{enumerate}
For $S \subseteq X$, put $1_S \defeq \set{1_x \suchthat x \in S}$ (the 
\emph{local unit} over $S$, or simply \emph{unit} if $S = X$).
\end{definition}

\begin{proposition}
For all $E \in \powerset(X^{\cross 2})$,
\[
    E \circ 1_S = (\pi_2 |_E)^{-1}(S)
\quad\text{and}\quad
    1_S \circ E = (\pi_1 |_E)^{-1}(S)
\]
\end{proposition}

\begin{remark}
We will refrain from calling $E \circ E'$ a ``product'' to avoid confusion with 
cartesian/categorical products (e.g., $X \cross Y$). The transpose 
$E^\transpose$ is often called the ``inverse'' and denoted $E^{-1}$; we avoid 
this terminology and notation since it is somewhat deceptive (though, 
admittedly, also rather suggestive). Our units $1_X$ are usually denoted 
$\Delta_X$ (and called the diagonal, for obvious reasons); our terminology is 
more representative of the ``algebraic'' role played by the unit (and the local 
units) and avoids confusion with the (related) diagonal map $\Delta_X \from X 
\to X \cross X$ (where $X \cross X$ is the cartesian/categorical product).
\end{remark}

The operations of addition and composition make $\powerset(X^{\cross 2})$ into 
a semiring: addition is commutative with identity $\emptyset$, multiplication 
is associative with identity $1_X$, multiplication distributes over addition, 
and $\emptyset \circ E = \emptyset = E \circ \emptyset$ for all $E$. Addition 
is idempotent in that $E + E = E$ for all $E$. Each $1_S$ is idempotent with 
respect to multiplication, i.e., $1_S \circ 1_S = 1_S$ for all $S$. The 
transpose is involutive, i.e., $(E^\transpose)^\transpose = E$ for all $E$, 
and, moreover, $(E + E')^\transpose = E^\transpose + (E')^\transpose$, $(E 
\circ E')^\transpose = (E')^\transpose \circ E^\transpose$, and 
$(1_S)^\transpose = 1_S$, for all $E$, $E'$, and $S$.

\begin{definition}[neighbourhoods and supports]
For any $S \subseteq X$ and $E \in \powerset(X^{\cross 2})$, put
\begin{align*}
    E \cdot S & \defeq \pi_1(E \circ 1_S) = \pi_1( (\pi_2 |_E)^{-1}(S) )
        && \text{(\emph{left $E$-neighbourhood of $S$})}
\shortintertext{and}
    S \cdot E & \defeq \pi_2(1_S \circ E) = \pi_2( (\pi_1 |_E)^{-1}(S) )
        && \text{(\emph{right $E$-neighbourhood of $S$})}.
\end{align*}
We also call $E \cdot X = \pi_1(E)$ the \emph{left support} of $E$ and $X \cdot 
E = \pi_2(E)$ the \emph{right support} of $E$.
\end{definition}

\begin{remark}
The notations $N_E(S) \defeq E_S \defeq E[S] \defeq E \cdot S$ and $E^S \defeq 
S \cdot E$ are common, though our notation is hopefully more suggestive of the 
relation between $E \cdot S$, $S \cdot E$ and the previously defined 
operations.
\end{remark}

\begin{proposition}%\label{prop:crs-easy}
For all $E$, $E'$, and $S$:
\begin{align*}
    (E + E') \cdot S & = E \cdot S \union E' \cdot S &
    & \text{and} &
    S \cdot (E + E') & = S \cdot E \union S \cdot E'; \\
    E \circ 1_{E' \cdot S} & = E \circ E' \circ 1_S &
    & \text{and} &
    1_{S \cdot E} \circ E' & = 1_S \circ E \circ E'; \\
    (E \circ E') \cdot S & = E \cdot (E' \cdot S) &
    & \text{and} &
    S \cdot (E \circ E') & = (S \cdot E) \cdot E'; \quad\text{and} \\
    E^\transpose \cdot S & = S \cdot E &
    & \text{and} &
    S \cdot E^\transpose & = E \cdot S.
\end{align*}
\end{proposition}

$\powerset(X^{\cross 2})$ and $\powerset(X)$ are partially ordered by 
inclusion. All of the above ``operations'' are monotonic with respect to these 
partial orders.

\begin{proposition}%\label{prop:crs-monotonic}
If $E_1, E_2, E'_1, E'_2 \in \powerset(X^{\cross 2})$ with $E_1 \subseteq E_2$ 
and $E'_1 \subseteq E'_2$, and $S_1, S_2 \subseteq X$ with $S_1 \subseteq S_2$, 
then:
\begin{align*}
    E_1 + E'_1 & \subseteq E_2 + E'_2, & &&
    E_1 \circ E'_1 & \subseteq E_2 \circ E'_2, \\
    (E_1)^\transpose & \subseteq (E_2)^\transpose, & &&
    1_{S_1} & \subseteq 1_{S_2}, \\
    E_1 \cdot S_1 & \subseteq E_2 \cdot S_2, &
    & \quad\text{and} &
    S_1 \cdot E_1 & \subseteq S_2 \cdot E_2.
\end{align*}
\end{proposition}

\subsection{Discrete properness}

Since our coarse spaces are essentially discrete, for now we only discuss 
properness for maps between discrete sets.

\begin{definition}\label{def:proper}
A set map $f \from Y \to X$ is \emph{proper} if the inverse image $f^{-1}(K)$ 
of every finite subset $K \subseteq X$ is again finite.
\end{definition}

If $Y$ is itself a finite set, then any $f \from Y \to X$ is automatically 
proper. We will use the following facts extensively (compare 
\cite{MR979294}*{\S{}10.1 Prop.~5}).

\begin{proposition}\label{prop:prop}
Consider the composition of (set) maps $Z \nameto{\smash{g}} Y 
\nameto{\smash{f}} X$:
\begin{enumerate}
\item\label{prop:prop:I} If $f$ and $g$ are proper, then $f \circ g$ is proper.
\item\label{prop:prop:II} If $f \circ g$ is proper, then $g$ is proper.
\item\label{prop:prop:III} If $f \circ g$ is proper and $g$ is surjective, then 
    $f$ (and $g$) are proper.
\end{enumerate}
Note that injectivity implies properness.
\end{proposition}

In \enumref{prop:prop:III} above, the hypothesis that $g$ be surjective can be 
weakened to the requirement that $Y \setminus g(Z)$ be a finite set. All 
restrictions (including range restrictions) of proper maps are again proper.

\subsection{The properness axiom and coarse structures}

\begin{definition}\label{def:prop-ax}
A set $E \in \powerset(X^{\cross 2})$ satisfies the \emph{properness axiom} if 
the restricted target and source maps (i.e., projections) $\pi_1 |_E, \pi_2 |_E 
\from E \to X$ (or, also restricting the ranges, $\pi_1 |_E^{E \cdot X}$, 
$\pi_2 |_E^{X \cdot E}$) are proper set maps.
\end{definition}

\begin{proposition}\label{prop:prop-ax}
For $E \in \powerset(X^{\cross 2})$, the following are equivalent:
\begin{enumerate}
\item\label{prop:prop-ax:I} $E$ satisfies the properness axiom;
\item\label{prop:prop-ax:II} $E \circ 1_S$ and $1_S \circ E$ are finite for all 
    finite $S \subseteq X$; and
\item\label{prop:prop-ax:III} $E \circ E'$ and $E' \circ E$ are finite for all 
    finite $E' \in \powerset(X^{\cross 2})$.
\end{enumerate}
\end{proposition}

\begin{proof}
\enumref{prop:prop-ax:I} $\iff$ \enumref{prop:prop-ax:II}: Immediate from $E 
\circ 1_S = (\pi_2 |_E)^{-1}(S)$ (and symmetrically).

\enumref{prop:prop-ax:II} $\iff$ \enumref{prop:prop-ax:III}: The reverse 
implication is clear. For the forward implication, if $E'$ is finite then $E' 
\cdot X$ is finite, and hence so too is
\[
    E \circ E' = E \circ E' \circ 1_X = E \circ 1_{E' \cdot X}
\]
(and symmetrically).
\end{proof}

\begin{corollary}\label{cor:prop-ax}
If $E \in \powerset(X^{\cross 2})$ satisfies the properness axiom, then $E 
\cdot S$ and $S \cdot E$ are finite for all finite $S \subseteq X$.
\end{corollary}

\begin{proof}
Use $E \cdot S \defeq \pi_1(E \circ 1_S)$ (and similarly symmetrically).
\end{proof}

\begin{remark}
The converse of the above Corollary holds since we are only considering pair 
groupoids: observe that
\[
    (\pi_1 |_E)^{-1}(S) \subseteq S \cross S \cdot E
\]
(and similarly symmetrically). However, the converse does not hold in general 
for coarse structures on groupoids.
\end{remark}

\begin{proposition}[``algebraic'' operations and the properness axiom]%
        \label{prop:prop-ax-alg}
If $E, E' \in \powerset(X^{\cross 2})$ satisfy the properness axiom, then $E + 
E'$, $E \circ E'$, $E^\transpose$, and all subsets of $E$ satisfy the 
properness axiom. Also, all singletons $\set{e}$, $e \in X^{\cross 2}$, and 
hence all finite subsets of $X^{\cross 2}$ satisfy the properness axiom, as 
does the unit $1_X$.
\end{proposition}

\begin{proof}
Clear, except possibly for $E \circ E'$; for this, use 
Proposition~\ref{prop:prop-ax}\enumref{prop:prop-ax:III} (and associativity of 
multiplication).
\end{proof}

If $T$, $T'$ are matrices over $X^{\cross 2}$ (with values in some ring) are 
supported on $E, E' \in \powerset(X^{\cross 2})$ satisfying the properness 
axiom, then the product $TT'$ is defined and has support contained in $E \circ 
E'$. The passage to rings of matrices motivates the following.

\begin{definition}\label{def:crs-sp}
A \emph{coarse structure} on $X$ is a subset $\calE_X \subseteq
\powerset(X^{\cross 2})$ such that:
\begin{enumerate}
\item each $E \in \calE_X$ satisfies the properness axiom;
\item $\calE_X$ is closed under the operations of addition, multiplication, 
    transpose, and the taking of subsets (i.e., if $E \subseteq E'$ and $E' \in 
    \calE_X$, then $E \in \calE_X$); and
\item for all $x \in X$, the singleton $\set{1_x}$ is in $\calE_X$.
\end{enumerate}
A \emph{coarse space} is a set $X$ equipped with a coarse structure $\calE_X$ 
on $X$. We denote such a coarse space by $|X|_{\calE_X}$ or simply $X$. The 
elements of $\calE_X$ are called \emph{entourages} (of $\calE_X$ or of $X$).
\end{definition}

\begin{example}[finite sets]
If $X$ is a finite set, then any coarse structure on $X$ must be unital. 
Moreover, there is only one connected coarse structure on $X$, namely the power 
set of $X$.
\end{example}

Here are two natural coarse structures which exist on any set.

\begin{definition}
The \emph{initial coarse structure} $\calE_{|X|_0}$ on $X$ is the minimum 
coarse structure on $X$. The \emph{terminal coarse structure} $\calE_{|X|_1}$ 
on a set $X$ is the maximum coarse structure. (We denote the corresponding 
coarse spaces by $|X|_0$ and $|X|_1$, respectively.)
\end{definition}

By Proposition~\ref{prop:prop-ax-alg}, $\calE_{|X|_1}$ simply consists of all 
the $E \in \powerset(X^{\cross 2})$ which satisfy the properness axiom. (Thus 
``$E \in \calE_{|X|_1}$'' is a convenient abbreviation for ``$E \in 
\powerset(X^{\cross 2})$ satisfies the properness axiom''.) Any coarse 
structure on $X$ is a subset of the terminal coarse structure (and obviously 
contains the initial coarse structure). More generally, we have the following.

\begin{proposition}
The intersection of any collection of coarse structures on $X$ (possibly 
infinite) is again a coarse structure on $X$.
\end{proposition}

\begin{definition}
The coarse structure $\langle \calE' \rangle_X$ on $X$ \emph{generated} by a 
subset $\calE' \subseteq \calE_{|X|_1}$ is the minimum coarse structure on $X$ 
which contains $\calE'$.
\end{definition}

Of course, $\langle \calE' \rangle_X$  is just the intersection of all the 
coarse structures on $X$ containing $\calE'$. Note that $\calE_{|X|_0} = 
\langle \emptyset \rangle_X$; more concretely, $\calE_{|X|_0}$ consists of all 
the finite local units $1_S$, $S \subseteq X$ finite.

Given two subsets $\calE', \calE'' \subseteq \calE_{|X|_1}$ (e.g., coarse 
structures on $X$), denote
\[
    \langle \calE', \calE'' \rangle_X
        \defeq \langle \calE' \union \calE'' \rangle_X.
\]
Observe that $\langle \calE', \calE'' \rangle_X$ contains both $\langle \calE' 
\rangle_X$ and $\langle \calE'' \rangle_X$. We use similar notation given three 
or more subsets of $\calE_{|X|_1}$ and, more generally, if $\set{\calE'_j 
\suchthat j \in J}$ ($J$ some index set) is a collection of subsets of 
$\calE_{|X|_1}$,
\[
    \langle \calE'_j \suchthat j \in J \rangle_X
        \defeq \bigl\langle \textstyle\bigunion_{j \in J} \calE'_j
                    \bigr\rangle_X.
\]

One can describe the coarse structure generated by $\calE'$ rather more 
concretely.

\begin{proposition}
If $\calE' \subseteq \calE_{|X|_1}$ contains all the singletons $\set{1_x}$, $x 
\in X$, and is closed under the ``algebraic'' operations of addition, 
multiplication, and transpose, then
\[
    \langle \calE' \rangle_X = \set{E \subseteq E' \suchthat E' \in \calE'}.
\]
\end{proposition}

\begin{corollary}\label{cor:crs-struct-gen}
For any $\calE' \subseteq \calE_{|X|_1}$, $\langle \calE' \rangle_X$ consists 
of the all subsets of the ``algebraic closure'' of the union
\[
    \calE' \union \set{\set{1_x} \suchthat x \in X}.
\]
\end{corollary}

Subsets of coarse spaces are naturally coarse spaces.

\begin{definition}
Suppose $X$ is a coarse space and $X' \subseteq X$ is a subset. Then
\[
    \calE_{X'} \defeq \calE_X |_{X'}
        \defeq \calE_X \intersect \powerset((X')^{\cross 2})
\]
is a coarse structure on $X'$, called the \emph{subspace coarse structure}. 
Call $X' \subseteq X$ equipped with the subspace coarse structure a (coarse) 
\emph{subspace} of $X$.
\end{definition}

\begin{example}[discrete metric spaces]\label{ex:disc-met}
Let $(X,d)$ be a discrete, proper metric space. ($X$ is metrically 
\emph{proper} if closed balls of $X$ are compact; thus $X$ is discrete and 
proper if and only if every metrically bounded subset is finite.) The 
\emph{($d$-)metric coarse space} $|X|_d$ (or just $|X|$ for short) has as 
entourages the $E \in \calE_{|X|_1} \subseteq \powerset(X^{\cross 2})$ such 
that
\begin{equation}\label{ex:disc-met:eq}
    \sup \set{d(x,x') \suchthat (x,x') \in E} < \infty.
\end{equation}
We may also allow $d(x,x') = \infty$ (for $x \neq x'$). In the senses defined 
below, $|X|_d$ is always unital but is connected if and only if $d(x,x') < 
\infty$ always. If $X' \subseteq X$, then the metric coarse structure on $X'$ 
obtained from the restriction of the metric $d$ is just the subspace coarse 
structure $\calE_{|X|_d} |_{X'}$.
\end{example}

\subsection{Unitality and connectedness}

\begin{definition}\label{def:uni-conn}
A coarse structure $\calE_X$ on $X$ is \emph{unital} if $1_X \in \calE_X$. 
$\calE_X$ is \emph{connected} if every singleton $\set{e}$, $e \in X^{\cross 
2}$, is an entourage of $\calE_X$. A pair of points $x, x' \in X$ are 
\emph{connected} (with respect to $\calE_X$, or in the coarse space $X$) if 
$\set{(x,x')} \in \calE_X$.
\end{definition}

Most treatments of coarse geometry assume both unitality and connectedness, but 
we will assume neither. Connectedness is a relatively benign assumption (see 
\S\ref{subsect:PCrs-conn}), but \emph{not} assuming unitality will be 
particularly crucial.

\begin{remark}\label{rmk:gpd-conn}
Connectedness in the general coarse groupoid case is more complicated, since 
there may be multiple arrows having the same target and source, and since a 
groupoid itself may not be connected. Let $\calE_\calG$ be a coarse structure 
on a groupoid $\calG$. There are several possible notions of connectedness:
\begin{enumerate}
\item The obvious translation of the above to groupoids is to say that 
    $\calE_\calG$ is (locally) \emph{connected} if all singletons $\set{e}$ 
    ($e$ an arrow in the groupoid) are entourages of $\calE_\calG$.
\item $\calE_\calG$ is \emph{globally connected} if it is (locally) connected 
    and $\calG$ is connected as a groupoid.
\setcounter{tempcounter}{\value{enumi}}
\end{enumerate}
Objects $x$, $x'$ are \emph{connected} if \emph{all} arrows $e$ with target $x$ 
and source $x'$ yield entourages $\set{e}$. Then $\calE_\calG$ is (locally) 
connected if and only if \emph{all groupoid-connected} pairs of objects are 
connected, and globally connected if and only if \emph{all} pairs of objects 
are connected. But there is also a weaker notion of connectedness: $x$, $x'$ 
are \emph{weakly connected} if there is \emph{some} arrow $e$ with target $x$ 
and source $x'$ such that $\set{e}$ is an entourage.
\begin{enumerate}
\setcounter{enumi}{\value{tempcounter}}
\item $\calE_\calG$ is (locally) \emph{weakly connected} if all 
    groupoid-connected objects $x$, $x'$ are weakly connected.
\item $\calE_\calG$ is \emph{globally weakly connected} if it is (locally) 
    weakly connected and $\calG$ is connected as a groupoid.
\end{enumerate}
When $\calG$ is a pair groupoid (i.e., in our case), all the above notions 
coincide.
\end{remark}

\begin{proposition}
The terminal structure on any set $X$ is always unital and connected.
\end{proposition}

The intersection of unital coarse structures on a set $X$ is again unital, and 
similarly for connected coarse structures. Thus, for any $\calE' \subseteq 
\calE_{|X|_1}$, there are \emph{unital}, \emph{connected}, and \emph{connected 
unital} coarse structures on $X$ generated by $\calE'$. These can be described 
rather simply as
\begin{align*}
    \langle \calE' \rangle_X^\TXTuni
        & \defeq \bigl\langle \calE', \set{1_X} \bigr\rangle_X, \\
    \langle \calE' \rangle_X^\TXTconn
        & \defeq \bigl\langle \calE', \set{\set{e}
                \suchthat e \in X^{\cross 2}} \bigr\rangle_X,
\quad\text{and} \\
    \langle \calE' \rangle_X^\TXTconnuni
        & \defeq \bigl\langle \calE', \set{\set{e}
                \suchthat e \in X^{\cross 2}}, \set{1_X} \bigr\rangle_X,
\end{align*}
respectively.

\begin{definition}
The \emph{initial unital}, \emph{initial connected}, or \emph{initial connected 
unital coarse structure} on a set $X$ is the minimum coarse structure having 
the given property or properties, respectively. Denote the resulting coarse 
spaces by $|X|_0^\TXTuni$, $|X|_0^\TXTconn$, and $|X|_0^\TXTconnuni$, 
respectively.
\end{definition}

Clearly, $\calE_{|X|_0^\TXTuni} = \langle \set{1_X} \rangle_X$, so a coarse 
structure on $X$ is unital if and only if it contains $\calE_{|X|_0^\TXTuni}$. 
Similarly for the other properties. Note in particular that 
$\calE_{|X|_0^\TXTconn}$ consists of all the finite subsets of $X^{\cross 2}$.

\begin{remark}
In the groupoid case, the intersection of (locally) connected coarse structures 
on a given groupoid (in the sense of Remark~\ref{rmk:gpd-conn}) is again 
(locally) connected, and so all of the above holds. However, the intersection 
of weakly connected coarse structures (see the same Remark) may not be weakly 
connected, so there may not be a minimum weakly connected coarse structure on a 
given groupoid.
\end{remark}

We get an obvious notion of \emph{unital subspace} of any coarse space $X$. 
Clearly, $X' \subseteq X$ is a unital subspace if and only if $1_{X'}$ is an 
entourage of $X$. (Bartels calls the set of unital subspaces of $X$ the 
``domain of $\calE_X$'' \cite{MR1988817}*{Def.~3.2}.) Slightly more is true.

\begin{proposition}\label{prop:uni-subsp}
A subspace $X' \subseteq X$ of a coarse space $X$ is unital if and only if it 
occurs as the left (or right) support of some entourage of $X$.
\end{proposition}

\begin{proof}
If $X'$ is a unital subspace, then $X' = 1_{X'} \cdot X$. Conversely, if $X' = 
E \cdot X$ for some $E \in \calE_X$, then $1_{X'} \subseteq E \circ 
E^\transpose$ must be an entourage of $X$.
\end{proof}

Similarly, we get a notion of \emph{connected subspace} of $X$.

\begin{definition}
A (connected) \emph{component} of a coarse space $X$ is a maximal connected 
subspace of $X$.
\end{definition}

\begin{proposition}
Any coarse space $X$ is partitioned, as a set, into (a disjoint union of) its 
connected components.
\end{proposition}

We caution this decomposition is not necessarily a coproduct (in the coarse or 
precoarse category); see Corollary~\ref{cor:PCrs-fin-components}.

\subsection{Local properness, preservation, and coarse maps}

Recall that any (set) map $f \from Y \to X$ induces a map (indeed, a groupoid 
morphism) $f^{\cross 2} \from Y^{\cross 2} \to X^{\cross 2}$. Insisting that 
$f$ be proper is too strong a requirement when $Y$ is a nonunital coarse space. 
We thus introduce the following weaker requirement.

\begin{definition}\label{def:loc-prop}
A map $f \from Y \to X$ is \emph{locally proper for $F \in \calE_{|Y|_1}$} if 
$E \defeq f^{\cross 2}(F) \in \calE_{|X|_1}$ and the restriction $f^{\cross 2} 
|_F \from F \to X^{\cross 2}$ (or $f^{\cross 2} |_F^E$) is a proper (set) map. 
If $Y$ is a coarse space, then $f$ is \emph{locally proper} (for $\calE_Y$) if 
it is locally proper for all $F \in \calE_Y$.
\end{definition}

Local properness only requires a coarse structure on the domain, so we cannot 
say that the composition of locally proper maps is again locally proper. 
Nonetheless, separating local properness from the following will be useful when 
we define push-forward coarse structures (below).

\begin{definition}
Suppose $X$ is a coarse space. A map $f \from Y \to X$ \emph{preserves $F \in 
\calE_{|Y|_1}$} (with respect to $\calE_X$) if $E \defeq f^{\cross 2}(F) \in 
\calE_X$. If $Y$ is also a coarse space, then $f$ \emph{preserves entourages} 
(of $\calE_Y$, with respect to $\calE_X$) if $f$ preserves every $F \in 
\calE_Y$.
\end{definition}

\begin{definition}
Suppose $X$ is a coarse space. A map $f \from Y \to X$ is \emph{coarse for $F 
\in \calE_{|Y|_1}$} if $f$ is locally proper for $F$ and if $f$ preserves $F$. 
If $Y$ is also a coarse space, then $f$ is \emph{coarse} (or is a \emph{coarse 
map}) if $f$ is coarse for every $F \in \calE_Y$, i.e., if $f$ is locally 
proper and preserves entourages.
\end{definition}

\begin{remark}
The definition of ``coarse map'' is slightly redundant: If $f$ preserves 
entourages, then $f^{\cross 2}(F) \in \calE_X \subseteq \calE_{|X|_1}$ (which 
is one of the stipulations of local properness).
\end{remark}

\begin{proposition}\label{prop:crs-map-comp}
Consider a composition of $Z \nameto{\smash{g}} Y \nameto{\smash{f}} X$, where 
$X$ and $Y$ are coarse spaces. If $f$, $g$ are locally proper and $g$ preserves 
entourages, then $f \circ g$ is locally proper.
\end{proposition}

\begin{corollary}
A composition of coarse maps is again a coarse map.
\end{corollary}

\subsection{Basic properties of maps}

We first concentrate on local properness.

\begin{proposition}\label{prop:loc-prop}
Suppose $f \from Y \to X$ is a set map, $F \in \calE_{|Y|_1}$, and $E \defeq 
f^{\cross 2}(F)$. The following are equivalent:
\begin{enumerate}
\item\label{prop:loc-prop:I} $f$ is locally proper for $F$;
\item\label{prop:loc-prop:II} the restrictions $f |_{F \cdot Y}$ and $f |_{Y 
    \cdot F}$ (or $f |_{F \cdot Y}^{E \cdot X}$ and $f |_{Y \cdot F}^{X \cdot 
    E}$) of $f$ to the left and right supports of $F$ are proper; and
\item\label{prop:loc-prop:III} $f^{-1}(S) \cdot F$ and $F \cdot f^{-1}(S)$ are 
    finite for all finite $S \subseteq X$.
\end{enumerate}
\end{proposition}

\begin{proof}
(We omit proofs of the symmetric cases.) Consider the diagram
\[\begin{CD}
    F @>{f^{\cross 2} |_F^E}>> E \\
    @V{\pi_1 |_F^{F \cdot Y}}VV @V{\pi_1 |_E^{E \cdot X}}VV \\
    F \cdot Y @>{f |_{F \cdot Y}^{E \cdot X}}>> E \cdot X
\end{CD}\quad.\]
Observe the following: the above diagram commutes, i.e.,
\[
    \pi_1 |_E^{E \cdot X} \circ f^{\cross 2} |_F^E
        = f |_{F \cdot Y}^{E \cdot X} \circ \pi_1 |_F^{F \cdot Y};
\]
the two maps emanating from $F$ are surjections; and $\pi_1 |_F^{F \cdot Y}$ is 
proper. We now apply Proposition~\ref{prop:prop} several times.

\enumref{prop:loc-prop:I} \textimplies{} \enumref{prop:loc-prop:II}: $f^{\cross 
2} |_F^E$ and $\pi_1 |_E^{E \cdot X}$ are proper, so their composition is 
proper. Since $\pi_1 |_F^{F \cdot Y}$ is surjective, $f |_{F \cdot Y}$ is 
proper.

\enumref{prop:loc-prop:II} \textimplies{} \enumref{prop:loc-prop:III}: $f |_{F 
\cdot Y}$ and $\pi_1 |_F^{F \cdot Y}$ are proper, so their composition is 
proper. Then
\begin{equation}\label{prop:loc-prop:pf:eq}\begin{split}
    f^{-1}(S) \cdot F & = \pi_2( (\pi_1 |_F)^{-1}(f^{-1}(S)) ) \\
        & = \pi_2( (f_{F \cdot Y} \circ \pi_1 |_F^{F \cdot Y})^{-1}(S) )
\end{split}\end{equation}
is finite if $S \subseteq X$ is finite.

\enumref{prop:loc-prop:III} \textimplies{} \enumref{prop:loc-prop:I}: By 
\eqref{prop:loc-prop:pf:eq} and since $\pi_2 |_F$ is proper, the composition 
$f_{F \cdot Y}^{E \cdot X} \circ \pi_1 |_F^{F \cdot Y}$ is proper. Hence 
$f^{\cross 2} |_F^E$ is proper and, since $f^{\cross 2} |_F^E$ is surjective, 
so is $\pi_1 |_E^{E \cdot X}$.
\end{proof}

\begin{corollary}
If a set map $f \from Y \to X$ is globally proper, then $f$ is locally proper 
for any $F \in \calE_{|Y|_1}$ (so $f$ is locally proper for any coarse 
structure on $Y$).
\end{corollary}

\begin{proof}
This follows from \enumref{prop:loc-prop:III} and Corollary~\ref{cor:prop-ax}.
\end{proof}

\begin{corollary}
If $X$ is a coarse space and $X' \subseteq X$ is a subspace, then the inclusion 
of $X'$ into $X$ is a coarse map. Thus the restriction of any coarse map to a 
subspace is a coarse map.
\end{corollary}

\begin{proof}
By definition of the subspace coarse structure, the inclusion map preserves 
entourages. The inclusion map is injective, hence (globally) proper, hence 
locally proper.
\end{proof}

\begin{corollary}\label{cor:loc-prop-uni}
Suppose $Y$ is a coarse space. A map $f \from Y \to X$ is locally proper if and 
only if the restriction of $f$ to every unital subspace of $Y$ is proper. Thus, 
for $Y$ unital, $f$ is locally proper if and only if $f$ is globally proper.
\end{corollary}

\begin{proof}
This follows from \enumref{prop:loc-prop:II} and 
Proposition~\ref{prop:uni-subsp}.
\end{proof}

For (discrete) unital coarse spaces, our notion of ``coarse map'' is just the 
classical notion. It also follows that local properness of a map $f \from Y \to 
X$ is a property which can be defined in terms of the unital subspaces of the 
coarse structure on $Y$. In particular, if $f$ is locally proper, then $f$ 
would also be locally proper for any coarse structure on $Y$ (possibly larger 
than $\calE_Y$) with the same unital subspaces.

\begin{remark}
One may take the \emph{definition} of local properness to be the 
characterization of the above Corollary, i.e., define $f \from Y \to X$ to be 
locally proper if $f$ is proper on every unital subspace of $Y$ (perhaps 
``unital properness'' would be a more apt term). Indeed, this is the form in 
which local properness appears in Bartels's definition of ``coarse map'' 
\cite{MR1988817}*{Def.~3.3}, and hence (modulo our coarse spaces not carrying 
topologies) our definition of ``coarse map'' is the same as Bartels's. More 
generally, one could remove coarse structures entirely, and define local 
properness for sets equipped with families of supports (i.e., of unital 
subspaces). However, we will not do so since we are mainly concerned with 
coarse maps, for which Definition~\ref{def:loc-prop} is most convenient.
\end{remark}

\begin{corollary}
Coarse maps send unital subspaces to unital subspaces, i.e., if $f \from Y \to 
X$ is a coarse map and $Y' \subseteq Y$ is a unital subspace, then the image 
$f(Y') \subseteq X$ is a unital subspace.
\end{corollary}

\begin{proposition}[``algebraic'' operations and local properness]%
        \label{prop:loc-prop-alg}
If $f \from Y \to X$ is locally proper for $F, F' \in \calE_{|Y|_1}$, then $f$ 
is locally proper for $F + F'$, $F \circ F'$, $F^\transpose$, and all subsets 
of $F$. Also, $f$ is locally proper for all singletons $\set{e}$, $e \in 
Y^{\cross 2}$, hence is locally proper for $\calE_{|Y|_0^\TXTconn} \supseteq 
\calE_{|Y|_0}$. (However, $f$ is locally proper for the unit $1_Y$ if and only 
if $f$ is globally proper.)
\end{proposition}

\begin{proof}
The only nontrivial assertion is that $f$ is locally proper for $F \circ F'$. 
By assumption, $f^{\cross 2}(F), f^{\cross 2}(F') \in \calE_{|X|_1}$ and, since
\[
    f^{\cross 2}(F \circ F') \subseteq f^{\cross 2}(F) \circ f^{\cross 2}(F'),
\]
$f^{\cross 2}(F \circ F')$ also satisfies the properness axiom, by 
Proposition~\ref{prop:prop-ax-alg}. We have a commutative diagram
\[\begin{CD}
    F \circ F' @>{f^{\cross 2} |_{F \circ F'}}>> X^{\cross 2} \\
    @V{\pi_1 |_{F \circ F'}^{(F \circ F') \cdot Y}}VV @V{\pi_1}VV \\
    (F \circ F') \cdot Y @>{f |_{(F \circ F') \cdot Y}}>> X
\end{CD}\quad.\]
By the same Proposition, $F \circ F' \in \calE_{|Y|_1}$, so $\pi_1 |_{F \circ 
F'}^{(F \circ F') \cdot Y}$ is proper. Since $(F \circ F') \cdot Y \subseteq F 
\cdot Y$ and $f |_{F \cdot Y}$ is proper by 
Proposition~\ref{prop:loc-prop}\enumref{prop:loc-prop:II}, $f |_{(F \circ F') 
\cdot Y}$ is proper. Hence the composition
\[
    f |_{(F \circ F') \cdot Y} \circ \pi_1 |_{F \circ F'}^{(F \circ F') \cdot Y}
        = \pi_1 \circ f^{\cross 2} |_{F \circ F'}
\]
is proper, so $f^{\cross 2} |_{F \circ F'}$ is proper by 
Proposition~\ref{prop:prop}\enumref{prop:prop:II}.
\end{proof}

\begin{corollary}\label{cor:loc-prop-gen}
If $f \from Y \to X$ is locally proper for all $F \in \calE' \subseteq 
\calE_{|Y|_1}$, then $f$ is locally proper for the coarse structure $\langle 
\calE' \rangle_Y$ on $Y$ generated by $\calE'$ (and for the connected coarse 
structure $\langle \calE' \rangle_Y^\TXTconn$ generated by $\calE'$).
\end{corollary}

\begin{proof}
This follows immediately from the Proposition and 
Corollary~\ref{cor:crs-struct-gen}.
\end{proof}

The same evidently does not hold for the unital (or connected unital) coarse 
structure generated by $\calE'$.

We now state some parallel results for preservation of entourages. Combining 
these with the above results for local properness, we get parallel results for 
coarseness of maps.

\begin{proposition}\label{``algebraic'' operations and preservation}%
        %\label{prop:preserve-alg}
Suppose $X$ is a coarse space. If $f \from Y \to X$ preserves $F, F' \in 
\calE_{|Y|_1}$, then $f$ preserves $F + F'$, $F \circ F'$, $F^\transpose$, and 
all subsets of $F$. Also, $f$ preserves all singletons $\set{1_y}$, $y \in Y$ 
(hence preserves $\calE_{|Y|_0}$); if $X$ is connected, $f$ preserves all 
singletons $\set{e}$, $e \in Y^{\cross 2}$ (hence preserves 
$\calE_{|Y|_0^\TXTconn}$); and if $X$ is unital, $f$ preserves $1_Y$.
\end{proposition}

\begin{proof}
The only (slightly) nontrivial one is $F \circ F'$, for which ones uses
\[
    f^{\cross 2}(F \circ F') \subseteq f^{\cross 2}(F) \circ f^{\cross 2}(F').
\]
\end{proof}

\begin{corollary}%\label{cor:preserve-gen}
Suppose $X$ is a coarse space. If $f \from Y \to X$ preserves $\calE' \subseteq 
\calE_{|Y|_1}$, then $f$ preserves the coarse structure $\langle \calE' 
\rangle_Y$ on $Y$ generated by $\calE'$. (If $X$ is also connected, then $f$ 
preserves $\langle \calE' \rangle_Y^\TXTconn$; if $X$ is unital, then $f$ 
preserves $\langle \calE' \rangle_Y^\TXTuni$; if $X$ is both, then $f$ 
preserves $\langle \calE' \rangle_Y^\TXTconnuni$.)
\end{corollary}

\begin{proposition}[``algebraic'' operations and coarseness]%
        \label{prop:coarse-alg}
Suppose $X$ is a coarse space. If $f \from Y \to X$ is coarse for $F, F' \in 
\calE_{|Y|_1}$, then $f$ is coarse for $F + F'$, $F \circ F'$, $F^\transpose$, 
and all subsets of $F$. Also, $f$ is coarse for all singletons $\set{1_y}$, $y 
\in Y$ (hence is coarse for $\calE_{|Y|_0}$); if $X$ is connected, $f$ is 
coarse for all singletons $\set{e}$, $e \in Y^{\cross 2}$ (hence is coarse for 
$\calE_{|Y|_0^\TXTconn}$); and if $X$ is unital and $f$ is proper, $f$ is 
coarse for $1_Y$.
\end{proposition}

\begin{corollary}\label{cor:coarse-gen}
Suppose $X$ and $Y$ are coarse spaces, $\calE' \subseteq \calE_{|Y|_1}$, and $f 
\from Y \to X$ is a set map.
\begin{enumerate}
\item If $\calE_Y = \langle \calE' \rangle_Y$, then $f$ is a coarse map if and 
    only if $f$ is coarse for all $F \in \calE'$.
\item If $\calE_Y = \langle \calE' \rangle_Y^\TXTconn$, then $f$ is a coarse 
    map if and only if $f$ is coarse for all $F \in \calE'$ and all $\set{e}$, 
    $e \in Y^{\cross 2}$.
\item If $\calE_Y = \langle \calE' \rangle_Y^\TXTuni$, then $f$ is a coarse map 
    if and only if $f$ is proper and $f$ is coarse for (or preserves) all $F 
    \in \calE'$.
\item If $\calE_Y = \langle \calE' \rangle_Y^\TXTconnuni$, then $f$ is a coarse 
    map if and only if $f$ is proper and $f$ is coarse for (or preserves) all 
    $F \in \calE'$ and all $\set{e}$, $e \in Y^{\cross 2}$.
\end{enumerate}
Note that requiring that $f$ be coarse for all $\set{e}$, $e \in Y^{\cross 2}$, 
is equivalent to requiring $f(y)$, $f(y')$ be connected for all $y, y' \in Y$.
\end{corollary}

If $f, f' \from Y \to X$ are (globally) proper maps, then certainly $f \cross 
f' \from Y \cross Y \to X \cross X$ is proper. The same also holds locally, and 
this will be essential later.

\begin{proposition}\label{prop:loc-prop-prod}
If (set) maps $f, f' \from Y \to X$ are locally proper for $F \in
\calE_{|Y|_1}$, then:
\begin{enumerate}
\item $E \defeq (f \cross f')(F) \subseteq X^{\cross 2}$ satisfies the 
    properness axiom; and
\item the restriction $(f \cross f') |_F^E \from F \to E$ is a proper map.
\end{enumerate}
\end{proposition}

\begin{proof}
Fix $F \in \calE_{|Y|_1}$, put $E \defeq (f \cross f')(F)$, and consider the 
commutative diagram
\[\begin{CD}
    F @>{(f \cross f') |_F^E}>> E \\
    @V{\pi_1 |_F^{F \cdot Y}}VV @V{\pi_1 |_E}VV \\
    F \cdot Y @>{f |_{F \cdot Y}}>> X
\end{CD}\quad.\]
The composition along the left and bottom is proper, and thus so is composition 
along the top and right. Consequently, $(f \cross f') |_F^E$ is proper. Since 
$(f \cross f') |_F^E$ is surjective, $\pi_1 |_E$ is proper and similarly for 
$\pi_2 |_E$.
\end{proof}

We have the following ``very'' local analogue of Proposition~\ref{prop:prop}. 
For a more general analogue, we will need push-forward coarse structures.

\begin{proposition}\label{prop:loc-prop-for-comp}
Consider the composition of (set) maps $Z \nameto{\smash{g}} Y 
\nameto{\smash{f}} X$, supposing that $G \in \calE_{|Z|_1}$ and putting $F 
\defeq g^{\cross 2}(G)$:
\begin{enumerate}
\item\label{prop:loc-prop-for-comp-I} If $g$ is locally proper for $G$ and $f$ 
    is locally proper for $F$, then $f \circ g$ is locally proper for $G$.
\item\label{prop:loc-prop-for-comp-II} If $f \circ g$ is locally proper for 
    $G$, then $g$ is locally proper for $G$.
\item\label{prop:loc-prop-for-comp-III} If $f \circ g$ is locally proper for 
    $G$, then $f$ is locally proper for $F$.
\end{enumerate}
\end{proposition}

\begin{proof}
Put $E \defeq f^{\cross 2}(F)$. We apply Proposition~\ref{prop:prop} to the 
commutative diagram
\[\begin{CD}
G @>{g^{\cross 2} |_G^F}>> F @>{f^{\cross 2} |_F^E}>> E \\
@V{\pi_1 |_G}VV @V{\pi_1 |_F}VV @V{\pi_1 |_E}VV \\
Z @>{g}>> Y @>{F}>> X
\end{CD}\quad.\]
\enumref{prop:loc-prop-for-comp-I} is clear. For 
\enumref{prop:loc-prop-for-comp-II} and \enumref{prop:loc-prop-for-comp-III}: 
If $f \circ g$ is locally proper for $G$, then $\pi_1 |_E$ and
\[
    (f \circ g)^{\cross 2} |_G^E
        = f^{\cross 2} |_F^E \circ g^{\cross 2} |_G^F
\]
are proper. By the latter, $g^{\cross 2} |_G^F$ is proper. $g^{\cross 2} |_G^F$ 
is surjective, so $f^{\cross 2} |_F^E$ is also proper. Then
\[
    \pi_1 |_E \circ f^{\cross 2} |_F^E = f \circ \pi_1 |_F
\]
is proper, so $\pi_1 |_F$ is proper.
\end{proof}

\subsection{Pull-back and push-forward coarse structures}

\begin{definition}
Suppose $X$ is a coarse space. The \emph{pull-back coarse structure} (of 
$\calE_X$) on $Y$ along (a set map) $f \from Y \to X$ is
\[
    f^* \calE_X \defeq \set{F \in \calE_{|Y|_1} \suchthat
            \text{$f$ is coarse for $F$}}.
\]
\end{definition}

By Proposition~\ref{prop:coarse-alg}, $f^* \calE_X$ is actually a coarse 
structure. If $X$ is connected, then $f^* \calE_X$ is connected. If $X$ is 
unital and $f$ is (globally) proper, then $f^* \calE_X$ is unital. The 
following are clear.

\begin{proposition}
If $X$ is a coarse space and $f \from Y \to X$ is a set map, then $f^* \calE_X$ 
is the maximum coarse structure on $Y$ which makes $f$ into a coarse map.
\end{proposition}

\begin{corollary}\label{cor:crs-factor-I}
If $f \from Y \to X$ is a coarse map, then $f$ factors as a composition of 
coarse maps
\[
    Y \nameto{\smash{\beta}} |Y|_{f^* \calE_X} \nameto{\smash{\utilde{f}}} X,
\]
where $\beta = \id_Y$ and $\utilde{f} = f$ as set maps.
\end{corollary}

More generally, if $\set{X_j \suchthat j \in J}$ ($J$ some index set) is a 
collection of coarse spaces and $\set{f_j \from Y \to X_j}$ is a collection of 
set maps, then
\[
    \calE \defeq \bigintersect_{j \in J} (f_j)^* \calE_{X_j}
\]
is the maximum coarse structure on $Y$ which makes all the $f_j$ into coarse 
maps. If $Y$ is a coarse space and the $f_j \from Y \to X_j$ are all coarse 
maps, then each $f_j$ factors as a composition of coarse maps $\utilde{f}_j 
\circ \beta$ in the obvious way. Moreover, if all the $X_j$ are connected, then 
$\calE$ is connected; if all the $X_j$ are unital and all the $f_j$ are 
(globally) proper, then $\calE$ is unital.

\begin{definition}
Suppose $Y$ is a coarse space. The \emph{push-forward coarse structure} (of 
$\calE_Y$) on $X$ along a \emph{locally proper} map $f \from Y \to X$ is
\[
    f_* \calE_Y \defeq \bigl\langle \set{f^{\cross 2}(F) \suchthat F \in
            \calE_Y} \bigr\rangle.
\]
We similarly define \emph{unital}, \emph{connected}, and \emph{connected unital 
push-forward coarse structures}.
\end{definition}

If $Y$ is connected and $f$ is surjective, then $f_* \calE_Y$ is connected. 
Similarly, if $Y$ is unital (hence $f$ globally proper) and $f$ is surjective, 
then $f_* \calE_Y$ is unital.

\begin{proposition}
If $Y$ is a coarse space and $f \from Y \to X$ is a locally proper map, then 
$f_* \calE_Y$ is the minimum coarse structure on $X$ which makes $f$ into a 
coarse map.
\end{proposition}

\begin{corollary}\label{cor:crs-factor-II}
If $f \from Y \to X$ is a coarse map, then $f$ factors as a composition of 
coarse maps
\[
    Y \nameto{\smash{\tilde{f}}} |X|_{f_* \calE_Y} \nameto{\smash{\alpha}} X
\]
where $\tilde{f} = f$ and $\alpha = \id_X$ as set maps.
\end{corollary}

Of course, there are obvious unital, connected, and connected unital versions 
of the above. For the unital versions one needs $f$ to be proper and $Y$ should 
probably be unital; for the connected versions, $Y$ should probably be 
connected.

More generally, if $\set{Y_j \suchthat j \in J}$ ($J$ some index set) is a 
collection of coarse spaces and $\set{f_j \from Y_j \to X}$ is a collection of 
locally proper maps, then
\[
    \calE \defeq \bigl\langle (f_j)_* \calE_{Y_j} \bigr\rangle
\]
is the minimum coarse structure on $X$ which makes all the $f_j$ into coarse 
maps. If $X$ is a coarse space and the $f_j \from Y_j \to X$ are all coarse 
maps, then each $f_j$ factors as $\alpha \circ \tilde{f}_j$. Again, there are 
unital, connected, and connected unital versions of this.

\begin{remark}
We emphasize that whereas one can pull back coarse structures along \emph{any} 
set map (or collection of set maps), one can only push forward coarse 
structures along \emph{locally proper} maps. If one wants all the coarse 
structures to be unital (and take unital, possibly connected, push-forwards), 
then one evidently requires all maps to be (globally) proper.
\end{remark}

It is easy to see what happens when one pushes a coarse structure forward and 
then pulls it back along the same map (or vice versa).

\begin{proposition}
If $Y$ is a coarse space and $f \from Y \to X$ is a locally proper map, then 
$\calE_Y \subseteq f^* f_* \calE_Y$.
\end{proposition}

\begin{proof}
$f$ is coarse as a map $Y \to |X|_{f_* \calE_Y}$. Applying 
Corollary~\ref{cor:crs-factor-I}, this map factors as $Y \nameto{\smash{\beta}} 
|Y|_{f^* f_* \calE_Y} \to |X|_{f_* \calE_Y}$ where $\beta$ is the identity as a 
set map.
\end{proof}

\begin{proposition}
If $X$ is a coarse space and $f \from Y \to X$ is any set map, then $f_* f^* 
\calE_X \subseteq \calE_X$.
\end{proposition}

\begin{proof}
Now $f$ is coarse as a map $|Y|_{f^* \calE_X} \to X$, to which we apply 
Corollary~\ref{cor:crs-factor-II}.
\end{proof}

Using push-forward coarse structures (and Corollary~\ref{cor:loc-prop-gen}), we 
can ``restate'' Proposition~\ref{prop:loc-prop-for-comp} as follows.

\begin{proposition}\label{prop:loc-prop-comp}
Consider the composition of (set) maps $Z \nameto{\smash{g}} Y 
\nameto{\smash{f}} X$, where $Z$ is a coarse space:
\begin{enumerate}
\item\label{prop:loc-prop-comp:I} If $g$ is locally proper and $f$ is locally 
    proper for the push-forward coarse structure $g_* \calE_Z$ on $Y$, then $f 
    \circ g$ is locally proper.
\item\label{prop:loc-prop-comp:II} If $f \circ g$ is locally proper, then $g$ 
    is locally proper.
\item\label{prop:loc-prop-comp:III} If $f \circ g$ is locally proper, then $f$ 
    is locally proper for the push-forward coarse structure $g_* \calE_Z$ on 
    $Y$.
\end{enumerate}
The above also hold with connected push-forward coarse structures in place of 
push-forward coarse structures. Also, that injectivity implies global 
properness implies local properness.
\end{proposition}

\begin{remark}
Applying the above Proposition with $Z \defeq |Z|_1$ having the terminal coarse 
structure, we get \enumref{prop:prop:I} and \enumref{prop:prop:II} of 
Proposition~\ref{prop:prop}. If $g$ is surjective, then the push-forward coarse 
structure $g_* \calE_{|Z|_1}$ is the terminal coarse structure $\calE_{|Y|_1}$ 
and we get \enumref{prop:prop:III} as well.
\end{remark}

%%%%%%%%%%%%%%%%%%%%%%%%%%%%%%%%%%%%%%%%%%%%%%%%%%%%%%%%%%%%%%%%%%%%%%%%%%%%%%%%

\section{The precoarse categories}\label{sect:PCrs}

We now define several categories of coarse spaces, whose arrows are coarse 
maps, and examine their properties. These \emph{precoarse categories} differ 
from the coarse categories, which are quotients of these categories (see 
\S\ref{sect:Crs}).

\subsection{Set and category theory}\label{subsect:set-cat}

We will be unusually careful with our set and category theoretic constructions. 
The following can mostly be ignored safely, though will be needed eventually 
for rigorous, ``canonical'' constructions (e.g., when we consider sets of 
``all'' modules over a coarse space).

Assuming the Grothendieck axiom that any set is contained in some universe, we 
first fix a universe $\calU$ (containing $\omega$). \emph{Small} (or 
$\calU$-small) objects are elements of $\calU$. A \emph{$\calU$-category} is 
one whose object set is a subset of $\calU$. A ($\calU$-)small category is one 
whose object set (hence morphism set and composition law) is in $\calU$. A 
small category is necessarily a $\calU$-category, but not vice versa. A 
$\calU$-category in turn is $\calU^+$-small, where $\calU^+$ denotes the 
smallest universe having $\calU$ as an element. A \emph{locally small} 
$\calU$-category is a $\calU$-category whose $\Hom$-sets $\Hom(\cdot,\cdot)$ 
are all small.

Recall the notion of quotient categories (from, e.g., \cite{MR1712872}*{Ch.~II 
\S{}8}): Given a category $\calC$ and an equivalence relation $\sim$ on each 
$\Hom$-set of $\calC$, there is a \emph{quotient category} $\calC/{\sim}$ and a 
\emph{quotient functor} $\Quotient \from \calC \to \calC/{\sim}$ satisfying the 
following universal property: For all functors $F \from \calC \to \calC'$ 
($\calC'$ any category, which can be taken to be $\calU$-small if $\calC$ is 
$\calU$-small) such that $f \sim f'$ ($f$, $f'$ in some $\Hom$-set of $\calC$) 
implies $F(f) \sim F(f')$, there is a unique functor $F' \from \calC/{\sim} \to 
\calC'$ such that $F = F' \circ \Quotient$. Moreover, if the equivalence 
relation $\sim$ is preserved under composition then, for all objects $X$, $Y$ 
of $\calC$, the set $\Hom_{\calC/{\sim}}(\Quotient(Y),\Quotient(X))$ is in 
natural bijection with the set of $\sim$-equivalence classes of 
$\Hom_{\calC}(Y,X)$.

As usual, $\CATSet$ denotes the category of small sets (and set maps). 
$\CATTop$ is the category of small topological spaces and continuous maps. 
Forgetful functors will be denoted by $\Forget$, with the source and target 
categories (the latter often being $\CATSet$) implied by context. For a 
category $\calC$ equipped with a forgetful functor to $\CATSet$, we denote the 
full subcategory of $\calC$ of nonempty objects (i.e., those $X$ with 
$\Forget(X) \neq \emptyset$) by $\CATne{\calC}$.

For the most part, henceforth $X$, $Y$, and $Z$ will be (small) coarse spaces, 
and $f \from Y \to X$ and $g \from Z \to Y$ coarse maps. $\setZplus \defeq 
\set{n \in \setZ \suchthat n \geq 0}$ is the set of nonnegative integers and 
similarly $\setRplus \defeq \coitvl{0,\infty}$ is the set of nonnegative real 
numbers.

\subsection{The precoarse categories}

\begin{definition}
The \emph{precoarse category} $\CATPCrs$ has as objects all (small) coarse 
spaces and as arrows coarse maps. The \emph{connected precoarse category} 
$\CATConnPCrs$ is full subcategory of $\CATPCrs$ consisting of the connected 
coarse spaces. Similarly define the \emph{unital precoarse category} 
$\CATUniPCrs$ and the \emph{connected unital precoarse category} 
$\CATConnUniPCrs$.
\end{definition}

\begin{remarks}
In many ways, the category $\CATne{\CATConnPCrs}$ of nonempty connected coarse 
spaces, i.e., coarse spaces with exactly one connected component, is more 
natural. Observe that that $\CATConnUniPCrs = \CATConnPCrs \intersect 
\CATUniPCrs$ is a full subcategory of the other three categories. (One might 
argue that the unital categories above are not the ``correct'' ones and further 
insist that the arrows in the unital categories should be ``unit 
preserving'', i.e., surjective as set maps. However, the above unital 
categories are the usual ones used in coarse geometry; see \S\ref{sect:top-crs} 
and especially Corollary~\ref{cor:UniCrs-RoeCrs-equiv}.)
\end{remarks}

We will analyze various properties of the categories $\CATPCrs$ and 
$\CATConnPCrs$ (which are better behaved than the others). In particular, we 
examine limits and colimits in these categories, which include as special cases 
products and coproducts, equalizers and coequalizers, and terminal and initial 
objects. (We use the standard terminology from category theory, topology, etc.: 
limits are also called ``projective limits'' or ``inverse limits'' and colimits 
are called ``inductive limits'' or ``direct limits'', though ``direct limits'' 
are often more specifically filtered colimits.)

Let us first recall some standard terminology (see, e.g., \cite{MR1712872}). 
Let $\calC$ be a category and suppose $\calF_X \from \calJ \to \calC$ ($\calJ$ 
a small, often finite, category) is a functor. A \emph{cone $\nu \from X \to 
\calF_X$ to $\calF_X$} consists of an $X \in \Obj(\calC)$ and arrows $\nu_j 
\from X \to X_j \defeq \calF_X(j)$, $j \in \Obj(\calJ)$, such that the 
triangles emanating from $X$ commute. A \emph{limit} in $\calC$ for $\calF_X$ 
is given by a cone $X \to \calF_X$ which is universal, i.e., a \emph{limiting 
cone}. Limits of $\calF_X$ in $\calC$ are unique up to natural isomorphism. 
Thus we will sometimes follow the customary abuses of referring to \emph{the} 
limit of $\calF_X$ and of referring to the object $X$ (often denoted $\OBJlim 
\calF_X$) as the limit with the $\nu_j$ understood. A functor $F \from \calC 
\to \calC'$ \emph{preserves limits} if whenever $\nu \from X \to \calF_X$ is a 
limiting cone in $\calC$, $F \circ \nu \from F(X) \to F \circ \calF_X$ is 
limiting in $\calC'$. Dually, one has \emph{cones from $\calF_Y$}, 
\emph{colimits}, \emph{colimiting cones}, and functors which \emph{preserve 
colimits}. All limits and colimits considered will be small. In particular, the 
category $\calJ$ and functors $\calF_X$ and $\calF_Y$ will be small.

First, we examine the relation between $\CATPCrs$ and $\CATConnPCrs$.

\subsection{\pdfalt{\maybeboldmath $\CATPCrs$ versus $\CATConnPCrs$}%
        {PCrs versus CPCrs}}\label{subsect:PCrs-conn}

Below, $I$ will always denote the inclusion $\CATConnPCrs \injto \CATPCrs$. 
Note that $I$ is fully faithful.

\begin{definition}
$\Connect \from \CATPCrs \to \CATConnPCrs$ is the functor defined as follows:
\begin{enumerate}
\item For a coarse space $X$, $\Connect(X)$ is just $X$ as a set, but with the 
    \emph{connected} coarse structure $\langle \calE_X \rangle_X^\TXTconn$ 
    generated by $\calE_X$.
\item For a coarse map $f \from Y \to X$, $\Connect(f) \from \Connect(Y) \to 
    \Connect(X)$ is the same as a set map as $f$ (which is coarse by 
    Corollary~\ref{cor:coarse-gen}).
\end{enumerate}
\end{definition}

The following is clear.

\begin{proposition}
$\Connect \circ I$ is the identity functor on $\CATConnPCrs$.
\end{proposition}

\begin{proposition}%\label{prop:PCrs-Conn-left-adj}
$\Connect \from \CATPCrs \to \CATConnPCrs$ is left adjoint to the inclusion 
functor.
\end{proposition}

The counit maps $Y \to I(\Connect(Y))$, $Y \in \Obj(\CATPCrs)$, of the above 
adjunction are just the identities as set maps. The unit maps $X = 
\Connect(I(X)) \to X$, $X \in \Obj(\CATConnPCrs)$, are the identity maps.

\begin{proof}
Since $\Connect \circ I$ is the identity, $\Connect$ induces natural maps
\[
    \Hom_{\CATPCrs}(Y, I(X)) \to \Hom_{\CATConnPCrs}(\Connect(Y), X)
\]
(for $Y$ possibly disconnected and $X$ connected), which are clearly bijections.
\end{proof}

\begin{corollary}\label{cor:PCrs-ConnPCrs-preserve}
$I \from \CATConnPCrs \injto \CATPCrs$ preserves limits and $\Connect \from 
\CATPCrs \to \CATConnPCrs$ preserves colimits. Moreover, if $\calF \from \calJ 
\to \CATConnPCrs$ is a functor and $\nu$ is a limiting cone to (or colimiting 
cone from) $I \circ \calF$ in $\CATPCrs$, then $\Connect \circ \nu$ is a 
limiting cone to (or colimiting cone from, respectively) $\calF = \Connect 
\circ I \circ \calF$ in $\CATConnPCrs$.
\end{corollary}

\begin{proof}
See, e.g., \cite{MR1712872}*{Ch.~V \S{}5} or \cite{MR0349793}*{16.4.6} for the 
first statement, and \cite{MR0349793}*{16.6.1} for the second.
\end{proof}

\subsection{Limits in the precoarse categories}

\begin{theorem}\label{thm:PCrs-lim}
$\CATPCrs$ has all nonzero limits (i.e., limits of functors $\calJ \to 
\CATPCrs$ for $\calJ$ nonempty). Moreover, the forgetful functor $\Forget \from 
\CATPCrs \to \CATSet$ preserves limits, and the limits of connected coarse 
spaces are connected. Consequently, the same hold with $\CATConnPCrs$ in place 
of $\CATPCrs$.
\end{theorem}

% I'm still mostly/slightly baffled by the remark in Mac Lane, top of p. 35, 
% that Set (and Grp, Top, etc.!) are not locally small.
It is actually easy to see that $\Forget \from \CATPCrs \to \CATSet$ preserves 
limits: $\Forget$ is naturally equivalent to the covariant $\Hom$-functor 
$\Hom_{\CATPCrs}(\ast,\cdot) \from \CATPCrs \to \CATSet$, where $\ast$ is any 
one-point coarse space, and thus preserves limits (see, e.g., 
\cite{MR1712872}*{Ch.~V \S{}4 Thm.~1}). Since I do not know a similar argument 
for colimits, let us proceed in ignorance of this.

\begin{proof}
Recall that $\CATSet$ has all limits. Given $\calF_X \from \calJ \to \CATPCrs$, 
fix a limiting set cone $\nu \from X \to \Forget \circ \calF_X$, so that $X$ is 
a set and $\nu_j \from X \to X_j \defeq \calF_X(j)$, $j \in \Obj(\calJ)$, are 
set maps. It suffices to put a coarse structure on $X$ so that we get a 
limiting cone $\nu \from X \to \calF_X$ in $\CATPCrs$ (with $X$ connected if 
all the $X_j$ are connected).

We need all the $\nu_j \from X \to X_j$ to become coarse maps. Taking the 
coarse structure on $X$ to be the intersection
\[
    \calE_X \defeq \bigintersect_{j \in \Obj(\calJ)} (\nu_j)^* \calE_{X_j}
\]
of pull-back coarse structures clearly makes this so. (Since pull-backs of 
connected coarse structures are connected and intersections of connected coarse 
structures are connected, $\calE_X$ is connected if all the $X_j$ are.) Since 
$\Forget$ is faithful, $\nu \from X \to \calF_X$ is a cone in $\CATPCrs$. We 
must show that it is universal.

Suppose $\mu \from Y \to \calF_X$ is another cone in $\CATPCrs$. Applying 
$\Forget$, we get a cone $\mu \from Y \to \Forget \circ \calF_X$ in $\CATSet$ 
(properly written $\Forget \circ \mu \from \Forget(Y) \to \Forget \circ 
\calF_X$). Since $\nu$ is universal in $\CATSet$, there is a set map $t \from Y 
\to X$ such that $\mu = \nu \circ t$ as cones in $\CATSet$. We must show that 
$t$ is actually a coarse map (uniqueness is clear).

First, since $\calJ$ is nonzero, there is some object $j_0 \in \Obj(\calJ)$; 
then $\mu_{j_0} = \nu_{j_0} \circ t$ (as set maps) is locally proper, so $t$ is 
locally proper 
(Proposition~\ref{prop:loc-prop-comp}\enumref{prop:loc-prop-comp:II}). Next, 
for each $j \in \Obj(\calJ)$ and $F \in \calE_Y$, $\nu_j$ is coarse for $E 
\defeq t^{\cross 2}(F)$ (which is in $\calE_{|X|_1}$, by local properness) and 
hence $E \in (\nu_j)^* \calE_{X_j}$: Since $\mu_j = \nu_j \circ t$ is locally 
proper for $F$, $\nu_j$ is locally proper for $E$ 
(Proposition~\ref{prop:loc-prop-for-comp}\enumref{prop:loc-prop-for-comp-II}), 
and also $\nu_j$ clearly preserves $E$.

For $\CATConnPCrs$, the assertions follow from 
Corollary~\ref{cor:PCrs-ConnPCrs-preserve}.
\end{proof}

The above proof gives a rather concrete description of limits in $\CATPCrs$ 
(and $\CATConnPCrs$), and in particular of products. The product 
$\pfx{\CATPCrs}\prod_{j \in J} X_j$ in $\CATPCrs$ (or in $\CATConnPCrs$) is 
just the set product (i.e., cartesian product) $X \defeq \pfx{\CATSet}\prod_{j 
\in J} X_j$ together with the entourages of $|X|_1$ which project properly to 
entourages of all the $X_j$.

The ``nonzero'' stipulation in Theorem~\ref{thm:PCrs-lim} is necessary.

\begin{proposition}\label{prop:PCrs-no-map-to}
For each coarse space $X$, there exists a (nonempty) connected, unital coarse 
space $Y$ such that there is no coarse map $Y \to X$.
\end{proposition}

\begin{proof}
Given $X$, take $Y \defeq |Y|_1$ to be an infinite set with cardinality 
strictly greater than the cardinality of $X$, equipped with the terminal coarse 
structure (which is connected and unital), e.g., $Y \defeq |\powerset(X) 
\disjtunion \setN|_1$. Then no locally proper map $Y \to X$ exists, since no 
globally proper map $Y \to X$ exists and $Y$ is unital (see 
Corollary~\ref{cor:loc-prop-uni}). (Note that the cardinality of sets in our 
universe $\calU$ is bounded above by some cardinal, namely by $\# \calU$, but 
no element of $\calU$ has this cardinality.)
\end{proof}

\begin{corollary}\label{cor:PCrs-no-term}
None of the precoarse categories ($\CATPCrs$, $\CATConnPCrs$, 
$\CATne{\CATConnPCrs}$, $\CATUniPCrs$, and $\CATConnUniPCrs$) has a terminal 
object.
\end{corollary}

The failure of existence of terminal objects in the precoarse categories is not 
just a failure of uniqueness of maps, but more seriously of existence. Thus we 
will also get the following on the coarse categories (which are quotients of 
the precoarse categories).

\begin{corollary}\label{cor:Crs-no-term}
No quotient of any of the above precoarse categories has a terminal object.
\end{corollary}

It is straightforward to show that the inclusion $\CATne{\CATConnPCrs} \injto 
\CATConnPCrs$ preserves limits, and moreover that a nonzero limit exists in 
$\CATne{\CATConnPCrs}$ if and only if the corresponding set limit is nonempty 
(but the example below shows that $\CATne{\CATConnPCrs}$ does \emph{not} have 
all nonzero limits). On the other hand unitality poses a fatal problem: The 
forgetful functor $\CATUniPCrs \to \CATSet$ still preserves limits, so a 
(nonzero) limit in $\CATUniPCrs$ can only exist when all the maps in the 
corresponding limiting set cone are proper (but this is often not the case, 
e.g., in the case of products).

\begin{example}
Let $X \defeq |\setZplus|_1$ (which is connected and nonempty), $f \defeq \id_X 
\from X \to X$ be the identity, and define $g \from X \to X$ by $g(x) \defeq 
x+1$. Then the equalizer of $f$ and $g$ in $\CATPCrs$ is the empty set.
\end{example}

To get ahead of ourselves (see \S\ref{sect:Crs}), note that though $f$ and $g$ 
are \emph{close}, the equalizer of $f$ and itself (which is just $X$ mapping 
identically to itself) is not \emph{coarsely equivalent} to the equalizer of 
$f$ and $g$. Indeed, one can obtain other inequivalent equalizers: e.g., the 
equalizer of $h \from X \to X$ where $h(x) \defeq \min \set{0,x-1}$ (which also 
close to $f$) and $f$ is $\set{0}$ (including into $X$). On the other hand, in 
the quotient coarse category $\CATCrs$, $[f] = [g] = [h]$ so the equalizer of 
any pair of these maps is $X$. Since limits in $\CATPCrs$ are not 
\emph{coarsely invariant}, they are of limited interest.

We also see that the quotient functor $\CATPCrs \to \CATCrs$ does not preserve 
equalizers, hence does not preserve limits. However, it \emph{does} preserve 
products. Using this and a method parallel to the one employed in the proof of 
Theorem~\ref{thm:PCrs-lim}, we will show that $\CATCrs$ also has all nonzero 
limits (which will, by definition, be coarsely invariant).

We will use products extensively. We take this opportunity to mention several 
canonical coarse maps which arise due to the existence of (nonzero) products 
(all objects are coarse spaces and arrows coarse maps):
\begin{enumerate}
\item For any $X$ and $Y$, there are \emph{projection maps} $\pi_X \from X 
    \cross Y \to X$ and $\pi_Y \from X \cross Y \to Y$.
\item For any $X$, there is a \emph{diagonal map} $\Delta_X \from X \to X 
    \cross X$.
\item For $f \from Y \to X$ and $f' \from Y' \to X'$, there is a \emph{product 
    map} $f \cross f' \from Y \cross Y' \to X \cross X'$.
\end{enumerate}
The above can all be generalized to larger (even infinite) products.

\subsection{Colimits in the precoarse categories}

\begin{theorem}\label{thm:PCrs-colim}
A colimit exists in $\CATPCrs$ if and only if all the maps from a corresponding 
colimiting set cone are locally proper. Moreover, the forgetful functor 
$\Forget \from \CATPCrs \to \CATSet$ preserves colimits. The same hold with 
$\CATConnPCrs$ in place of $\CATPCrs$.
\end{theorem}

\begin{proof}
This proof is basically dual to the proof of Theorem~\ref{thm:PCrs-lim}, only 
with the added onus of showing the ``only if''. The reason for the local 
properness requirement is that coarse structures can only be pushed forward 
along locally proper maps (whereas they can be pulled back along all maps).

Recall that $\CATSet$ has all colimits. Given $\calF_Y \from \calJ \to 
\CATPCrs$, fix a colimiting set cone $\nu \from \Forget \circ \calF_Y \to Y$, 
so that $Y$ is a set and $\nu_j \from Y_j \defeq \calF_Y(j) \to Y$, $j \in 
\Obj(\calJ)$, are set maps. Suppose all of the $\nu_j$ are locally proper. 
Taking the coarse structure on $Y$ to be
\[
    \calE_Y \defeq \langle (\nu_j)_* \calE_{Y_j}
            \suchthat j \in \Obj(\calJ) \rangle_Y,
\]
we clearly get a cone $\nu \from \calF_Y \to Y$ in $\CATPCrs$; we must prove 
that it is universal.

Suppose $\mu \from \calF_Y \to X$ is another cone in $\CATPCrs$. Then there is 
a canonical set map $t \from Y \to X$ such that $\mu = t \circ \nu$ as cones in 
$\CATSet$. We must show that $t$ is coarse (again uniqueness is clear). 
Entourages $(\nu_j)^{\cross 2}(F)$, $F \in \calE_{Y_j}$, $j \in \Obj(\calJ)$, 
generate $\calE_Y$. $t$ is locally proper for each such entourage (using $\mu_j 
= t \circ \nu_j$ and 
Proposition~\ref{prop:loc-prop-for-comp}\enumref{prop:loc-prop-for-comp-III}) 
and clearly preserves each such entourage. Thus $t$ is coarse, as required.

If the $\nu_j$ are \emph{not} all locally proper, we must show that $\calF_Y$ 
does not have a colimit (in $\CATPCrs$); in fact, we show something stronger, 
that there is no cone from $\calF_Y$ in $\CATPCrs$. We proceed by 
contradiction, so suppose that $\nu_{j_0}$ is not locally proper ($j_0 \in 
\Obj(\calJ)$ fixed) and suppose $\mu \from \calF_Y \to X$ is a cone in 
$\CATPCrs$. Again there must be a set map $t \from Y \to X$ such that $\mu = t 
\circ \nu$ as set cones. But then $\mu_{j_0} = t \circ \nu_{j_0}$ is locally 
proper, which implies that $\nu_{j_0}$ is locally proper 
(Proposition~\ref{prop:loc-prop-comp}\enumref{prop:loc-prop-comp:II}) which is 
a contradiction.

To get the asserted colimits in $\CATConnPCrs$, simply apply 
Corollary~\ref{cor:PCrs-ConnPCrs-preserve}. To show that $\CATConnPCrs$ has no 
more colimits than $\CATPCrs$ (i.e., $\calF_Y \from \calJ \to \CATConnPCrs$ has 
a colimit in $\CATConnPCrs$ only if $I \circ \calF_Y \from \calJ \to \CATPCrs$ 
has a colimit in $\CATPCrs$), it is probably simplest to modify the above 
proof.
\end{proof}

The following are clear.

\begin{corollary}
$\CATPCrs$ and $\CATConnPCrs$ have all coproducts.
\end{corollary}

\begin{corollary}
The empty coarse space is the (unique) initial object in $\CATPCrs$ and in 
$\CATConnPCrs$.
\end{corollary}

\begin{corollary}\label{cor:PCrs-fin-components}
Any coarse space with only \emph{finitely} many connected components is 
(isomorphic in $\CATPCrs$ to) the coproduct in $\CATPCrs$ of its connected 
components.
\end{corollary}

The above Corollary does not necessarily hold for coarse spaces with infinitely 
many connected components. One may say, more generally, that any coarse space 
whose unital subspaces have only finitely many connected components is the 
coproduct of its connected components.

We get concrete descriptions of coproducts in $\CATPCrs$ and in $\CATConnPCrs$. 
The coproduct $\pfx{\CATPCrs}\coprod_{j \in J} Y_j$ in $\CATPCrs$ is just the 
set coproduct (i.e., disjoint union) $Y \defeq \pfx{\CATSet}\coprod_{j \in J} 
Y_j$ with entourages finite unions of entourages of the $Y_j$ (included into 
$Y$). The corresponding coproduct in $\CATConnPCrs$ is the same as a set, but 
one may also take an additional union with an arbitrary finite subset of 
$Y^{\cross 2}$.

The inclusion $\CATne{\CATConnPCrs} \injto \CATConnPCrs$ preserves colimits. 
$\CATne{\CATConnPCrs}$ does not have a zero colimit (i.e., initial object), but 
otherwise has a colimit if the corresponding colimit exists in $\CATConnPCrs$, 
in which case the two colimits coincide; note that a nonzero colimit of 
nonempty sets is nonempty. Unitality does not pose a problem for colimits: 
Theorem~\ref{thm:PCrs-colim} also holds with $\CATUniPCrs$ in place of 
$\CATPCrs$ (and $\CATConnUniPCrs$ in place of $\CATConnPCrs$). In the proof, 
one simply takes the unital coarse structure
\[
    \langle (\nu_j)_* \calE_{Y_j} \suchthat j \in \Obj(\calJ) \rangle_Y^\TXTuni
\]
instead. Of course, in the unital cases, one may substitute ``(globally) 
proper'' for ``locally proper''.

The ``locally proper'' hypothesis is necessary, as the following shows.

\begin{example}
Let $X \defeq |\setZplus|_1$, $f \from X \to X$ be the identity, and define $g 
\from X \to X$ by $g(x) \defeq \min \set{0,x-1}$. The coequalizer of $f$ and 
$g$ in $\CATSet$ is the one-point set $\ast$; since $X$ is unital, $f$ and $g$ 
do not have a coequalizer in $\CATPCrs$.
\end{example}

Again, to get ahead of ourselves, we see that coequalizers in $\CATPCrs$ are 
not coarsely invariant. Even though $f$ is close to $g$ and the coequalizer of 
$f$ and itself is just $X$, $f$ and $g$ do not have a coequalizer in 
$\CATPCrs$. In the quotient category $\CATCrs$, there are no problems: the 
coequalizer of $[f]$ and $[g]$ is $X$, as expected.

The quotient functor $\CATPCrs \to \CATCrs$ does not preserve coequalizers or 
colimits in general. However, it \emph{does} preserve coproducts, and we will 
use these to show that in fact $\CATCrs$ has \emph{all} colimits (which are 
evidently coarsely invariant). In particular, $\CATCrs$ has all coequalizers, 
which contrasts with the situation in $\CATPCrs$ (recall that having all 
coproducts and all coequalizers would imply having all colimits).

%%%%%%%%%%%%%%%%%%%%%%%%%%%%%%%%%%%%%%%%%%%%%%%%%%%%%%%%%%%%%%%%%%%%%%%%%%%%%%%%

\section{The coarse categories}\label{sect:Crs}

\subsection{Closeness of maps}

In classical (unital) coarse geometry, two maps $f, f' \from Y \to X$ are 
\emph{close} if $(f \cross f')(1_Y)$ is an entourage of $X$. Closeness is an 
equivalence relation on maps $Y \to X$, but note that it does not involve the 
coarse structure on $Y$ at all! In the nonunital case, we must modify the 
definition, lest closeness not even be reflexive (e.g., take $Y \defeq X$ 
nonunital and $f \defeq f' \defeq \id_X$).

\begin{definition}
Coarse maps $f, f' \from Y \to X$ are \emph{close} (write $f \closeequiv f'$) 
if $(f \cross f')(F) \in \calE_X$ for all $F \in \calE_Y$.
\end{definition}

\begin{proposition}
Closeness of coarse maps $Y \to X$ is an equivalence relation (on the 
$\Hom$-set $\Hom_{\CATPCrs}(Y,X)$).
\end{proposition}

\begin{proof}
Reflexivity follows since coarse maps preserve entourages. Symmetry follows by 
taking transposes. Transitivity: Suppose $f, f', f'' \from Y \to X$ are coarse 
maps with $f \closeequiv f'$ and $f' \closeequiv f''$. For any $F \in \calE_Y$,
\[
    (f \cross f'')(F)
        \subseteq (f \cross f')(1_{F \cdot Y}) \circ (f' \cross f'')(F)
\]
is an entourage of $X$ since $1_{F \cdot Y} \in \calE_Y$ (since $1_{F \cdot Y} 
\subseteq F \circ F^\transpose$).
\end{proof}

Like local properness, closeness is also determined ``on'' unital subspaces of 
the domain. Thus for unital coarse spaces, our notion of closeness is just the 
classical one.

\begin{proposition}\label{prop:close-uni}
Coarse maps $f, f' \from Y \to X$ are close if and only if, for every unital 
subspace $Y' \subseteq Y$, $f |_{Y'}$ and $f' |_{Y'}$ are close (i.e., $(f 
\cross f')(1_{Y'}) \in \calE_X)$. Thus, for $Y$ unital, $f$ and $f'$ are close 
if and only if $(f \cross f')(1_Y) \in \calE_X$.
\end{proposition}

\begin{proof}
(\textimplies): Immediate.

(\textimpliedby): For $F \in \calE_Y$, $Y' \defeq F \cdot Y \union Y \cdot F$ 
is a unital subspace of $Y$, and $F \in \calE_{Y'}$. Then
\[
    (f \cross f')(F) = (f |_{Y'} \cross f' |_{Y'})(F) \in \calE_X,
\]
as required.
\end{proof}

We have not used local properness at all, so we can actually define closeness 
for maps which preserve entourages (but are not necessarily locally proper). 
However, we will not need this.

The following observation is rather important.

\begin{proposition}\label{prop:term-close}
Suppose $f, f' \from Y \to X$ are coarse maps. If $X = |X|_1$ has the terminal 
coarse structure, then $f$ and $f'$ are close.
\end{proposition}

Thus if $X = |X|_1$, then for any coarse space $Y$ there is \emph{at most} one 
(but possibly no) closeness class of coarse map $Y \to X$.

\begin{proof}
This follows immediately from Proposition~\ref{prop:loc-prop-prod}.
\end{proof}

\subsection{The coarse categories}

Closeness, an equivalence relation on the $\Hom$-sets of $\CATPCrs$, yields a 
quotient category
\[
    \CATCrs \defeq \CATPCrs/{\closeequiv}
\]
(see \S\ref{subsect:set-cat}), which we call the \emph{coarse category}, 
together with a quotient functor $\Quotient \from \CATPCrs \to \CATCrs$. We may 
similarly define quotients $\CATConnCrs$, $\CATUniCrs$, and $\CATConnUniCrs$ of 
$\CATConnPCrs$, $\CATUniPCrs$, and $\CATConnUniPCrs$, respectively. These 
latter categories are full subcategories of $\CATPCrs$, so their quotients are 
full subcategories of $\CATCrs$.

The following is clear.

\begin{proposition}
Closeness is respected by composition: If $f, f' \from Y \to X$ and $g, g' 
\from Z \to Y$ are coarse maps with $f \closeequiv f'$ and $g \closeequiv g'$, 
then $f \circ g \closeequiv f' \circ g'$.
\end{proposition}

This allows us to describe the arrows of $\CATCrs$ as closeness equivalence 
classes of coarse maps. Denote such classes by $[f]_\TXTclose \from Y \to X$ 
(or simply $[f]$ for brevity), where $f$ is usually taken to be a 
representative map $Y \to X$, i.e., $\Quotient(f) = [f]$. However, we will use 
the notation $[f] \from Y \to X$ for arrows $Y \to X$ in $\CATCrs$ even when we 
do not have a particular $f$ in mind.

The notion of isomorphism in $\CATCrs$ is weaker than in $\CATPCrs$. A coarse 
map $f \from Y \to X$ is a \emph{coarse equivalence} if $[f]$ is an isomorphism 
in $\CATCrs$. In other words, $f$ is a coarse equivalence if and only if there 
is a coarse map $g \from X \to Y$ so that the two possible compositions are 
close to the identities (i.e., $[f \circ g] = [\id_X]$ and $[g \circ f] = 
[\id_Y]$).

A functor $F \from \CATPCrs \to \calC$, $\calC$ any category, is \emph{coarsely 
invariant} if $f \closeequiv f'$ implies $F(f) = F(f')$. Any coarsely invariant 
$F$ induces a functor $[F] \from \CATCrs \to \calC$ with $F = [F] \circ 
\Quotient$. Coarsely invariant functors send coarse equivalences to 
isomorphisms. For functors $F \from \CATPCrs \to \CATPCrs$ (or with codomain 
one of the other precoarse categories), we abuse terminology and also say that 
$F$ is \emph{coarsely invariant} if $\Quotient \circ F \from \CATPCrs \to 
\CATCrs$ is coarsely invariant in the previous (stronger) sense. Such a 
coarsely invariant functor $F \from \CATPCrs \to \CATPCrs$ induces a functor 
$[F] \from \CATCrs \to \CATCrs$; if $F \from \CATPCrs \to \CATConnPCrs$, then 
$[F] \from \CATCrs \to \CATConnCrs$; etc.

\subsection{\pdfalt{\maybeboldmath $\CATCrs$ versus $\CATConnCrs$}%
        {Crs versus CCrs}}

The relation between the quotient categories $\CATCrs$ and $\CATConnCrs$ is 
essentially the same as that between $\CATPCrs$ and $\CATConnPCrs$ for the 
following reasons, which are easy to check.

\begin{proposition}
The functors $I \from \CATConnPCrs \injto \CATPCrs$ and $\Connect \from 
\CATPCrs \to \CATConnPCrs$ are coarsely invariant, hence induce functors $[I] 
\from \CATConnCrs \to \CATCrs$ and $[\Connect] \from \CATCrs \to \CATConnCrs$, 
respectively. In fact, $[I]$ is just the inclusion and is fully faithful. 
Again, $[\Connect] \circ [I]$ is the identity functor (now on $\CATConnCrs$), 
and $[\Connect]$ is left adjoint to $[I]$.
\end{proposition}

Consequently, we get the following (exact) analogues of 
Corollary~\ref{cor:PCrs-ConnPCrs-preserve}.

\begin{corollary}%\label{cor:Crs-ConnCrs-preserve}
$[I]$ preserves limits and $[\Connect]$ preserves colimits. If $\nu$ is a 
limiting cone to (or colimiting cone from) $[I] \circ \calF$, where $\calF 
\from \calJ \to \CATConnCrs$, then $[\Connect] \circ \nu$ is a limiting cone to 
(or colimiting cone from, respectively) $\calF = [\Connect] \circ [I] \circ 
\calF$.
\end{corollary}

\begin{remark}
Evidently, $I$ and $\Connect$ ``commute'' with the quotient functors 
$\Quotient$ ($\CATPCrs \to \CATCrs$ and its restriction $\CATConnPCrs \to 
\CATConnCrs$) in that
\[
    \Quotient \circ I = [I] \circ \Quotient
\quad\text{and}\quad
    \Quotient \circ \Connect = [\Connect] \circ \Quotient.
\]
The quotient functors give a map of adjunctions (see, e.g., 
\cite{MR1712872}*{Ch.~IV \S{}7}) from $(\Connect,I)$ to $([\Connect],[I])$.
\end{remark}

\begin{remark}
$[I]$ is fully faithful, but $[\Connect]$ is neither full nor faithful (even 
though $\Connect$ is faithful, though also not full): e.g., consider
\[
    \Hom_{\CATCrs}(\ast, \ast \copro \ast)
\quad\text{and}\quad
    \Hom_{\CATCrs}(\Connect(|\setZplus|_1 \copro |\setZplus|_1),
            |\setZplus|_1 \copro |\setZplus|_1),
\]
respectively.
\end{remark}

\subsection{\pdfalt{\maybeboldmath $\CATConnCrs$ versus $\CATne{\CATConnCrs}$}%
        {CCrs versus CCrs\textcaret{}x}}

After passing to the quotients by closeness, the situation with respect to 
\emph{nonempty} connected coarse spaces is greatly improved. Below, we work in 
$\CATConnPCrs$ or its quotient $\CATConnCrs$ (or the nonempty subcategories), 
so all coarse spaces will be connected. Let $I \from \CATne{\CATConnPCrs} 
\injto \CATConnPCrs$ denote the inclusion; it is coarsely invariant, hence 
induces $[I] \from \CATne{\CATConnCrs} \injto \CATConnCrs$, which is also the 
inclusion and which is fully faithful. Again, the inclusion functors 
``commute'' with the quotient functors.

\begin{definition}
Fix a one-point coarse space $\ast$. Define a functor
\[
    \AddPt \from \CATConnPCrs \to \CATne{\CATConnPCrs}
\]
as follows:
\begin{enumerate}
\item For a coarse space $X$, $\AddPt(X) \defeq X \copro_{\CATConnPCrs} \ast$ 
    (coproduct in $\CATConnPCrs$).
\item For a coarse map $f \from Y \to X$, $\AddPt(f) \defeq f 
    \copro_{\CATConnPCrs} \id_\ast$.
\end{enumerate}
\end{definition}

(It is easy to construct the functor $\AddPt$ concretely, and all functors 
satisfying the above are naturally equivalent.) The following are all easy to 
verify.

\begin{proposition}
$\AddPt \from \CATConnPCrs \to \CATne{\CATConnPCrs}$ is coarsely invariant and 
hence induces a functor $[\AddPt] \from \CATConnCrs \to \CATne{\CATConnCrs}$.
\end{proposition}

$\AddPt$ is not terribly useful, but $[\AddPt]$ is.

\begin{proposition}%\label{prop:ConnCrs-AddPt-left-adj}
$[\AddPt] \circ [I]$ is naturally equivalent to the identity on 
$\CATne{\CATConnCrs}$. Moreover, $[\AddPt] \from \CATConnCrs \to 
\CATne{\CATConnCrs}$ is left adjoint to $[I]$.
\end{proposition}

It follows that $[\AddPt]$ is naturally equivalent to a functor $[\AddPt]'$ 
such that $[\AddPt]' \circ [I]$ is equal to the identity functor. It is easy to 
give a natural equivalence $\Id_{\CATne{\CATConnCrs}} \to [\AddPt] \circ [I]$: 
for each (nonempty, connected) $X$, the canonical inclusion $\iota_X \from X 
\to X \copro \ast$ is a coarse equivalence hence an isomorphism $[\iota_X] 
\from X \to \AddPt(X)$ in $\CATne{\CATConnCrs}$.

\begin{corollary}%\label{cor:ConnCrs-neConnCrs-preserve}
$[I]$ preserves limits and $[\AddPt]$ preserves colimits. If $\nu$ is a 
limiting cone to (or, a colimiting cone from) $[I] \circ \calF$, where $\calF 
\from \calJ \to \CATne{\CATConnCrs}$, then $[\AddPt] \circ \nu$ is a limiting 
cone to (or, respectively, a colimiting cone from) $[\AddPt] \circ [I] \circ 
\calF$ (or $\calF$ after applying a natural equivalence).
\end{corollary}

\subsection{Limits in the coarse categories}

We first prove our assertion that nonzero products in the nonunital coarse 
categories are just images (under the quotient functor) of products in the 
precoarse categories. We then show the nonunital coarse categories also have 
all equalizers of pairs of arrows. It then follows by a standard construction 
that the nonunital coarse categories have all nonzero limits.

\begin{proposition}\label{prop:Crs-prod}
Suppose $\set{X_j \suchthat j \in J}$ ($J$ some index set) is a nonzero 
collection of coarse spaces. The product of the $X_j$ in $\CATCrs$ (or in 
$\CATConnCrs$ or $\CATne{\CATConnCrs}$, as appropriate) is just the coarse 
space
\[
    X \defeq \pfx{\CATPCrs}\prod_{j \in J} X_j
\]
(product in $\CATPCrs$) together with the projections $[\pi_j] \from X \to 
X_j$, $j \in J$ (closeness classes of the projections). Thus $\CATCrs$ (and 
$\CATConnCrs$ and $\CATne{\CATConnCrs}$) have all nonzero products.
\end{proposition}

Recall, from Corollary~\ref{cor:Crs-no-term}, that none of the quotient coarse 
categories has a zero product, i.e., terminal object.

\begin{proof}
The cone $\pi$ in $\CATPCrs$ maps (via the quotient functor) to a cone $[\pi] 
\defeq \Quotient \circ \pi$ in $\CATCrs$; we must prove universality. Suppose 
$Y$ is a coarse space and $[\mu_j] \from Y \to X_j$, $j \in J$, is a collection 
of arrows in $\CATCrs$. Choosing representative coarse maps $\mu_j \from Y \to 
X_j$, we get (since the cone $\pi$ is universal) a natural coarse map $t \from 
Y \to X$ such that $\mu_j = \pi_j \circ t$ for all $j$. Of course, this implies 
$[\mu_j] = [\pi_j] \circ [t]$.

We must show that this $[t]$ is unique (hence does not depend on our choice of 
representatives $\mu_j$). Suppose $[t'] \from Y \to X$ is a class such that 
$[\mu_j] = [\pi_j] \circ [t']$ for all $j$. Choose a representative $t'$. 
Suppose $F \in \calE_Y$, and put $E \defeq (t \cross t')(F)$, which is in 
$\calE_{|X|_1}$ by Proposition~\ref{prop:loc-prop-prod}. For each $j$, we have 
that $\mu'_j \defeq \pi_j \circ t' \closeequiv \mu_j = \pi_j \circ t$, and 
hence
\[
    (\pi_j)^{\cross 2}(E) = ((\pi_j \circ t) \cross (\pi_j \circ t'))(F)
\]
is in $\calE_{X_j}$. Moreover, since $(\mu_j \cross \mu'_j) |_F = 
(\pi_j)^{\cross 2} |_E \circ (t \cross t') |_E^F$ is proper (by the same 
Proposition), $\pi_j$ is locally proper for $E$. Thus $E \in (\pi_j)^* 
\calE_{X_j}$ for all $j$, so $E \in \calE_X$. Hence $t$ is close to $t'$.

For $\CATConnCrs$ and $\CATne{\CATConnCrs}$, it suffices to recall that nonzero 
products of connected coarse spaces are connected, and nonzero products of 
nonempty coarse spaces are nonempty.
\end{proof}

\begin{remark}
For obvious reasons, we cannot usually obtain products in the unital coarse 
categories using the above construction. However, unlike in $\CATPCrs$, this 
does not imply the nonexistence of products. In certain cases (see, e.g., 
Remark~\ref{rmk:term-unital-prod}), the (nonunital) product above will be 
coarsely equivalent to a unital coarse space which is a product in 
$\CATUniCrs$. I do not know, in general, which products exist in $\CATUniCrs$.
\end{remark}

Next, we examine equalizers in the (nonunital) coarse categories. Unlike 
products, equalizers in the coarse categories are not usually ``the same'' as 
equalizers in the precoarse categories.

\begin{definition}
Suppose $X$ is a coarse space and $f, f' \from Y \to X$ are set maps ($Y$ some 
set). $f$ and $f'$ are \emph{pointwise connected} if $f(y)$ is connected to 
$f'(y)$ for all $y \in Y$. $f$ and $f'$ are \emph{close for $F \in f^* \calE_X 
\intersect (f')^* \calE_X$} if $(f \cross f')(F) \in \calE_X$.
\end{definition}

Of course, if $X$ is connected, all set maps into $X$ are pointwise connected. 
If $f \from Y \to X$ is a coarse map, then any coarse map close to $f$ is 
pointwise connected to $f$.

\begin{lemma}
Suppose $X$ is a coarse space and $f, f' \from Y \to X$ are set maps. If $f$ 
and $f'$ are close for both $F, F' \in \calE_{|Y|_1}$, then $f$ and $f'$ are 
close for $F + F'$, $F \circ F'$, $F^\transpose$, and all subsets of $F$.
\end{lemma}

\begin{proof}
Again, the only (slightly) nontrivial one is $F \circ F'$:
\[
    (f \cross f')(F \circ F') \subseteq f^{\cross 2}(F) \circ (f \cross f')(F')
            \in \calE_X.
\]
\end{proof}

\begin{definition}
Suppose $X$ is a coarse space and $f, f' \from Y \to X$ are pointwise connected 
coarse maps. The \emph{equalizing pull-back coarse structure} $(f,f')^* 
\calE_X$ (on $Y$ along $f$ and $f'$) is
\[
    (f,f')^* \calE_X \defeq \set{F \in f^* \calE_X \intersect (f')^* \calE_X
                \suchthat \text{$f$ and $f'$ are close for $F$}}.
\]
\end{definition}

Pointwise-connectedness is important: it guarantees that the singletons 
$\set{1_y}$, $y \in Y$, are in $(f,f')^* \calE_X$. It then follows from the 
Lemma that $(f,f')^* \calE_X$ is a coarse structure on $Y$.

\begin{definition}
Suppose $f, f' \from Y \to X$ are coarse maps. Define the \emph{equalizer} of 
$[f]$ and $[f']$ is
\[
    \OBJequalizer_{[f],[f']}
        \defeq \set{y \in Y \suchthat \text{$f(y)$ is connected to $f'(y)$}}
        \subseteq Y
\]
with the coarse structure
\[
    \calE_{\OBJequalizer_{[f],[f']}} \defeq \calE_Y |_{\OBJequalizer_{[f],[f']}}
            \intersect (f |_{\OBJequalizer_{[f],[f']}},
                f' |_{\OBJequalizer_{[f],[f']}})^* \calE_X
\]
(where $\calE_Y |_{\OBJequalizer_{[f],[f']}}$ is the subspace coarse structure 
on $\OBJequalizer_{[f],[f']} \subseteq Y$), together with closeness class of 
the inclusion map
\[
    \equalizer_{[f,[f']} \from \OBJequalizer_{[f],[f']} \to Y
\]
(which is coarse).
\end{definition}

Clearly, the restrictions $f |_{\OBJequalizer_{[f],[f']}}$ and $f' 
|_{\OBJequalizer_{[f],[f']}}$ are pointwise connected, so 
$\calE_{\OBJequalizer_{[f],[f']}}$ really is a coarse structure. Also, the 
above definition does not depend on order, i.e., $\OBJequalizer_{[f],[f']} = 
\OBJequalizer_{[f'],[f]}$.

\begin{lemma}
Suppose $f, f' \from Y \to X$ are coarse maps. The equalizer of $[f]$ and 
$[f']$ is coarsely invariant in the sense that $\OBJequalizer_{[f],[f']}$ and 
$[\equalizer_{[f],[f']}]$ (indeed, $\equalizer_{[f],[f']}$) only depend on the 
closeness class of $f$ and $f'$ (hence the notation).
\end{lemma}

\begin{proof}
Suppose $e, e' \from Y \to X$ are close to $f$, $f'$, respectively. Then, for 
all $y \in Y$, $e(y)$ is connected to $f(y)$ and $e'(y)$ is connected to 
$f'(y)$; for $y \in \OBJequalizer_{[f],[f']} \subseteq Y$, $f(y)$ is connected 
to $f'(y)$ hence $e(y)$ is connected to $e'(y)$. Thus the \emph{set} 
$\OBJequalizer_{[f],[f']}$ is coarsely invariant.

It remains to show that the coarse structure $\calE_{\OBJequalizer_{[f],[f']}}$ 
is also coarsely invariant. Observe that
\[
    \calE_{\OBJequalizer_{[f],[f']}}
        = \set{F \in \calE_Y |_{\OBJequalizer_{[f],[f']}} \suchthat
                (f \cross f')(F) \in \calE_X}
\]
and $\calE_Y |_{\OBJequalizer_{[f],[f']}} = \calE_Y 
|_{\OBJequalizer_{[f],[f']}}$. If $F \in \calE_{\OBJequalizer_{[f],[f']}}$, 
then
\[
    (e \cross e')(F) \subseteq (e \cross f)(1_{F \cdot Y})
        \circ (f \cross f')(F) \circ (f' \cross e)(1_{Y \cdot F})
\]
is in $\calE_X$, and so $F \in \calE_{\OBJequalizer_{[e],[e']}}$; the reverse 
inclusion follows symmetrically.
\end{proof}

\begin{proposition}\label{prop:Crs-equal}
The equalizer of $[f], [f'] \from Y \to X$ really is (in the categorical sense) 
the equalizer of $[f]$ and $[f']$ in $\CATCrs$ (or in $\CATConnCrs$ or 
$\CATne{\CATConnCrs}$, as appropriate), hence the terminology. Thus $\CATCrs$ 
(and $\CATConnCrs$ and $\CATne{\CATConnCrs}$) have all equalizers of pairs of 
arrows.
\end{proposition}

\begin{proof}
Fix representative coarse maps $f$ and $f'$, and suppose $g \from Z \to Y$ is a 
coarse map such that $f \circ g \closeequiv f' \circ g$. Then clearly the (set) 
image of $g$ is contained in $\OBJequalizer_{[f],[f']}$, and indeed
\[
    \tilde{g} \defeq g |^{\OBJequalizer_{[f],[f']}}
            \from Z \to \OBJequalizer_{[f],[f']}
\]
is clearly coarse with $g = \equalizer_{[f],[f']} \circ \tilde{g}$ (hence $[g] 
= [\equalizer_{[f],[f']}] \circ [\tilde{g}]$).

We must prove uniqueness of $[\tilde{g}]$. Suppose $\tilde{g}' \from Z \to 
\OBJequalizer_{[f],[f']}$ is a coarse map with $g \closeequiv 
\equalizer_{[f],[f']} \circ \tilde{g}' \eqdef g'$. Then, for all $G \in 
\calE_Z$,
\[
    F \defeq (\tilde{g} \cross \tilde{g}')(G) = (g \cross g')(G) \in \calE_Y
\]
and, since $f \circ g' \closeequiv f \circ g \closeequiv f' \circ g \closeequiv 
f' \circ g'$, we have
\[
    (f \cross f')(F) = ((f \circ g') \cross (f' \circ g'))(F) \in \calE_X,
\]
so $F \in \calE_{\OBJequalizer_{[f],[f']}}$. Hence $\tilde{g}$ is close to 
$\tilde{g}'$, as required.

If $X$ and $Y$ are connected, then $\OBJequalizer_{[f],[f']}$ is clearly 
connected. Moreover, if $X$ is connected, then $\OBJequalizer_{[f],[f']} = Y$ 
as a set and hence is nonempty if $Y$ is nonempty.
\end{proof}

\begin{remark}
The above construction does not work in the unital coarse categories because 
the equalizing pull-back coarse structures are not in general unital (and one 
cannot ``unitalize'' them and still have the required properties). Again, this 
does not imply the nonexistence of equalizers in $\CATUniCrs$, and I do not 
know which equalizers exist in $\CATUniCrs$.
\end{remark}

\begin{remark}
When $X$ is a coarse space and $f, f' \from Y \to X$ are just set maps, one can 
take
\[
    \calE_Y \defeq f^* \calE_X \intersect (f')^* \calE_X
\]
and apply the above Proposition. If $g \from Z \to Y$ is another set map, one 
can then take $\calE_Z \defeq g^* \calE_Y$. Then $f \circ g$ is close to $f' 
\circ g$ if and only if $g$ factors through the equalizer of $[f]$ and $[f']$.
\end{remark}

We have now shown that the nonunital coarse categories have all nonzero 
products and all equalizers. It follows, using a standard argument, that these 
categories have all nonzero limits. For completeness, we give this argument.

\begin{theorem}\label{thm:Crs-lim}
The nonunital coarse categories $\CATCrs$, $\CATConnCrs$, and 
$\CATne{\CATConnCrs}$ have all nonzero limits.
\end{theorem}

\begin{proof}
Let $\calC$ be one of the above categories and suppose $\calF_X \from \calJ \to 
\calC$ ($\calJ$ nonzero, small) is a functor, putting $X_j \defeq \calF_X(j)$ 
for $j \in \Obj(\calJ)$ as usual. If $\calJ$ has no arrows (i.e., $\Map(\calJ) 
= \emptyset$), then $\pfx{\calC}\OBJlim \calF_X$ is just a product, and we are 
done.

Otherwise, let
\begin{align*}
    Y & \defeq \pfx{\calC}\prod_{\mathclap{j \in \Obj(\calJ)}} \; X_j
\shortintertext{and}
    X & \defeq \,\pfx{\calC}\prod_{\mathclap{u \in \Map(\calJ)}}\;\,
            X_{\target(u)}.
\end{align*}
We have two collections of arrows $[f_u], [f'_u] \from Y \to X_{\target(u)}$, 
$u \in \Map(\calJ)$:
\[
    [f_u] \defeq [\pi_{\target(u)}]
\quad\text{and}\quad
    [f'_u] \defeq \calF_X(u) \circ [\pi_{\source(u)}].
\]
By the universal property of products, these collections of arrows give rise to 
canonical arrows $[f] \from Y \to X$ and $[f'] \from Y \to X$, respectively. 
Put
\[
    \pfx{\calC}\OBJlim \calF_X \defeq \OBJequalizer_{[f],[f']},
\]
with the cone $[nu] \from \pfx{\calC}\OBJlim \calF_X \to \calF_X$ being defined 
by $[\nu_j] \defeq [\pi_j] \circ [\equalizer_{[f],[f']}]$ for $j \in 
\Obj(\calJ)$. It is easy to check that $[\nu]$ is indeed a limiting cone.
\end{proof}

\begin{remark}
It follows from the above proof that, as a set, one can always take the limit 
$\OBJlim \calF_X$ to be a subset of the product (set) $\prod_{j \in 
\Obj(\calJ)} X_j$. When all the coarse spaces $X_j$ are connected (i.e., in 
$\CATConnCrs$), one can take $\OBJlim \calF_X$ to be, as a set, exactly the set 
product. Moreover, the proof actually gives a concrete description of limits in 
the coarse categories. If all the $X_j$ are connected, the coarse structure on
\[
    Y \defeq \OBJlim \calF_X
        \defeq \pfx{\CATSet}\prod_{\mathclap{j \in \Obj(\calJ)}}\; X_j
\]
consists of all $F \in \calE_{|Y|_1}$ such that, for all arrows $u \in 
\Map(\calJ)$ and (all) representative coarse maps $f_u \from X_{\source(u)} \to 
X_{\target(u)}$ of $\calF_X(u)$:
\begin{enumerate}
\item $((f_u \circ \pi_{\source(u)}) \cross \pi_{\target(u)}) |_F$ is proper; 
    and
\item $((f_u \circ \pi_{\source(u)}) \cross \pi_{\target(u)})(F)$ is an 
    entourage of $X_{\target(u)}$.
\end{enumerate}
(By taking $u$ to be the identity arrow of $j \in \Obj(\calJ)$, one gets that 
the $\pi_j$ are coarse for $F$.)
\end{remark}

\subsection{Entourages as subspaces of products}\label{ent-subsp-prod}

Is there a relation between entourages of a coarse space $X$, which are subsets 
of $X^{\cross 2} \defeq X \cross X$, and the product coarse space $X \cross X$? 
We first need a coarse space $\Terminate(X)$ which we will discuss more 
thoroughly in \S\ref{subsect:Crs-Term}: For any $X$, $\Terminate(X) \defeq X$ 
as a set, with coarse structure
\[
    \calE_{\Terminate(X)} \defeq \set{E \in \calE_{|X|_1}
            \suchthat 1_{E \cdot X}, 1_{X \cdot E} \in \calE_X}.
\]
Note that if $X$ is unital, $\Terminate(X) = |X|_1$.

The following will be useful later in conjunction with various universal 
properties, as well as generalized coarse quotients (which we intend to study 
in \cite{crscat-quot}).

\begin{proposition}\label{prop:ent-prod}
Suppose $X$ is a coarse space. If $E \in \calE_{\Terminate(X)}$, then $E$ can 
be considered as a unital subspace $|E|$ of the product coarse space $X \cross 
X$. If in fact $E \in \calE_X$, then the restricted projections $\pi_1 |_{|E|}, 
\pi_2 |_{|E|} \from |E| \to X$ are close. Conversely, any unital subspace $|E| 
\subseteq X \cross X$ determines a subset $E \in \calE_{\Terminate(X)} 
\subseteq \calE_{|X|_1}$; if $\pi_1 |_{|E|}$, $\pi_2 |_{|E|}$ are close, then 
$E \in \calE_X$.
\end{proposition}

\begin{proof}
If $E \in \calE_{\Terminate(X)}$, then $1_{|E|}$ is an entourage of $X \cross 
X$: certainly $1_{|E|} \in \calE_{|X \cross X|_1}$, and $(\pi_1 \cross 
\pi_1)(1_{|E|}) = 1_{E \cdot X}$ and $(\pi_2 \cross \pi_2)(1_{|E|}) = 1_{X 
\cdot E}$ are entourages of $X$. If $E \in \calE_X$, then $(\pi_1 \cross 
\pi_2)(1_{|E|}) = E \in \calE_X$; since $|E|$ is unital, it follows that the 
restricted projections are close.

Conversely, suppose $|E| \subseteq X \cross X$ is a unital subspace. Then the 
restricted projections $\pi_1 |_E = \pi_1 |_{|E|}$ and $\pi_2 |_E = \pi_2 
|_{|E|}$ are proper, so $E \in \calE_{|X|_1}$. Since $\pi_1$ maps unital 
subspaces of $X \cross X$ to unital subspaces of $X$ and $\pi_1(|E|) = E \cdot 
X$, the left support, and symmetrically the right support, of $E$ is a unital 
subspace of $X$, and so $E \in \calE_{\Terminate(X)}$. If the restricted 
projections are close, then $E = (\pi_1 \cross \pi_2)(1_{|E|}) \in \calE_X$.
\end{proof}

\subsection{Colimits in the coarse categories}

We now do the same for coproducts, coequalizers, and thus colimits in the 
coarse categories.

\begin{proposition}
Suppose that $\calC$ is one of the coarse categories $\CATCrs$, $\CATConnCrs$, 
$\CATUniCrs$, or $\CATConnUniCrs$, that $\calP\calC$ is the corresponding 
precoarse category, and that $\set{Y_j \suchthat j \in J}$ ($J$ some index set) 
is a collection of coarse spaces in $\calC$ (or $\calP\calC$). The coproduct of 
the $Y_j$ in $\calC$ is just the coarse space
\[
    Y \defeq \pfx{\calP\calC}\coprod_{j \in J} Y_j
\]
(coproduct in $\calP\calC$) together with the ``inclusions'' $[\iota_j] \from 
Y_j \to Y$, $j \in J$ (closeness classes of the inclusions). If instead $\calC 
= \CATne{\CATConnCrs}$, then the same holds except when $J = \emptyset$, in 
which case the coproduct is any one-point coarse space. Thus all the coarse 
categories have all coproducts.
\end{proposition}

\begin{proof}
We have shown (or at least mentioned, in the unital cases) the existence of the 
corresponding coproduct cone $\iota$ in the corresponding precoarse category, 
leaving aside the special case of $\calC = \CATne{\CATConnCrs}$ and $J = 
\emptyset$ (which is easily handled). The quotient functor yields a cone 
$[\iota]$ in the coarse category $\calC$; we must show that it is universal.

Suppose $X$ is a coarse space and $[\mu_j] \from Y_j \to X$, $j \in J$, is a 
collection of arrows in $\calC$. Choosing representative coarse maps $\mu_j$, 
we get a natural coarse map $t \from Y \to X$ such that $\mu_j = t \circ 
\iota_j$ (and hence $[\mu_j] = [\pi_j \circ [t]$) for all $j$. We must show 
this $[t]$ is unique. Suppose $t' \from Y \to X$ is such that $\mu_j 
\closeequiv t' \circ \iota_j$ for all $j$. The coarse structure on the 
precoarse coproduct $Y$ is generated by $F \defeq (\iota_j)^{\cross 2}(F_j)$, 
$F_j \in \calE_{Y_j}$, $j \in J$, and so to show $t \closeequiv t'$ it is 
enough to show that $(t \cross t')(F) \in \calE_X$ for such $F$. But
\[
    (t \cross t')(F) = ((t \circ \iota_j) \cross (t' \circ \iota_j))(F_j)
\]
is in $\calE_X$ since $t \circ \iota_j = \mu_j \closeequiv t' \circ \iota_j$, 
as required.
\end{proof}

Next, coequalizers: Unlike coproducts, coequalizers in the coarse categories 
differ from coequalizers in the precoarse categories; in particular, they 
always exist.

\begin{definition}
Suppose $Y$ is a coarse space and $f, f' \from Y \to X$ ($X$ some set) are 
locally proper maps. The \emph{coequalizing push-forward coarse structure} 
$(f,f')_* \calE_Y$ (on $X$ along $f$ and $f'$) is
\[
    (f,f')_* \calE_Y \defeq \langle f_* \calE_Y, (f')_* \calE_Y,
            \set{(f \cross f')(F) \suchthat F \in \calE_Y} \rangle_Y.
\]
(We may similarly define connected, unital, and connected unital versions.)
\end{definition}

By Proposition~\ref{prop:loc-prop-prod}, the sets $(f \cross f')(F)$ satisfy 
the properness axiom. The coequalizing push-forward coarse structure makes $f$ 
and $f'$ \emph{close} coarse maps, and is the minimum coarse structure on $X$ 
for which this is true.

\begin{definition}
Suppose $f, f' \from Y \to X$ are coarse maps. The \emph{coequalizer} of $[f]$ 
and $[f']$ is $\OBJcoequalizer_{[f],[f']} \defeq X$ equipped the coarse 
structure
\[
    \calE_{\OBJcoequalizer_{[f],[f']}}
        \defeq \langle \calE_X, (f,f')_* \calE_Y \rangle_X,
\]
together with the closeness class of ``identity'' map
\[
    \coequalizer_{[f],[f']} \from X \to \OBJcoequalizer_{[f],[f']}
\]
(which is a coarse map).
\end{definition}

Observe that if $X$ is unital so too is the coequalizer, and similarly if $X$ 
is connected.

\begin{lemma}
Suppose $f, f' \from Y \to X$ are coarse maps. The coequalizer of $[f]$ and 
$[f']$ is coarsely invariant (hence the notation).
\end{lemma}

\begin{proof}
Suppose $e, e' \from Y \to X$ are close to $f$, $f'$, respectively. Observe 
that, since $f_* \calE_Y, (f')_* \calE_Y \subseteq \calE_X$,
\[
    \calE_{\OBJcoequalizer_{[f],[f']}} = \langle \calE_X,
            \set{(f \cross f')(F) \suchthat F \in \calE_Y} \rangle_X
\]
and similarly for $e$ and $e'$. Thus it suffices to show
\[
    \set{(e \cross e')(F) \suchthat F \in \calE_Y}
        \subseteq \calE_{\OBJcoequalizer_{[f],[f']}}
\]
and similarly symmetrically. But if $F \in \calE_Y$, then
\[
    (e \cross e')(F) \subseteq (e \cross f)(1_{F \cdot Y})
            \circ (f \cross f')(F) \circ (f' \cross e')(1_{Y \cdot F})
\]
is in $\calE_{\OBJcoequalizer_{[f],[f']}}$, as required.
\end{proof}

\begin{proposition}
The coequalizer of $[f], [f'] \from Y \to X$ really is (in the categorical 
sense) the coequalizer of $[f]$ and $[f']$ in $\CATCrs$ (or in $\CATConnCrs$, 
$\CATne{\CATConnCrs}$, $\CATUniCrs$, or $\CATConnUniCrs$, as appropriate), 
hence the terminology. Thus $\CATCrs$ (and the other coarse categories) have 
all coequalizers of pairs of arrows.
\end{proposition}

\begin{proof}
Fix representative coarse maps $f$ and $f'$, and suppose $g \from X \to W$ is a 
coarse map such that $g \circ f \closeequiv g \circ f'$. Let $\utilde{g} \from 
\OBJcoequalizer_{[f],[f']} \to W$ be the same, as a set map, as $g$; then 
clearly $g = \utilde{g} \circ \coequalizer_{[f],[f']}$, and hence $[g] = 
[\utilde{g}] \circ [\coequalizer_{[f],[f']}]$, assuming $\utilde{g}$ is 
actually a coarse map. To show that $\utilde{g}$ is coarse, it suffices to show 
that $\utilde{g}$ coarse for sets $E \defeq (f \cross f')(F)$, $F \in \calE_Y$. 
Since
\[
    ((g \circ f) \cross (g \circ f')) |_F
        = g^{\cross 2} |_E \circ (f \cross f') |_F^E
\]
is proper (Proposition~\ref{prop:loc-prop-prod}), it follows that 
$\utilde{g}^{\cross 2} |_E = g^{\cross 2} |_E$ is proper, hence $\utilde{g}$ is 
locally proper for $E$. Since $g \circ f$ and $g \circ f'$ are close, it 
follows that $\utilde{g}$ preserves $E$.

Uniqueness of $[\utilde{g}]$: Suppose $\utilde{g}' \from 
\OBJcoequalizer_{[f],[f']} \to W$ is a coarse map such that $g \closeequiv 
\utilde{g}' \circ \coequalizer_{[f],[f']}$. To show that $\utilde{g}$ is close 
to $\utilde{g}'$, we must show that $(\utilde{g} \cross \utilde{g}')(E) \in 
\calE_W$ for all $E \in \calE_{\OBJcoequalizer_{[f],[f']}}$. Clearly, this is 
the case for $E \in \calE_X \subseteq \calE_{\OBJcoequalizer_{[f],[f']}}$, so 
it suffices to show this for $E = (f \cross f')(F)$ for some $F \in \calE_Y$. 
The map $g' \defeq \utilde{g}' \circ \coequalizer_{[f],[f']}$ is close to $g$, 
hence $g \circ f \closeequiv g' \circ f'$. Therefore,
\[
    (\utilde{g} \cross \utilde{g}')((f \cross f')(F))
        = ((g \circ f) \cross (g' \circ f'))(F),
\]
is in $\calE_W$, as required.

As previously noted, if $X$ is connected, unital, and/or nonempty, then 
$\OBJequalizer_{[f],[f']}$ has the corresponding property or properties, so the 
above actually proves the result in all the coarse categories.
\end{proof}

Since the coarse categories have all coproducts and coequalizers, we 
immediately get the following.

\begin{theorem}\label{thm:Crs-colim}
The coarse categories $\CATCrs$, $\CATConnCrs$, $\CATne{\CATConnCrs}$, 
$\CATUniCrs$, and $\CATConnUniCrs$ have all colimits.
\end{theorem}

\subsection{The termination functor}\label{subsect:Crs-Term}

For essentially set theoretic reasons, $\CATCrs$ does not have a terminal 
object (Corollary~\ref{cor:Crs-no-term}). However, for many purposes, one can 
find a suitable substitute. We begin with some general definitions which are 
applicable in any category $\calC$.

\begin{definition}
In $\calC$, an object $\tilde{X}$ \emph{terminates} an object $X$ if:
\begin{enumerate}
\item there is a (unique) arrow $\tau_X \from X \to \tilde{X}$; and
\item for all $Y \in \Obj(\calC)$, there is at most one arrow $Y \to \tilde{X}$.
\end{enumerate}
I.e., $\tilde{X}$ is terminal in the full subcategory of $\calC$ consisting of 
$X$ and all objects mapping to $\tilde{X}$. $\tilde{X}$ \emph{universally 
terminates} $X$ if it is the smallest object terminating $X$ (i.e., for all 
$\tilde{X}'$ terminating, $X$ there is an arrow $\tilde{X} \to \tilde{X}'$).
\end{definition}

If $\tilde{X}$ terminates $X$, then for all $Y$ and pairs of arrows $f, g \from 
Y \to X$, $\tau_X \circ f = \tau_X \circ g$. Two objects universally 
terminating $X$ are canonically and uniquely isomorphic. If $\tilde{X}$ 
terminates any object, then it universally terminates itself.

In a category with a terminal object $1$, the product of any object $Y$ and $1$ 
is just $Y$. The following generalizes this.

\begin{proposition}\label{prop:term-id}
If there is some arrow $f \from Y \to X$ in $\calC$ and $\tilde{X}$ 
terminates $X$ in $\calC$, then $Y$ is the (categorical) product of 
$\tilde{X}$ and $Y$ (in $\calC$).
\end{proposition}

\begin{proof}
The two ``projections'' from $Y$ are $\pi_{\tilde{X}} \defeq \tau_X \circ f 
\from Y \to \tilde{X}$ and $\pi_Y \defeq \id_Y \from Y \to Y$. Suppose $Z \in 
\Obj(\calC)$ is equipped with arrows $p_{\tilde{X}} \from Z \to \tilde{X}$ and 
$p_Y \from Z \to Y$. Both these arrows factor through $p_Y$: evidently $p_Y = 
\pi_Y \circ p_Y$, but also $p_{\tilde{X}} = \pi_{\tilde{X}} \circ p_Y$ since 
there is only one arrow $Z \to \tilde{X}$.
\end{proof}

If $\calC$ is known to have products (of pairs of objects), we can restate the 
above Proposition in the following way: Whenever there is an arrow $f \from Y 
\to X$ and $\tilde{X}$ terminates $X$, the projection $\pi_Y \from 
\tilde{X} \cross Y \to Y$ is an isomorphism. Moreover, one the inverse 
isomorphism is given by the composition
\[
    Y \nameto{\smash{\Delta_Y}} Y \cross Y
        \nameto{\smash{(\tau_X \circ f) \cross \id_Y}} \tilde{X} \cross Y.
\]

\begin{definition}
A \emph{termination functor} on $\calC$ is a functor $\calC \to \calC$ 
(temporarily denoted $X \mapsto \tilde{X}$) which sends each $X$ to an object 
$\tilde{X}$ terminating $X$; such a functor is \emph{universal} if $\tilde{X}$ 
always universally terminates $X$.
\end{definition}

The following is implied: Whenever there is an arrow $f \from Y \to X$, there 
is a unique arrow $\tilde{Y} \to \tilde{X}$ (namely $\tilde{f}$). Note that 
universality is meant in the ``pointwise'' sense, and we do not assert 
universality as a termination functor. Universal termination functors are 
unique up to natural equivalence. Also observe that universal termination 
functors are idempotent up to natural equivalence.

\begin{example}
If $\calC$ has a terminal object $1$, then $1$ terminates all objects, and $X 
\mapsto 1$ is a termination functor (not necessarily universal). In 
$\CATpt{\CATSet}$ or $\CATpt{\CATTop}$ (pointed sets or topological spaces, 
respectively), the functor $X \mapsto \ast$, where $\ast$ is any one-point set 
or space, is a universal termination functor. More generally, in any category 
$\calC$ with a zero object $0$ (i.e., $0$ is initial and terminal), $X \mapsto 
0$ is a universal termination functor.
\end{example}

\begin{example}\label{ex:Set-Top-univ-term}
In $\CATSet$ or $\CATTop$, the functor given by
\[
    X \mapsto \begin{cases}
            \emptyset & \text{if $X = \emptyset$, or} \\
            \ast & \text{if $X \neq \emptyset$,}
        \end{cases}
\]
is a universal termination functor.
\end{example}

\begin{example}
In $\CATCrs$ (and our various full subcategories), $|X|_1$ terminates any 
coarse space $X$ (Proposition~\ref{prop:term-close}). However, $X \mapsto 
|X|_1$ does not define a functor on $\CATPCrs$ (or $\CATCrs$). E.g., for any 
set $X$, there is always a (unique) coarse map from $|X|_0$ to a one-point 
coarse space $\ast$, but no coarse map $|X|_1 = |\,|X|_0\,|_1 \to \ast$ when 
$X$ is infinite. The problem is that coarse maps from $|X|_1$ must be globally 
proper; in the unital categories this is not a problem, so $X \mapsto |X|_1$ 
does define a coarsely invariant functor $\CATUniPCrs \to \CATUniPCrs$ (for 
example). The induced functor on unital coarse category $\CATUniCrs$ is a 
universal termination functor. We wish to generalize this to all of $\CATCrs$.
\end{example}

We recall the definition of the coarse space $\Terminate(X)$ (for $X$ a coarse 
space) from \S\ref{ent-subsp-prod}, and extend $\Terminate$ to a functor in the 
obvious way.

\begin{definition}
For any coarse space $X$, $\Terminate(X)$ is the coarse space which is just $X$ 
as a set with coarse structure
\[
    \calE_{\Terminate(X)} \defeq \set{E \in \calE_{|X|_1}
            \suchthat 1_{E \cdot X}, 1_{X \cdot E} \in \calE_X};
\]
$\tau_X \from X \to \Terminate(X)$ is the ``identity'' map. If $f \from Y \to 
X$, $\Terminate(f) \from \Terminate(Y) \to \Terminate(X)$ is the same as $f$ as 
a set map.
\end{definition}

Observe the following:
\begin{enumerate}
\item $E \subseteq X^{\cross 2}$ is an entourage of $\Terminate(X)$ if and only 
    if $E$ satisfies the properness axiom (i.e., $E \in \calE_{|X|_1}$) and the 
    left and right supports of $E$ are unital subspaces of $X$.
\item $\Terminate(X)$ has the same unital subspaces as $X$ and is the maximum 
    coarse structure on $X$ with this property. (Consequently, if $X$ is 
    unital, $\Terminate(X) = |X|_1$. It also follows that $\Terminate$ is 
    idempotent, and hence so too is the induced functor $[\Terminate]$; see 
    below.)
\end{enumerate}

\begin{proposition}
$\Terminate$ is a coarsely invariant functor $\CATPCrs \to \CATPCrs$. The 
induced functor $[\Terminate] \from \CATCrs \to \CATCrs$ is a universal 
termination functor.
\end{proposition}

\begin{proof}
That $\Terminate(f)$ is a coarse map follows from the above observations, and 
hence $\Terminate(f)$ is a functor. Moreover, using the above observations, we 
see that, for all $X$, all coarse maps to $\Terminate(X)$ are close. In 
particular, this implies first that $\Terminate$ is coarsely invariant and 
second that $[\Terminate]$ is a termination functor on $\CATCrs$.

It only remains to show universality. Suppose $\tilde{X}$ terminates $X$, so 
there is a unique $[t] \from X \to \tilde{X}$, represented by a coarse map $t$, 
say. It suffices to show that there is a coarse map $t' \from \Terminate(X) \to 
\tilde{X}$; since $\tilde{X}$ terminates $X$ in $\CATCrs$, uniqueness of $[t']$ 
follows, as does the equality $[t] = [t'] \circ [\tau_X]$.

Take $t' \defeq t \from \Terminate(X) = X \to \tilde{X}$ as a set map. Local 
properness of $t'$ follows from the above observations and 
Proposition~\ref{prop:loc-prop}\enumref{prop:loc-prop:II}. To see that $t'$ 
preserves entourages, we use Proposition~\ref{prop:ent-prod}: If $E \in 
\calE_{\Terminate(X)}$, consider the unital subspace $|E|$ of the product 
coarse space $X \cross X$. Since $\tilde{X}$ terminates $X$, $t \circ \pi_1 
|_{|E|} \closeequiv t \circ \pi_2 |_{|E|}$, and hence
\[
    ((t \circ \pi_1 |_{|E|}) \cross (t \circ \pi_2 |_{|E|}))(1_{|E|})
        = (t')^{\cross 2}(E)
\]
is an entourage of $\tilde{X}$, as required.
\end{proof}

In the above proof, one could instead consider the map $\Terminate(t) \from 
\Terminate(X) \to \Terminate(\tilde{X})$, and show that $\Terminate(\tilde{X}) =
\tilde{X}$.

\begin{remark}
$\Terminate$ restricts to (coarsely invariant) endofunctors on the other 
precoarse categories, and hence $[\Terminate]$ restricts to universal 
termination functors on the other coarse categories. (The proof of the above 
Proposition requires only unital coarse spaces $|E|$ and not actually the 
nonunital products $X \cross X$, and hence works even in the unital cases.) Of 
course, in the unital cases, $\Terminate$ is just the functor $X \mapsto 
|X|_1$.
\end{remark}

By applying Proposition~\ref{prop:term-id}, we immediately get the following, 
which will play a crucial role in the development of exponential objects in the 
coarse categories \cite{crscat-II}.

\begin{corollary}\label{cor:Crs-term-id}
If there is a coarse map $Y \to \Terminate(X)$, where $X$ and $Y$ are coarse 
spaces, then
\[
    \pi_Y \from \Terminate(X) \cross Y \to Y
\]
is a coarse equivalence. The maps
\[
    D_\tau \defeq (\tau \cross \id_Y) \circ \Delta_Y
            \from Y \to \Terminate(X) \cross Y,
\]
where $\tau \from Y \to \Terminate(X)$ is any coarse map (they are all close), 
are coarsely inverse to $\pi_Y$. Hence, if there is a coarse map $Y \to 
\Terminate(X)$, then $Y \cong \Terminate(X) \cross Y$ canonically in $\CATCrs$ 
(or in $\CATConnCrs$ or $\CATne{\CATConnCrs}$). In the case $Y \defeq X$, we 
get that $\pi_X \from \Terminate(X) \cross X \to X$ and $D_X \defeq D_{\tau_X} 
\from X \to \Terminate(X) \cross X$ are coarsely inverse coarse equivalences, 
so $X \cong \Terminate(X) \cross X$ canonically in $\CATCrs$ (or in 
$\CATConnCrs$ or $\CATne{\CATConnCrs}$).
\end{corollary}

\begin{remark}\label{rmk:term-unital-prod}
For any set $X$, $\Terminate(|X|_1) = |X|_1$, so $|X|_1 \cross |X|_1$ is 
(canonically) coarsely equivalent to $|X|_1$. While $|X|_1$ is always unital, 
$|X|_1 \cross |X|_1$ is unital only when $X$ is finite. In particular, 
unitality is \emph{not} coarsely invariant. It also follows easily that $|X|_1$ 
is actually the product of $|X|_1$ with itself in the unital coarse category 
$\CATUniCrs$. More generally, for any coarse space $X$, the product of $X$ and 
$|X|_1$ in $\CATUniCrs$ is just $X$. (As previously mentioned, $\CATUniCrs$ has 
some products of infinite spaces, even though the natural construction of the 
corresponding products in $\CATCrs$ are nonunital.)
\end{remark}

\subsection{Monics and images}

\begin{example}
Pull-back coarse structures are not coarsely invariant. That is, suppose $f, f' 
\from Y \to X$ are coarse maps. Even if $f \closeequiv f'$, it may not be the 
case that $f^* \calE_X = (f')^* \calE_X$. To see this, take $Y \defeq 
|\setN|_0^\TXTconn$, $X \defeq |\setN|_1$, $f$ to be the ``identity'' map (as a 
set map), and $f'$ to be a constant map. Then $f^* \calE_X = \calE_{|Y|_1}$ 
whereas $(f')^* \calE_X = \calE_Y$.
\end{example}

\begin{proposition}
If $f, f' \from Y \to X$ are coarse maps with $f \closeequiv f'$, then
\[
    \calE_{\Terminate(Y)} \intersect f^* \calE_X
        = \calE_{\Terminate(Y)} \intersect (f')^* \calE_X.
\]
\end{proposition}

\begin{proof}
We prove inclusion $\subseteq$; containment $\supseteq$ follows symmetrically. 
Suppose $F \in \calE_{\Terminate(Y)} \intersect f^* \calE_X$. Since $F \in 
\Terminate(Y)$, $f'$ is locally proper for $F$ 
(Proposition~\ref{prop:loc-prop}\enumref{prop:loc-prop:II}). It only remains to 
show that $(f')^{\cross 2}(F) \in \calE_X$. But
\[\begin{split}
    (f')^{\cross 2}(F) \subseteq (f' \cross f)(1_{F \cdot Y})
            \circ f^{\cross 2}(F) \circ (f \cross f')(1_{Y \cdot F}) \in \calE_X
\end{split}\]
since $f \closeequiv f'$ (and the left and right supports of $F$ are unital 
subspaces of $Y$) and $f^{\cross 2}(F) \in \calE_X$.
\end{proof}

\begin{definition}
Suppose $[f] \from Y \to X$. The \emph{coarsely invariant pull-back coarse 
structure} $[f]^* \calE_X$ on $Y$ (along $[f]$) is given by
\[
    [f]^* \calE_X \defeq \calE_{\Terminate(Y)} \intersect f^* \calE_X
\]
(where $f \from Y \to X$ is any representative coarse map).
\end{definition}

\begin{proposition}\label{prop:Crs-factor-I}
If $[f] \from Y \to X$ is represented by a coarse map $f$, then $[f]$ factors 
as
\[
    Y \nameto{\smash{[\beta]}} |Y|_{[f]^* \calE_X}
        \nameto{\smash{[\utilde{f}]}} X,
\]
where $\beta = \id_Y$ and $\utilde{f} = f$ as set maps (i.e., $\calE_Y 
\subseteq [f]^* \calE_X$). Moreover, $[\utilde{f}]$ depends only on $[f]$ (and 
not on the particular $f$) and is unique in the above factorization.
\end{proposition}

\begin{proof}
The factorization follows immediately from Corollary~\ref{cor:crs-factor-I}. We 
now show that $f \closeequiv f'$ implies $\utilde{f} \closeequiv \utilde{f}'$ 
(noting that $[f]^* \calE_X = [f']^* \calE_X$). If $F \in [f]^* \calE_X$, then
\[\begin{split}
    (f \cross f')(F) & = (f \cross f')(F \circ 1_{Y \cdot F}) \\
        & \subseteq f^{\cross 2}(F) \circ (f \cross f')(1_{Y \cdot F})
\end{split}\]
is in $\calE_X$ since $f^{\cross 2}(F) \in \calE_X$ and $1_{Y \cdot F} \in 
\calE_Y$ so $(f \cross f')(1_{Y \cdot F}) \in \calE_X$ as $f \closeequiv f'$. 
Uniqueness: If $[f] = [g] \circ [\beta]$, where $[g] \from |Y|_{[f]^* \calE_X} 
\to X$ and $g$ is any representative, then $f \closeequiv g \circ \beta$, so 
$\utilde{f} \closeequiv (g \circ \beta)\utilde{\mathstrut} = g$.
\end{proof}

\begin{proposition}\label{prop:Crs-monic}
$[f] \from Y \to X$ is monic in $\CATCrs$ if and only if $\calE_Y = [f]^* 
\calE_X$ (i.e., if and only if $Y = |Y|_{[f]^* \calE_X}$).
\end{proposition}

\begin{proof}
Fix a representative coarse map $f \from Y \to X$ and let $Y 
\nameto{\smash{\beta}} |Y|_{[f]^* \calE_X} \nameto{\smash{\utilde{f}}} X$ be 
the canonical factorization.

(\textimplies): Suppose there exists some $F \in [f]^* \calE_X \setminus 
\calE_Y$. Consider $|F|$ as a unital subspace of the product $Y \cross Y$, with 
projections $\pi_1 |_{|F|}, \pi_2 |_{|F|} \from |F| \to Y$. Then
\[
    (\pi_1 |_{|F|} \cross \pi_2 |_{|F|})(1_{|F|}) = F,
\]
so $\pi_1 |_{|F|}$ is not close to $\pi_2 |_{|F|}$, but $\beta \circ \pi_1 
|_{|F|}$ is close to $\beta \circ \pi_2 |_{|F|}$. Hence $[\pi_1 |_{|F|}] \neq 
[\pi_2 |_{|F|}]$ but
\[
    [f] \circ [\pi_1 |_{|F|}] = [\utilde{f}] \circ [\beta] \circ [\pi_1 |_{|F|}]
        = [\utilde{f}] \circ [\beta] \circ [\pi_2 |_{|F|}]
        = [f] \circ [\pi_2 |_{|F|}],
\]
so $[f]$ is not monic.

(\textimpliedby): Suppose $g, g' \from Z \to Y$ are coarse maps such that $[f] 
\circ [g] = [f] \circ [g']$. Then, for each $G \in \calE_Z$,
\[
    ((f \circ g) \cross (f \circ g'))(G)
        = f^{\cross 2}((g \cross g')(G)) \in \calE_X.
\]
But then $(g \cross g')(G) \in [f]^* \calE_X = \calE_Y$, so $[g] = [g']$, as 
required.
\end{proof}

\begin{corollary}
For any $[f] \from Y \to X$, the canonical arrow
\[
    [\utilde{f}] \from |Y|_{[f]^* \calE_X} \to X
\]
is monic in $\CATCrs$.
\end{corollary}

\begin{definition}\label{def:Crs-image}
Suppose $[f] \from Y \to X$. Denote $\OBJim [f] \defeq |Y|_{[f]^* \calE_X}$ and 
$\im [f] \defeq [\utilde{f}] \from \OBJim [f] \to X$, where $[\utilde{f}]$ is 
defined as above. We will also sometimes write $[f](Y) \defeq \OBJim [f]$.
\end{definition}

Despite the notation, $[f](Y)$ should not be considered as a subspace of $X$ 
(however, see Proposition~\ref{prop:Crs-subsp-images} and the discussion which 
precedes it).

\begin{theorem}
For any $[f] \from Y \to X$, the subobject of $X$ represented by the arrow $\im 
[f] \from \OBJim [f] \monto X$ is the (categorical) image of $[f]$ in 
$\CATCrs$.
\end{theorem}

\begin{proof}
Suppose $[f]$ also factors as $Y \nameto{\smash{[h]}} Z \nameto{\smash{[g]}} X$ 
where $[g]$ is monic, so that $\calE_Z = [g]^* \calE_X$. We must show that 
there is a unique $[\underline{h}] \from \OBJim [f] \to Z$ such that $\im [f] = 
[g] \circ [\underline{h}]$.

Pick a representative coarse map $h \from Y \to Z$, and put $\underline{h} 
\defeq h$ as a set map $\OBJim [f] = Y \to Z$. First, $\underline{h}$ is a 
coarse map: Local properness is equivalent to properness when restricted to 
unital subspaces (Corollary~\ref{cor:loc-prop-uni}); since $\OBJim [f]$ and $Y$ 
have the same unital subspaces (and $\underline{h} = h$ as set maps), 
$\underline{h}$ is locally proper. Reasoning similarly, for any $F \in 
\calE_{\OBJim [f]} = [f]^* \calE_X$, $\underline{h}^{\cross 2}(F)$ is in 
$\calE_{\Terminate(Z)}$. Then, since $\calE_Z = [g]^* \calE_X$, it follows that 
$\underline{h}$ is coarse. From the uniqueness assertion of 
Proposition~\ref{prop:Crs-factor-I}, we get that $[g] \circ [\underline{h}] = 
\im [f]$. Uniqueness of $[\underline{h}]$: If $[h'] \from \OBJim [f] \to Z$ and 
$[g] \circ [h'] = \im [f] = [g] \circ [h]$, then $[h] = [h']$ since $[g]$ is 
monic.
\end{proof}

\subsection{Epis and coimages}

For rather trivial reasons, push-forward coarse structures are not coarsely 
invariant. Recall that coarse structures are semirings, which gives rise to an 
obvious notion of ideals.

\begin{definition}\label{def:ideal}
Suppose $\calE_X$ is a coarse structure on a set $X$. A subset $\calE \subseteq 
\calE_X$ is an \emph{ideal} of $\calE_X$ if it is a coarse structure on $X$ 
such that $E \circ E', E' \circ E \in \calE$ for all $E \in \calE$, $E' \in 
\calE_X$. Note that any intersection of ideals is again an ideal. The 
\emph{ideal} $\lAngle \calE \rAngle_X$ (of $\calE_X$ generated by $\calE$) is 
the smallest ideal of $\calE_X$ which contains $\calE$.
\end{definition}

\begin{proposition}
Suppose $f, f' \from Y \to X$ are coarse maps with $f \closeequiv f'$. Then
\[
    \lAngle f_* \calE_Y \rAngle_X = \lAngle (f')_* \calE_Y \rAngle_X.
\]
\end{proposition}

\begin{proof}
Elements $E \in \lAngle f_* \calE_Y \rAngle_X$ are exactly subsets
\[
    E \subseteq E' \circ f^{\cross 2}(F) \circ E'' \union E'''
\]
for some $F \in \calE_Y$ and some $E', E'', E''' \in \calE_X$ with $E'''$ 
finite. But then
\[
    E \subseteq (E' \circ (f \cross f')(1_{F \cdot Y})) \circ (f')^{\cross 2}(F)
            \circ ((f' \cross f)(1_{Y \cdot F}) \circ E'') \union E'''
\]
is in $\lAngle (f')_* \calE_Y \rAngle_X$ (and symmetrically) as required.
\end{proof}

\begin{definition}
Suppose $[f] \from Y \to X$. The \emph{coarsely invariant push-forward coarse 
structure} $[f]_* \calE_Y$ on $X$ (along $[f]$) is given by
\[
    [f]_* \calE_Y \defeq \lAngle f_* \calE_Y \rAngle_X
\]
(where $f \from Y \to X$ is any representative coarse map).
\end{definition}

Despite the obvious parallels with coarsely invariant pull-backs, the coarsely 
invariant push-forward $[f]_* \calE_Y$ depends very little on $\calE_Y$. In 
fact, it depends only on the set of unital subspaces of $Y$ (recall from 
Proposition~\ref{prop:close-uni} that closeness is entirely determined on the 
unital subspaces). Thus we have the following.

\begin{proposition}
For any $[f] \from Y \to X$,
\[
    [f]_* \calE_Y = (\im [f])_* \calE_{\OBJim [f]}.
\]
\end{proposition}

\begin{proof}
Recall that $\OBJim [f] \defeq |Y|_{[f]^* \calE_X}$, where $[f]^* \calE_X 
\defeq \calE_{\Terminate(Y)} \intersect f^* \calE_X$ (for any representative 
map $f$) and $\im [f] \defeq f$ as a set map. Since $\calE_Y \subseteq [f]^* 
\calE_X$, $[f]_* \calE_Y \subseteq (\im [f])_* \calE_{\OBJim [f]}$. For the 
opposite inclusion, it suffices to show that, for $F \in [f]^* \calE_X$,
\[
    E \defeq f^{\cross 2}(F) \in \lAngle f_* \calE_Y \rAngle_X;
\]
but $F \cdot Y$ is a unital subspace of $Y$ (hence $1_{F \cdot Y} \in \calE_Y$) 
and $f^{\cross 2}(F) \in \calE_X$, so
\[
    E = f^{\cross 2}(1_{F \cdot Y}) \circ E \in \lAngle f_* \calE_Y \rAngle_X,
\]
as required.
\end{proof}

Suppose $[f] \from Y \to X$, represented by a coarse map $f$. Denote
\[
    X_{[f]} \defeq \set{x \in X
        \suchthat \text{$x$ is connected to some $x' \in f(Y)$}} \subseteq X,
\]
a subspace of $X$. It is easy to see that $X_{[f]}$ really only depends on the 
closeness class $[f]$, as the notation indicates. (If $X$ is connected, then of 
course $X_{[f]} = X$.)

The subspace $X_{[f]} \subseteq X$ contains the set image of $f$ (and indeed of 
any coarse map close to $f$), and hence we may take the range restriction $f 
|^{X_{[f]}}$ which is evidently a coarse map $Y \to X_{[f]}$. It is easy to see 
that the closeness class $[f |^{X_{[f]}}]$ only depends on the closeness class 
$[f]$, and hence we also temporarily denote
\[
    [f] |^{X_{[f]}} \defeq [f |^{X_{[f]}}] \from Y \to X_{[f]}.
\]
Now, we may coarsely invariantly push $\calE_Y$ forward along $[f] |^{X_{[f]}}$ 
to get a coarse space $|X_{[f]}|_{([f] |^{X_{[f]}})_* \calE_Y}$. We get the 
following.

\begin{proposition}\label{prop:Crs-factor-II}
If $[f] \from Y \to X$ is represented by a coarse map $f$, then $[f]$ factors 
as
\[
    Y \nameto{\smash{[\tilde{f}]}} |X_{[f]}|_{([f] |^{X_{[f]}})_* \calE_Y}
        \nameto{\smash{[\alpha]}} X,
\]
where $\tilde{f} = f |^{X_{[f]}}$ and $\alpha$ is the inclusion as set maps 
(thus $([f] |^{X_{[f]}})_* \calE_Y \subseteq \calE_X$). Moreover, $[\tilde{f}]$ 
depends only on $[f]$ (and not $f$) and is unique in the above factorization.
\end{proposition}

\begin{proof}
Nearly all the assertions are clear from the definitions, 
Corollary~\ref{cor:crs-factor-II}, and the previous remarks. We show that $f 
\closeequiv f'$ implies $\tilde{f} \closeequiv \tilde{f}'$: If $F \in \calE_Y$, 
then
\[
    (\tilde{f} \cross \tilde{f}')(F) = (f \cross f')(F)
        \subseteq f^{\cross 2}(F) \circ (f \cross f')(1_{Y \cdot F})
\]
is in $([f] |^{X_{[f]}})_* \calE_Y$ since $f^{\cross 2}(F) \in (f 
|^{X_{[f]}})_* \calE_Y$ and $(f \cross f')(1_{Y \cdot F}) \in \calE_X 
|_{X_{[f]}}$. Uniqueness: If $[f] = [\alpha] \circ [g]$, where $[g] \from Y \to 
|X_{[f]}|_{([f] |^{X_{[f]}})_* \calE_Y}$ and $g$ is any representative, then $f 
\closeequiv \alpha \circ g$, so $\tilde{f} \closeequiv (\alpha \circ 
g)\tilde{\mathstrut} = g$.
\end{proof}

\begin{proposition}\label{prop:Crs-epi}
$[f] \from Y \to X$ is epi in $\CATCrs$ if and only if $X_{[f]} = X$ and $[f]_* 
\calE_Y = \calE_X$ (i.e., if and only if $|X_{[f]}|_{([f] |^{X_{[f]}})_* 
\calE_Y} = X$).
\end{proposition}

\begin{proof}
(\textimplies): Consider the push-out square
\[\begin{CD}
    Y @>{[f]}>> X \\
    @V{[f]}VV @V{[e_1]}VV \\
    X @>{[e_2]}>> X \copro_Y X
\end{CD}\]
(in $\CATCrs$). Fix a representative coarse map $f \from Y \to X$. As a set, 
one may take $X \copro_Y X \defeq X_1 \disjtunion X_2$ (disjoint union of sets) 
where $X_1 \defeq X_2 \defeq X$, with coarse structure
\[
    \langle \calE_{X_1}, \calE_{X_2},
        \set{(f_1 \cross f_2)(F) \suchthat F \in \calE_Y}
            \rangle_{X_1 \disjtunion X_2},
\]
where $\calE_{X_j} \defeq \calE_X \subseteq \powerset((X_j)^{\cross 2})$ and 
$f_j \defeq f \from Y \to X = X_j$, for $j = 1, 2$. As set maps, one may take 
$e_1$, $e_2$ to be the two inclusions.

If $X_{[f]} \neq X$, then there exists $x_0 \in X$ not connected to any $f(y)$, 
$y \in Y$. The entourage $\set{1_{x_0}} \in \calE_X$ then shows that $e_1$ is 
not close to $e_2$, hence $[e_1] \neq [e_2]$ while $[e_1] \circ [f] = [e_2] 
\circ [f]$ so $[f]$ is not epi. Similarly, if $E \in \calE_X \setminus [f]_* 
\calE_Y$, then one can show that $(e_1 \cross e_2)(E)$ is not an entourage of 
$X \copro_Y X$, hence again $[f]$ is not epi.

(\textimpliedby): It suffices to show that $|X_{[f]}|_{([f] |^{X_{[f]}})_* 
\calE_Y} = X$ implies that $[e_1] = [e_2]$ in the push-out square considered 
above. If $|X_{[f]}|_{([f] |^{X_{[f]}})_* \calE_Y} = X$, then every entourage 
of $[f]_* \calE_Y$ is a subset of an entourage of the form $E_1 \circ f^{\cross 
2}(F) \circ E_2$ for $F \in \calE_Y$ and $E_1, E_2 \in \calE_X$. Thus if $[f]_* 
\calE_Y = \calE_X$, given $E \in \calE_X$ choose $F$, $E_1$, and $E_2$ so that 
$E \subseteq E_1 \circ f^{\cross 2}(F) \circ E_2$, and then
\[
    (e_1 \cross e_2)(E) \subseteq E_1 \circ (f_1 \cross f_2)(F) \circ E_2
\]
(where we now consider $E_j \in \calE_{X_j} = \calE_X$ for $j = 1, 2$) is an 
entourage of $X \copro_Y X$. Thus $e_1$ is close to $e_2$ as required.
\end{proof}

\begin{corollary}
For any $[f] \from Y \to X$, the canonical arrow
\[
    [\tilde{f}] \from Y \to |X_{[f]}|_{([f] |^{X_{[f]}})_* \calE_Y}
\]
is epi in $\CATCrs$.
\end{corollary}

\begin{corollary}\label{cor:Crs-epi-crs-structs}
Suppose $\calE$, $\calE'$ are coarse structures on a set $X$ with $\calE' 
\subseteq \calE$. If every unital subspace of $|X|_{\calE}$ is a unital 
subspace of $|X|_{\calE'}$, then the class $[q]$ of the ``identity'' map
\[
    q \from |X|_\calE' \to |X|_\calE
\]
is epi in $\CATCrs$.
\end{corollary}

\begin{proof}
Trivially, $(|X|_\calE)_{[q]} = |X|_\calE$. We have that
\[
    [q]_* \calE' = \lAngle \calE' \rAngle_{|X|_{\calE}}
\]
is an ideal of $\calE$; we must prove equality, so suppose $E \in \calE$. Then 
$1_{E \cdot X}$ is in $\calE$ hence also in $\calE'$, so
\[
    E = 1_{E \cdot X} \circ E
\]
is in $[q]_* \calE'$, as required.
\end{proof}

\begin{definition}\label{def:Crs-coimage}
Suppose $[f] \from Y \to X$. Denote $\OBJcoim [f] \defeq |X_{[f]}|_{([f] 
|^{X_{[f]}})_* \calE_Y}$ and $\coim [f] \defeq [\tilde{f}] \from Y \to \OBJcoim 
[f]$, where $[\tilde{f}]$ is defined as above.
\end{definition}

\begin{theorem}
For any $[f] \from Y \to X$, the quotient object of $Y$ represented by the 
arrow $\coim [f] \from Y \surto \OBJcoim [f]$ is the (categorical) coimage of 
$[f]$ in $\CATCrs$.
\end{theorem}

\begin{proof}
Suppose $[f]$ also factors as $Y \nameto{\smash{[h]}} Z \nameto{\smash{[g]}} X$ 
where $[h]$ is epi, so that $Z_{[h]} = Z$ and $\calE_Z = [h]_* \calE_Y$. We 
must show that there is a unique $[\bar{g}] \from Z \to \OBJcoim [f]$ such that 
$\coim [f] = [\bar{g}] \circ [h]$.

Pick representative coarse maps $g \from Z \to X$ and $h \from Y \to Z$. We may 
then take $f \defeq g \circ h$ as a representative for $[f]$. Since $Z_{[f]} = 
Z$ and $[g] \circ [h] = [f]$, it follows that $g$ has set image contained in 
$X_{[f]}$. Thus we may put $\bar{g} \defeq g |^{X_{[f]}}$ as a set map $Z \to 
X_{[f]} = \OBJcoim [f]$. $\bar{g}$ is a coarse map: It is locally proper since 
$g = \alpha \circ \bar{g}$ is locally proper. Since $\calE_Z = [h]_* \calE_Y$, 
every entourage of $Z$ is contained in one of the form $G_1 \circ h^{\cross 
2}(F) \circ G_2$, for $F \in \calE_Y$, $G_1, G_2 \in \calE_Z$. For such an 
entourage,
\[
    \bar{g}^{\cross 2}(G_1 \circ h^{\cross 2}(F) \circ G_2)
        \subseteq g^{\cross 2}(G_1) \circ (g \circ h)^{\cross 2}(F)
                \circ g^{\cross 2}(G_2)
\]
is in $([f] |^{X_{[f]}})_* \calE_Y$ since $g^{\cross 2}(G_1), g^{\cross 2}(G_2) 
\in \calE_X$ (and $g$ has set image in $X_{[f]}$) and $(g \circ h)^{\cross 
2}(F) = f^{\cross 2}(F)$. Thus $\bar{g}$ is coarse. From the uniqueness 
assertion of Proposition~\ref{prop:Crs-factor-II} (or, since $\tilde{f} = 
\bar{g} \circ h$), we get that $\coim [f] = [\bar{g}] \circ [h]$. Uniqueness of 
$[\bar{g}]$ follows immediately from the hypothesis that $[h]$ is epi.
\end{proof}

\subsection{Monic and epi arrows}

I do not know if $\CATCrs$ is a \emph{balanced} category, i.e., whether every 
arrow in $\CATCrs$ which is both monic and epi is an isomorphism (the converse 
is always true, of course). To show that a monic and epi $[f] \from Y \to X$ is 
an isomorphism one must show that there is an inverse $[f]^{-1} \from X \to Y$. 
When $X$ is unital, this is fairly straightforward (see below), but I do not 
know how to prove it when $X$ is not.

\begin{theorem}\label{thm:Crs-unital-balanced}
If $[f] \from Y \to X$ is monic and epi in $\CATCrs$ and $X$ is a unital coarse 
space, then $[f]$ is an isomorphism in $\CATCrs$.
\end{theorem}

\begin{proof}
Fix a representative coarse map $f \from Y \to X$. Since $[f]$ is epi, by 
Proposition~\ref{prop:Crs-epi}, $X_{[f]} = X$ and
\[
    [f]_* \calE_Y \defeq \lAngle f_* \calE_Y \rAngle_X = \calE_X.
\]
Then every entourage $E_0 \in \calE_X$ is contained in one of the form $E_1 
\circ E_2 \circ E_3$, where $E_1, E_3 \in \calE_X$ and $E_2 \in f_* \calE_Y$. 
Every $E_2 \in f_* \calE_Y$ is contained in an entourage of the form
\[
    \bigl(f^{\cross 2}(F_2^1) \circ \dotsb \circ f^{\cross 2}(F_2^N)\bigr)
        \union \bigunion_{j \in J} (K_j \cross K'_j),
\]
where $F_2^1, \dotsc, F_2^N \in \calE_Y$ (some $N \geq 0$), $J$ is the set of 
connected components of $X$, and $K_j$, $K'_j$ are finite subsets of $j$ for 
each $j \in J$. Since $X_{[f]} = X$ (and $f^{\cross 2}(F_2^k) \in \calE_X$ for 
$k = 2, \dotsc, N$), it follows that every $E \in \calE_X$ is contained in a 
some entourage
\[
    E_0 \circ f^{\cross 2}(F_0) \circ E'_0,
\]
where $E_0, E'_0 \in \calE_X$ and $F_0 \in \calE_Y$.

We specialize the above discussion to the case $E = 1_X$ which is in $\calE_X$ 
by unitality. Fix $E_0, E'_0 \in \calE_X$ and $F_0 \in \calE_Y$, so that $1_X 
\subseteq E_0 \circ f^{\cross 2}(F_0) \circ E'_0$. Define a set map $e \from X 
\to Y$ as follows. For each $x \in X$, there are $x', x'' \in X$ and $y', y'' 
\in Y$ such that $(x,x') \in E_0$, $(x'',x) \in E'_0$, $f(y') = x'$, $f(y'') = 
x''$, and $(y',y'') \in F_0$; choosing such a $y'' \in Y$ in particular, put 
$e(x) \defeq y''$.

We must verify that (any) $e \from X \to Y$ as constructed above is a coarse 
map. Local properness: $X$ is unital, so $e$ is locally proper if and only if 
it is proper. For any $y \in Y$, $e^{-1}(\set{y}) \subseteq (E_0 \circ 
f^{\cross 2}(F_0)) \cdot \set{f(y)}$ is finite, since $E_0 \circ f^{\cross 
2}(F_0) \in \calE_X \subseteq \calE_{|X|_1}$ satisfies the properness axiom. 
$e$ preserves entourages: Fix $E \in \calE_X$ and put $F \defeq e^{\cross 
2}(E)$. Since $[f]$ is monic, by Proposition~\ref{prop:Crs-monic},
\[
    \calE_Y = [f]^* \calE_X \defeq \calE_{\Terminate(Y)} \intersect f^* \calE_X.
\]
Since $e$ is (locally) proper, $F$ satisfies the properness axiom; since the 
image of $e$ is contained in the unital subspace $Y \cdot F_0$ of $Y$, it then 
follows that $F \in \calE_{\Terminate(Y)}$ and hence also that $f$ is locally 
proper for $F$. To show that $F \in f^* \calE_X$, it only remains to show that 
$f^{\cross 2}(F) \in \calE_X$: Since
\[
    G_0 \defeq (\id_X \cross (f \circ e))(1_X)
        \subseteq E_0 \circ f^{\cross 2}(F_0)
\]
is in $\calE_X$,
\[
    f^{\cross 2}(F)
        = (f \circ e)^{\cross 2}(E)
        \subseteq (G_0)^\transpose \circ E \circ G_0
\]
is also in $\calE_X$.

Since $G_0 \in \calE_X$, we also get that $f \circ e$ is close to $\id_X$, 
i.e., $[f \circ e] = [f] \circ [e]$ is the identity arrow $[\id_X]$ of $X$ in 
$\CATCrs$. Since $[e]$ is monic (and $[f] \circ [e] \circ [f] = [f] = [f] \circ 
[\id_Y]$), it also follows that $[e] \circ [f] = [\id_Y]$. Thus $[e] = 
[f]^{-1}$, as required.
\end{proof}

\begin{corollary}
If $[f] \from Y \to X$ is monic and epi in $\CATCrs$ and $X$ is coarsely 
equivalent to a unital coarse space, then $[f]$ is an isomorphism in $\CATCrs$.
\end{corollary}

The problem with the above Corollary, of course, is that I do not know when a 
coarse space is coarsely equivalent to a unital one. If $\iota \from X' \injto 
X$ is the inclusion of a subspace of $X$ into $X$, then $[\iota]$ is monic (and 
$\OBJim [\iota] = X'$), so $\coim [\iota] \from X' \to \OBJcoim [\iota]$ is 
both monic and epi. (If $X$ is connected and $X'$ nonempty, $\OBJcoim [\iota]$ 
is just the set $X$ equipped with the coarse structure of entourages in 
$\calE_X$ ``supported near $X'$''.) However, I do not know when $\coim [\iota]$ 
is a coarse equivalence.

More generally, for any $[f] \from Y \to X$, the natural arrow $Y \to \OBJim 
[f]$ is epi (either use Proposition~\ref{prop:Crs-epi}, or the fact that 
$\CATCrs$ has equalizers and, e.g., \cite{MR0202787}*{Ch.~I Prop.~10.1}) and 
hence there is a natural epi arrow $[\gamma] \from \OBJim [f] \to \OBJcoim [f]$ 
through which $\im [f] \from \OBJim [f] \to X$ factors; as $\im [f]$ is monic, 
$[\gamma]$ must also be monic. (One may dually show that the natural arrow 
$\OBJcoim [f] \to X$ is monic, but this yields the same arrow $[\mu]$.) Of 
course, I do not know when $[\mu]$ is an isomorphism. But when it is an 
isomorphism, one can, in a coarsely invariant way, describe the image of $[f]$ 
as a subset of $X$ with a certain coarse structure. This would be an appealing 
``generalization'' of the following, which is not coarsely invariant in the 
desired sense.

\begin{proposition}\label{prop:Crs-subsp-images}
If $f \from Y \to X$ is a coarse map and $Y$ is unital, then $\OBJim [f] = 
f(Y)$ (where $f(Y)$ is the subspace of $X$ determined by the set image of $f$) 
as subobjects of $X$ in $\CATCrs$.
\end{proposition}

\begin{proof}
If $Y$ is unital, $X' \defeq f(Y)$ is also unital. The range restriction $f 
|^{X'} \from Y \to X$ is a coarse map, and $[f]^* \calE_X = [f |^{X'}]^* 
\calE_X$ hence $\OBJim [f] = \OBJim [f |^{X'}]$. Using this equality, we get 
$\im [f] = [\iota] \circ \im [f |^{X'}]$, where $\iota \from X' \injto X$ is 
the inclusion. But it is easy to check that $\calE_{X'} \defeq \calE_X |_{X'} = 
[f |^{X'}]_* \calE_Y$, so $[f |^{X'}]$ is epi. Hence $\im [f |^{X'}] \from 
\OBJim [f] = \OBJim [f |^{X'}] \to X'$ is both monic and epi, hence an 
isomorphism by Theorem~\ref{thm:Crs-unital-balanced}.
\end{proof}

\subsection{Quotients of coarse spaces}\label{subsect:Crs-quot}

We now discuss a notion of quotient coarse spaces in $\CATCrs$. The quotient 
spaces below are not the most general possible; rather, they appear to be a 
special case of a more general notion (of quotients by \emph{coarse equivalence 
relations}). However, I have not fully explored the more general notion, and so 
I leave it to a future paper.

Suppose $\calC$ is a category with zero object $0$ (e.g., an abelian category), 
i.e., $0$ is both initial and terminal. Given an arrow $f \from Y \to X$ (often 
taken to be monic) in $\calC$, a standard way of defining the quotient, denoted 
$X/f(Y)$, is as the push-out $X \copro_Y 0$ (assuming it exists); i.e., 
$X/f(Y)$ fits into a push-out square
\[\begin{CD}
    Y @>f>> X \\
    @VVV @VVV \\
    0 @>>> X/f(Y)
\end{CD}\quad.\]
The quotient $X/f(Y)$ comes equipped with an arrow $X \to X/f(Y)$ and, in the 
above case, also an arrow $0 \to X/f(Y)$.

In an abelian category, $X/f(Y)$ is by definition just the cokernel of $f$. If 
$\calC = \CATpt{\CATSet}$ or $\CATpt{\CATTop}$ (pointed sets or spaces), then 
one has $0 = \ast$ (a one-point set/space) and $X/f(Y)$ is (isomorphic to) just 
$X$ with the image of $f$ collapsed to the base point. The situation in 
$\CATSet$ or $\CATTop$ is slightly more complicated: If $Y \neq \emptyset$, one 
can again take the push-out $X/f(Y) \defeq X \copro_Y \ast$. However, if $Y = 
\emptyset$, then $X/f(Y) \cong X$; one should instead take $X/f(Y) \defeq X 
\copro_Y \emptyset$. In other words, one takes $X/f(Y) \defeq X \copro_Y 
\tilde{Y}$, where $\tilde{Y}$ universally terminates $Y$ (see 
Example~\ref{ex:Set-Top-univ-term}). This is exactly what we do in the coarse 
categories.

\begin{definition}
Suppose $[f] \from Y \to X$ (in $\CATCrs$). The \emph{quotient coarse space} 
$X/[f](Y)$ is the push-out $X \copro_Y \Terminate(Y)$ in $\CATCrs$, i.e., 
$X/[f](Y)$ fits into a push-out square
\[\begin{CD}
    Y @>{[f]}>> X \\
    @V{[\tau_Y]}VV @V{[q]}VV \\
    \Terminate(Y) @>{[f]/[f]}>> X/[f](Y)
\end{CD}\quad.\]
If $Y \subseteq X$ is a subspace, we will write $X/[Y] \defeq X/[\iota(Y)]$, 
where $\iota \from Y \injto X$ is the inclusion.
\end{definition}

The justification for our notation is the following.

\begin{proposition}
For any $[f] \from Y \to X$, the quotient coarse space $X/[f](Y)$ and the 
natural map $X \to X/[f](Y)$ only depend on the image of $[f]$.
\end{proposition}

\begin{proof}
$[f]$ factorizes canonically as
\[
    Y \nameto{\smash{[\beta]}} [f](Y) \nameto{\smash{\im [f]}} X
\]
(Proposition~\ref{prop:Crs-factor-I} and Definition~\ref{def:Crs-image}). Thus 
$X/[f](Y)$ is also the colimit of the diagram
\[\begin{CD}
    Y @>{[\beta]}>> [f](Y) @>{\im [f]}>> X \\
    @V{[\tau_Y]}VV @V{[\tau_{[f](Y)}]}VV \\
    \Terminate(Y) @>{[\Terminate(\beta)]}>> \Terminate([f](Y)),
\end{CD}\]
and hence also of the cofinal subdiagram obtained by deleting $Y$ and 
$\Terminate(Y)$.
\end{proof}

The coarse categories have all push-outs and we have seen how to describe them 
concretely; the standard construction would take $X/[f](Y)$ to be, as a set, 
the disjoint union of $X$, $Y$, and $\Terminate(Y)$. Taking a representative 
coarse map $f \from Y \to X$, we have two ``smaller'' descriptions of the 
quotient:
\begin{enumerate}
\item Take $X/[f](Y) \defeq X \disjtunion \Terminate(Y)$ as a set with the 
    coarse structure generated by $\calE_X$, $\calE_{\Terminate(Y)}$, and 
    $\set{(f \cross \tau_Y)(F) \suchthat F \in \calE_Y}$, where we consider 
    $\Terminate(Y)$ and $X$ as subsets of $X \disjtunion \Terminate(Y)$. (This 
    is a particular instance of a ``smaller'' construction of push-outs in 
    $\CATCrs$.)
\item Take $X/[f](Y) \defeq X$ as a set with the coarse structure generated by 
    $\calE_X$ and $f_* \calE_{\Terminate(Y)}$, where we treat $f$ as a set map 
    $\Terminate(Y) = Y \to X$.
\end{enumerate}

Using the second description above and applying 
Corollary~\ref{cor:Crs-epi-crs-structs} (the left and right supports of 
entourages in $f_* \calE_{\Terminate(Y)}$ are already unital subspaces of $X$), 
we immediately get the following.

\begin{proposition}
For any $[f] \from Y \to X$, $X/[f](Y)$ is a quotient of $X$ in the categorical 
sense (i.e., the natural map $[q] \from X \to X/[f](Y)$ is epi).
\end{proposition}

\subsection{Restricted coarse categories}\label{subsect:rest-Crs}

The lack of restriction on the size of coarse spaces (other than that imposed 
by the choice of universe) may be somewhat bothersome, and moreover prevent 
$\CATCrs$ from having a terminal object. It is tempting to restrict the 
cardinality of coarse spaces, i.e., consider the full subcategory of $\CATCrs$ 
of the coarse spaces of cardinality at most $\kappa$, for some fixed, small 
(probably infinite) cardinal $\kappa$. This is not the correct thing to do: 
First, one would no longer have all small limits and colimits (though as long 
as $\kappa$ is infinite one have all finite limits and colimits). Second, and 
more importantly, it would bar constructions involving the set of (set) 
functions $Y \to X$ ($\#X, \#Y \leq \kappa$) which will be important in 
\cite{crscat-II}.

A better way to proceed is to consider the full subcategory of $\CATCrs$ of 
coarse spaces $X$ for which there exists a coarse map $X \to R$, where $R 
\defeq \Terminate(R_0)$ for some fixed $R_0$. (Of particular interest is the 
case when $R_0$ is a unital coarse space of some infinite cardinality $\kappa$, 
in which case $R = |R_0|_1$ only depends on $\kappa$ up to coarse equivalence.)

We will first discuss this in full generality, using terminology from the 
beginning of \S\ref{subsect:Crs-Term}. In the following, suppose $\calC$ is 
some category and that is some object which $\tilde{X}$ terminates any object 
(e.g., itself) in $\calC$.

\begin{definition}
The \emph{$\tilde{X}$-restriction} $\calC_{\preceq \tilde{X}}$ of $\calC$ is 
the full subcategory of $\calC$ consisting of all the objects $Y$ in $\calC$ 
such that there exists some (unique) arrow $Y \to \tilde{X}$.
\end{definition}

In other words, $\calC_{\preceq \tilde{X}}$ consists of all objects which are 
terminated by $\tilde{X}$. Equivalently, one may consider the comma category 
$(\calC \CATover \tilde{X})$. It is easy to check that the range restricted 
projection functor $(\calC \CATover \tilde{X}) \to \calC_{\preceq \tilde{X}}$ 
is an isomorphism of categories.

Let $I \from \calC_{\preceq \tilde{X}} \to \calC$ denote the inclusion functor. 
When a \emph{nonzero} limit $\calC_{\preceq \tilde{X}}$ already exists in 
$\calC$, the limits are the same. More precisely, we have the following.

\begin{proposition}\label{prop:restricted-lim}
Suppose $\calF \from \calJ \to \calC_{\preceq \tilde{X}}$, where $\calJ$ is 
nonempty. If the limit $\pfx{\calC}\OBJlim (I \circ \calF)$ exists, then
\[
    \pfx{\calC_{\preceq \tilde{X}}}\OBJlim \calF
        = \pfx{\calC}\OBJlim (I \circ \calF);
\]
i.e., the limit of $\calF$ in $\calC_{\preceq \tilde{X}}$ exists and any 
limiting cone in $\calC$ gives a limiting cone in $\calC_{\preceq \tilde{X}}$.
\end{proposition}

\begin{proof}
The nonemptiness of $\calJ$ ensures that the object $\pfx{\calC}\OBJlim (I 
\circ \calF)$ is in $\calC_{\preceq \tilde{X}}$ (since it must map to some 
object of $\calC_{\preceq \tilde{X}}$, hence to $\tilde{X}$). The rest follows 
easily, since the inclusion functor $I$ is fully faithful.
\end{proof}

The following is trivial.

\begin{proposition}
$\tilde{X}$ is a terminal object (i.e., zero limit) in $\calC_{\preceq 
\tilde{X}}$.
\end{proposition}

Thus $\calC_{\preceq \tilde{X}}$ has all the limits that $\calC$ does (to the 
extent that this makes sense), but also has a terminal object, which $\calC$ 
may not have. However, $\calC$ may have a terminal object which is not 
isomorphic to $\tilde{X}$ (in which case $\calC_{\preceq \tilde{X}}$ is a 
proper subcategory of $\calC$), so the inclusion functor $I$ may not preserve 
limits.

The result dual to Proposition~\ref{prop:restricted-colim} is true without the 
nonemptiness criterion.

\begin{proposition}\label{prop:restricted-colim}
Suppose $\calF \from \calJ \to \calC_{\preceq \tilde{X}}$. If the colimit 
$\pfx{\calC}\OBJcolim (I \circ \calF)$ exists, then
\[
    \pfx{\calC_{\preceq \tilde{X}}}\OBJcolim \calF
        = \pfx{\calC}\OBJcolim (I \circ \calF).
\]
\end{proposition}

\begin{proof}
If $\pfx{\calC}\OBJcolim (I \circ \calF)$ exists, then it maps to $\tilde{X}$ 
since there is a (unique) cone $\calJ \to \tilde{X}$; thus the colimiting cone 
is actually in $\calC_{\preceq \tilde{X}}$ and is universal since $I$ is fully 
faithful.
\end{proof}

Now, we return to our coarse context. Suppose $R \defeq \Terminate(R_0)$ for 
some coarse space $R_0$. The \emph{$R$-restricted coarse category} 
$\CATCrs_{\preceq R}$ is, as the notation indicates, the $R$-restriction of 
$\CATCrs$. We similarly get $R$-restricted connected and connected, nonempty 
coarse categories. We refer to the above categories collectively (i.e., for all 
$R$ and the various cases) as the \emph{restricted coarse categories}.

\begin{theorem}
The restricted coarse categories have all (small) limits and colimits.
\end{theorem}

\begin{proof}
This follows immediately from Theorems \ref{thm:Crs-lim} 
and~\ref{thm:Crs-colim}, and Propositions \ref{prop:restricted-lim} 
and~\ref{prop:restricted-colim}.
\end{proof}

One can also check that all the earlier facts on monics and images, epis and 
coimages, quotients, etc. hold in the restricted coarse categories.

%%%%%%%%%%%%%%%%%%%%%%%%%%%%%%%%%%%%%%%%%%%%%%%%%%%%%%%%%%%%%%%%%%%%%%%%%%%%%%%%

\section{Topology and coarse spaces}\label{sect:top-crs}

Our coarse spaces are discrete, as opposed to the more standard definition of 
\emph{proper coarse spaces} which allows coarse spaces to carry topologies and 
thus has different properness requirements (see the works of Roe, e.g., 
\cite{MR2007488}*{Def.~2.22}). Our aim here is not to provide a general 
discussion of topological coarse spaces but to provide a means from going from 
Roe's \emph{proper coarse spaces} to our (discrete) coarse spaces.

We will use the terms \emph{compact} and \emph{locally compact} in the sense of 
Bourbaki \cite{MR1712872}*{Ch.~I \S{}9}, including the Hausdorff condition; in 
fact, all spaces will be Hausdorff unless otherwise stated. Throughout, $X$ and 
$Y$ will denote paracompact, locally compact topological spaces. Recall that a 
subset $K$ of a space $X$ is \emph{relatively compact} if it is contained in 
some compact subspace of $X$. (If $X$ is Hausdorff, $K$ is relatively compact 
if and only if $\overline{K}$ is compact.)

\subsection{Roe coarse spaces}\label{subsect:Roe-crs-sp}

We will diverge from the standard terminology to avoid confusion with our 
previously defined terms. \emph{Roe coarse spaces} will be what are usually 
called proper coarse spaces. Let us recall these definitions (compare 
Definitions \ref{def:prop-ax} and~\ref{def:crs-sp}).

\begin{definition}[see, e.g., \cite{MR2007488}*{Def.~2.1}]%
        \label{def:Roe-prop-ax}
A subset $E \subseteq X^{\cross 2}$ satisfies the \emph{Roe properness axiom} 
if $E \cdot K$ and $K \cdot E$ are relatively compact subsets of $X$ for all 
(relatively) compact $K \subseteq X$.
\end{definition}

\begin{definition}[see, e.g., \cite{MR2007488}*{Def.~2.22}]%
        \label{def:Roe-crs-sp}
A \emph{Roe coarse structure} on $X$ is a subset $\calR_X \subseteq
\powerset(X^{\cross 2})$ such that:
\begin{enumerate}
\item each $E \in \calR_X$ satisfies the Roe properness axiom;
\item $\calR_X$ is closed under the operations of addition, multiplication, 
    transpose, and the taking of subsets;
\item\label{def:Roe-crs-sp:III} if $K \subseteq X$ is \emph{bounded} in the 
    sense that $K^{\cross 2} \in \calR_X$, then $K$ is relatively compact; and
\item\label{def:Roe-crs-sp:IV} there is a neighbourhood (with respect to the 
    product topology on $X^{\cross 2}$ of the unit (i.e., diagonal) $1_X$ which 
    is in $\calR_X$.
\end{enumerate}
A \emph{Roe coarse space} is a paracompact, locally compact space $X$ equipped 
with a Roe coarse structure $\calR_X$ on $X$.
\end{definition}

\enumref{def:Roe-crs-sp:IV} implies Roe coarse spaces are always unital (in the 
obvious sense; see Definition~\ref{def:uni-conn}) and that any Roe coarse space 
$X$ has an open cover $\calU \subseteq \powerset(X)$ which is \emph{uniformly 
bounded} in the sense that $\bigunion_{U \in \calU} U^{\cross 2}$ is in 
$\calR_X$. Paracompactness implies that this cover can be taken to be locally 
finite. The local compactness requirement is redundant, since it is implied by 
\enumref{def:Roe-crs-sp:III} and \enumref{def:Roe-crs-sp:IV}.

\begin{definition}\label{def:top-prop}
A continuous map $f \from Y \to X$ between locally compact spaces is 
\emph{topologically proper} if $f^{-1}(K)$ is compact for every compact $K 
\subseteq X$. More generally, also say that a (not necessarily continuous) map 
$f \from Y \to X$ between locally compact spaces is \emph{topologically proper} 
if $f^{-1}(K)$ is relatively compact for every relatively compact $K \subseteq 
X$.
\end{definition}

\begin{definition}[see, e.g., \cite{MR2007488}*{Def. 2.21 and~2.14}]
A (not necessarily continuous) map $f \from Y \to X$ between Roe coarse spaces 
is a \emph{Roe coarse map} if it is topologically proper and \emph{preserves 
entourages} in the sense that $f^{\cross 2}(F) \in \calR_X$ for all $F \in 
\calR_Y$. (Roe coarse maps are usually called \emph{proper coarse maps}.) Roe 
coarse maps $f, f' \from Y \to X$ are \emph{close} if $(f \cross f')(1_Y) \in 
\calR_X$ (or equivalently if $(f \cross f')(F) \in \calR_X$ for all $F \in 
\calR_Y$).
\end{definition}

We get an obvious \emph{Roe precoarse category} $\CATRoePCrs$ with objects all 
(small) Roe coarse spaces and arrows Roe coarse maps, and a quotient \emph{Roe 
coarse category} $\CATRoeCrs$ with the same objects but whose arrows are 
closeness classes of Roe coarse maps. \emph{Roe coarse equivalences} are Roe 
coarse maps which represent isomorphisms in $\CATRoeCrs$.

\subsection{Discretization of Roe coarse spaces}\label{subsect:Disc}

We now provide a way of passing from Roe coarse spaces to our (discrete) coarse 
spaces.

\begin{definition}
A set $E \in \powerset(X^{\cross 2})$ satisfies the \emph{topological 
properness axiom} (with respect to the topology of $X$) if, for all compact 
subspaces $K \subseteq X$, $(\pi_1 |_E)^{-1}(K)$ and $(\pi_2 |_E)^{-1}(K)$ are 
finite.
\end{definition}

Since all our spaces are Hausdorff hence $\text{T}_{\text{1}}$, the topological 
properness axiom implies the (discrete) properness axiom 
(Definition~\ref{def:prop-ax}).

The following is easy to check.

\begin{proposition}
A set $E \in \powerset(X^{\cross 2})$ satisfies the topological properness 
axiom if and only if $E$ is a (closed) discrete subset of $X^{\cross 2}$ and 
the restricted projections $\pi_1 |_E, \pi_2 |_E \from E \to X$ are 
topologically proper maps.
\end{proposition}

\begin{remark}\label{rmk:top-crs-sp}
We provide only a means from passing from Roe coarse spaces to our coarse 
spaces and not a complete discussion of ``topological coarse spaces'' since the 
topological properness axiom does not encompass the axioms of 
Definition~\ref{def:Roe-crs-sp} (\enumref{def:Roe-crs-sp:IV} in particular). We 
would like not just a direct translation of Roe's definition to our setting, 
but a proper generalization: First, we would like to allow nonunital 
topological coarse spaces. Second, we do not want to impose local compactness 
for two (possibly related) reasons:
\begin{inparaenum}
\item The ``topological coarse category'' should have all colimits (including 
    infinite ones). In particular, we are interested in ``large'' simplicial 
    complexes which may not be locally finite.
\item We wish to be able to analyze Hilbert space and other Banach spaces 
    directly as coarse spaces. This seems especially relevant as methods 
    involving uniform (i.e., coarse) embeddings into such spaces have gained 
    prominence in recent years (e.g., in \cite{MR1728880}, Yu shows that the 
    Coarse Baum--Connes Conjecture is true for metric spaces of bounded 
    geometry which uniformly embed in Hilbert space).
\end{inparaenum}

Instead of requiring that spaces be paracompact and locally compact, we 
should probably require that spaces be \emph{compactly generated} (i.e., be 
weak Hausdorff $k$-spaces). The topological properness axiom makes sense for 
such spaces (weak Hausdorffness still implies the $\text{T}_{\text{1}}$ 
condition), but the problem of translating axioms \enumref{def:Roe-crs-sp:III} 
and \enumref{def:Roe-crs-sp:IV} becomes more complicated. Moreover, in the 
compactly generated case, there are different, inequivalent definitions for 
``topological properness'' (whereas they all agree in the locally compact case; 
see, e.g., \cite{MR1712872}*{Ch.~I \S{}10}), though perhaps one could still use 
Definition~\ref{def:top-prop} verbatim. In that case, the above Proposition 
remains true so long as $X^{\cross 2}$ is given the categorically appropriate 
topology, namely the $k$-ification of the standard product topology. We leave 
these problems to a future paper \cite{crscat-III}.
\end{remark}

Compare the following, which is easy, with Proposition~\ref{prop:prop-ax-alg}.

\begin{proposition}
If $E, E' \in \powerset(X^{\cross 2})$ satisfy the topological properness 
axiom, then $E + E'$, $E \circ E'$, $E^\transpose$, and all subsets of $E$ 
satisfy the topological properness axiom. Also, all singletons $\set{e}$, $e 
\in X^{\cross 2}$, and hence all finite subsets of $X^{\cross 2}$ satisfy the 
properness axiom. Consequently,
\[
    \calE_{|X|_\tau}
        \defeq \set{E \in \calE_{|X|_1} \subseteq \powerset(X^{\cross 2})
            \suchthat \text{$E$ satisfies the topological properness axiom}}
\]
is a coarse structure on the set $X$ (in the sense of 
Definition~\ref{def:crs-sp}).
\end{proposition}

\begin{definition}\label{def:Disc-Roe}
The \emph{discretization} of a Roe coarse space $X$ is the coarse space 
$\Disc(X) \defeq X$ as a set with the coarse structure
\[
    \calE_{\Disc(X)} \defeq \calR_X \intersect \calE_{|X|_\tau}
\]
(consisting of all elements of $\calR_X$ which satisfy the topological 
properness axiom).
\end{definition}

It is easy to check that $\calE_{\Disc(X)}$ is in fact a coarse structure on 
the set $X$. Unless $X$ is discrete, the coarse space $\Disc(X)$ is not unital, 
even though the Roe coarse space $X$ is.

\begin{proposition}
If $f \from Y \to X$ is a Roe coarse map, then the set map $\Disc(f) \defeq f$ 
is coarse as a map $\Disc(Y) \to \Disc(X)$.
\end{proposition}

\begin{proof}
The only thing to check is that if $f$ (not necessarily continuous) is 
topologically proper and $F \subseteq Y^{\cross 2}$ satisfies the topological 
properness axiom, then $E \defeq f^{\cross 2}(F) \subseteq X^{\cross 2}$ also 
satisfies the topological properness axiom. This follows since
\[
    E \cdot K \subseteq f(F \cdot f^{-1}(K))
\]
and $f$ is topologically proper (and similarly symmetrically).
\end{proof}

Since, trivially, $\Disc(f \circ g) = \Disc(f) \circ \Disc(g)$, we get the 
following.

\begin{corollary}
$\Disc$ is a functor from the Roe precoarse category $\CATRoePCrs$ to the 
precoarse category $\CATPCrs$.
\end{corollary}

$\Disc$ is coarsely invariant in the following way, which yields a canonical 
functor $[\Disc] \from \CATRoeCrs \to \CATCrs$ between the closeness quotients. 
(We continue to write $\Disc(X)$ instead of $[\Disc](X)$ for Roe coarse 
spaces.)

\begin{proposition}
If Roe coarse maps $f, f' \from Y \to X$ are close, then
\[
    \Disc(f), \Disc(f') \from \Disc(Y) \to \Disc(X)
\]
are close coarse maps.
\end{proposition}

\begin{proof}
The result follows easily from the following fact (which is also easy): If $f, 
f' \from Y \to X$ are topologically proper and $F \subseteq Y^{\cross 2}$ 
satisfies the topological properness axiom, then $(f \cross f')(F) \subseteq 
X^{\cross 2}$ also satisfies the topological properness axiom.
\end{proof}

\begin{corollary}
If $f \from Y \to X$ is a Roe coarse equivalence, then $\Disc(f) \from \Disc(Y) 
\to \Disc(X)$ is a coarse equivalence.
\end{corollary}

\subsection{Properties of the discretization functors}

Let $\CATDiscRoePCrs \subseteq \CATRoePCrs$ and $\CATDiscRoeCrs \subseteq 
\CATRoeCrs$ be the full subcategories of \emph{discrete} Roe coarse spaces 
(call them the \emph{discrete Roe precoarse} and \emph{coarse categories}, 
respectively). On the discrete subcategories, $\Disc$ and $[\Disc]$ are fully 
faithful.

\begin{proposition}\label{prop:DiscRoePCrs-Disc-fullfaith}
If $X$, $Y$ are Roe coarse spaces with $Y$ discrete, then the map
\[
    \Disc_{Y,X} \from \Hom_{\CATRoePCrs}(Y,X)
        \to \Hom_{\CATPCrs}(\Disc(Y),\Disc(X))
\]
is a bijection. Hence, in particular, the restriction of $\Disc$ to 
$\CATDiscRoePCrs$ (which actually maps into $\CATUniPCrs$) is a fully faithful 
functor.
\end{proposition}

\begin{proof}
$\Disc_{Y,X}$ is trivially injective, so it only remains to show surjectivity. 
Suppose $f \from \Disc(Y) \to \Disc(X)$ is a coarse map. If $K \subseteq X$ is 
relatively compact, then
\[
    f^{-1}(K) = f^{-1}(f^{\cross 2}(1_Y) \cdot K)
\]
is finite: since $Y$ is discrete, $\Disc(Y)$ is unital so $f^{\cross 2}(1_Y) 
\in \calE_{\Disc(X)}$ satisfies the topological properness axiom (so $f^{\cross 
2}(1_Y) \cdot K$ is finite) and $f$ is (discretely) globally proper. Thus $f$ 
is topologically proper. Since $Y$ is discrete, $\calE_{\Disc(Y)} = \calR_Y$, 
so $f$ preserves entourages of $\calR_Y$ (of course, $\calE_{\Disc(X)} 
\subseteq \calR_X$). Thus $f$ is Roe coarse as a map $Y \to X$.
\end{proof}

The unrestricted functor $\Disc \from \CATRoePCrs \to \CATPCrs$ is \emph{not} 
full.

\begin{example}\label{ex:RoePCrs-Disc-notfull}
Let $X \defeq \setRplus$ equipped with the Euclidean metric Roe coarse 
structure (see \S\ref{subsect:prop-met}), and $Y \defeq \setRplus \union 
\set{\infty}$ be the one-point compactification of $\setRplus$ equipped with 
the unique Roe coarse structure $\calR_Y \defeq \powerset(Y^{\cross 2})$ (which 
is also the metric Roe coarse structure for any metric which metrizes $Y$ 
topologically). Define $f \from Y \to X$ by
\[
    f(t) \defeq \begin{cases}
            t & \text{if $t \in \setRplus$, and} \\
            0 & \text{if $t = \infty$.}
        \end{cases}
\]
Then $f$ is actually coarse as a map $\Disc(Y) \to \Disc(X)$. However, clearly 
$f$ does not preserve entourages of $\calR_Y$, hence does \emph{not} define a 
Roe coarse map $Y \to X$. As a map $\Disc(Y) \to \Disc(X)$, $f$ is close to any 
constant map $\Disc(Y) \to \Disc(X)$ (sending all of $Y$ to some fixed element 
of $X$); every such constant map \emph{does} define a Roe coarse map $Y \to X$.
\end{example}

\begin{proposition}\label{prop:DiscRoeCrs-Disc-fullfaith}
If $X$, $Y$ are Roe coarse spaces with $Y$ discrete, then the map
\[
    [\Disc]_{Y,X} \from \Hom_{\CATRoeCrs}(Y,X)
        \to \Hom_{\CATCrs}(\Disc(Y),\Disc(X))
\]
is a bijection. Hence the restriction of $[\Disc]$ to $\CATDiscRoeCrs$ (which 
actually maps into $\CATUniCrs$) is fully faithful.
\end{proposition}

\begin{proof}
By the previous Proposition, $[\Disc]_{Y,X}$ is surjective, so it only remains 
to show injectivity. Suppose $f, f' \from Y \to X$ are Roe coarse maps. If 
$\Disc(f)$ is close to $\Disc(f')$, then since $\Disc(Y)$ is unital,
\[
    (f \cross f')(1_Y) = (\Disc(f) \cross \Disc(f'))(1_Y)
            \in \calE_{\Disc(X)} \subseteq \calR_X,
\]
so $f$ is close to $f'$, as required.
\end{proof}

If $X' \subseteq X$ is a \emph{closed} subspace of a Roe coarse space, then the 
obvious \emph{Roe subspace coarse structure} $\calR_{X'} \defeq \calR_X |_{X'} 
\defeq \calR_X \intersect \powerset((X')^{\cross 2})$ is actually Roe coarse 
structure on $X'$ (this is not the case if $X'$ is not closed), which makes 
$X'$ into a \emph{Roe coarse subspace} of $X$. The inclusion of any Roe coarse 
subspace into the ambient space is a Roe coarse map. The following result is 
well known.

\begin{proposition}\label{prop:Roe-disc-subsp}
For any Roe coarse space $X$, there is a (closed) discrete Roe coarse subspace 
$X' \subseteq X$ such that the inclusion $\iota \from X' \to X$ is a Roe coarse 
equivalence.
\end{proposition}

\begin{proof}
Fix a locally finite, uniformly bounded open cover $\calU$ of $X$ by nonempty 
sets. For each $U \in \calU$, pick a point $x'_U \in U$ and put $X' \defeq 
\set{x'_U \suchthat U \in \calU}$. Since $\calU$ is locally finite, it is easy 
to check that $X'$ is closed and discrete.

Invoking the Axiom of Choice, fix a map $\kappa \from X \to X'$ such that, for  
all $x \in X$, $\kappa(x) \in U$ for some $U \in \calU$ such that $x \in U$. We 
may also ensure that $\kappa(x') = x'$ for all $x' \in X'$. $\kappa$ is 
topologically proper: For any $x' \in X'$,
\[
    \kappa^{-1}(\set{x'})
        \subseteq \bigunion_{\substack{U \in \calU \suchthat \\ x' \in U}} U
\]
which is a finite union of relatively compact sets, hence 
$\kappa^{-1}(\set{x'})$ is relatively compact (this suffices to show 
topological properness since $X'$ is discrete). $\kappa$ preserves entourages 
of $X$: Put
\[
    E_\calU \defeq \bigunion_{U \in \calU} U^{\cross 2} \in \calR_X;
\]
for any $E \in \calR_X$,
\[
    \kappa^{\cross 2}(E) \subseteq E_\calU \circ E \circ E_\calU \in \calR_X,
\]
hence $\kappa^{\cross 2}(E) \in \calR_X |_{X'}$, as required. Thus $\kappa$ is 
a Roe coarse map.

Trivially, $\kappa \circ \iota = \id_{X'}$. Finally, $\iota \circ \kappa$ is 
close to $\id_X$: Letting $E_\calU$ be as above, we have
\[
    (\kappa \cross \id_X)(1_X) \subseteq E_\calU \in \calR_X,
\]
as required.
\end{proof}

\begin{remark}
Though we do not so insist, Roe coarse maps are sometimes required to be Borel 
(see, e.g., \cite{MR1451755}*{Def.~2.2}). In that case, the map $\kappa$ used 
in the above proof may not suffice. However, if one insists that all Roe coarse 
spaces be, e.g., second countable, then one can construct $\kappa$ to be Borel. 
Thus, as long as one so constrains the allowable Roe coarse spaces, the above 
Proposition remains true.
\end{remark}

\begin{corollary}
The inclusion functor $\CATDiscRoeCrs \injto \CATRoeCrs$ is fully faithful and 
in fact an equivalence of categories.
\end{corollary}

\begin{theorem}\label{thm:RoeCrs-Disc-fullfaith}
The functor $[\Disc] \from \CATRoeCrs \to \CATCrs$ is fully faithful.
\end{theorem}

\begin{proof}
This is immediate upon combining Propositions 
\ref{prop:DiscRoeCrs-Disc-fullfaith} and~\ref{prop:Roe-disc-subsp}.
\end{proof}

Every unital coarse space (in our sense) becomes a Roe coarse space when it is 
given the discrete topology, with coarse maps between unital coarse spaces 
becoming Roe coarse maps. Thus $\CATUniPCrs$ and $\CATDiscRoePCrs$ are 
isomorphic as categories, and hence so too are $\CATUniCrs$ and 
$\CATDiscRoeCrs$.

\begin{corollary}\label{cor:UniCrs-RoeCrs-equiv}
Our unital coarse category $\CATUniCrs$ is equivalent to the Roe coarse 
category $\CATRoeCrs$, with the functor which sends a unital coarse space to 
the ``identical'' discrete Roe coarse space an equivalence of categories.
\end{corollary}

%%%%%%%%%%%%%%%%%%%%%%%%%%%%%%%%%%%%%%%%%%%%%%%%%%%%%%%%%%%%%%%%%%%%%%%%%%%%%%%%

\section{Examples and applications}\label{sect:ex-appl}

As stated in the Introduction, we will not discuss even the standard 
applications of coarse geometry. We will first discuss a couple of basic 
examples which we will need later, namely proper metric spaces and continuous 
control, and then briefly examine a few things which arise from the categorical 
point of view (some of which are not obviously possible in standard, unital 
coarse geometry).

\subsection{Proper metric spaces}\label{subsect:prop-met}

Suppose that $(X,d) \defeq (X,d_X)$ is a proper metric space (i.e., its closed 
balls are compact). We wish to produce a coarse space from $X$; we have already 
discussed the discrete case in Example~\ref{ex:disc-met}, and what follows is a 
generalization of that.

There is a well known way to produce a Roe coarse space $|X|_d^\TXTRoe$ from 
$(X,d)$ (noting that properness implies local compactness, and metrizability 
implies paracompactness), taking the Roe coarse structure to be consist of the 
$E \subseteq X^{\cross 2}$ satisfying inequality \eqref{ex:disc-met:eq} of 
Example~\ref{ex:disc-met} (see, e.g., \cite{MR2007488}*{Ex.~2.5}). One can then 
apply the discretization functor $\Disc$ to this Roe coarse space to obtain the 
\emph{($d$-)metric coarse space} $|X| \defeq |X|_d$. More directly and entirely 
equivalently, $|X|_d$ has as entourages the $E \in \calE_{|X|_\tau}$ (i.e., the 
$E$ satisfying the topological properness axiom) which also satisfy the same 
inequality \eqref{ex:disc-met:eq}. As in the discrete case, we may also allow 
$d(x,x') = \infty$ (for $x \neq x'$), and $|X|_d$ is connected if and only if 
$d(x,x') < \infty$ for all $x$ and $x'$. If $X' \subseteq X$ is a closed 
(topological) subspace, then the restriction of $d$ to $X'$ makes $X'$ into a 
proper metric space; the subspace coarse structure on $X'$ is the same as the 
coarse structure coming from the restricted metric.

Suppose $(Y,d_Y)$ is another proper metric space. A (not necessarily 
continuous) map $f \from Y \to X$ is Roe coarse as a map $|Y|_{d_Y}^\TXTRoe \to 
|X|_{d_X}^\TXTRoe$ if and only if it is topologically proper and
\begin{equation}\label{subsect:prop-met:Roe-crs:eq}
    \sup \set{d_X(f(y),f(y'))
            \suchthat \text{$y, y' \in Y$ and $d_Y(y,y') \leq r$}} < \infty
\end{equation}
for every $r \geq 0$. Since $X$, $Y$ are proper metric spaces, $f$ is 
topologically proper if and only if it is \emph{metrically proper} in the sense 
that inverse images of metrically bounded subsets of $X$ are metrically bounded 
in $Y$. Roe coarse maps $f, f' \from |Y|_{d_Y}^\TXTRoe \to |X|_{d_X}^\TXTRoe$ 
are close if and only if
\begin{equation}\label{subsect:prop-met:Roe-close:eq}
    \sup \set{d_X(f(y),f'(y)) \suchthat y \in Y} < \infty.
\end{equation}

We must warn that there may be a map $f \from Y \to X$ which is coarse (in our 
sense) as a map $|Y|_{d_Y} \to |X|_{d_X}$, yet does not satisfy 
\eqref{subsect:prop-met:Roe-crs:eq}. Similarly, there may be coarse maps $f, f' 
\from |Y|_{d_Y} \to |X|_{d_X}$ which are close but do not satisfy 
\eqref{subsect:prop-met:Roe-close:eq}. Example~\ref{ex:RoePCrs-Disc-notfull}, 
which shows that $\Disc$ is not full, exhibits both phenomena. In the former 
case, Theorem~\ref{thm:RoeCrs-Disc-fullfaith} shows that every coarse map $f' 
\from |Y|_{d_Y} \to |X|_{d_X}$ is close to some coarse map $f \from |Y|_{d_Y} 
\to |X|_{d_X}$ which satisfies \ref{subsect:prop-met:Roe-crs:eq}. (The 
corresponding statement in the latter case is trivial.) Alternatively, one may 
avoid both ``problems'' by considering only discrete, proper metric spaces 
(Proposition~\ref{prop:DiscRoePCrs-Disc-fullfaith}); every proper metric space 
is Roe coarsely equivalent to a discrete one 
(Proposition~\ref{prop:Roe-disc-subsp}).

\begin{remark}[see, e.g., \cite{MR2007488}*{\S{}1.3}]\label{rmk:lsLip-qisom}
If $X$ and $Y$ are proper \emph{length spaces}, then one can characterize the 
Roe coarse maps, and indeed Roe coarse equivalences, $Y \to X$ a bit more 
strictly: A map $f \from Y \to X$ (not necessarily continuous) is Roe coarse if 
and only if it is (metrically/topologically) proper and \emph{large-scale 
Lipschitz} in the sense that there exist constants $C > 0$ and $R \geq 0$ such 
that
\[
    d_X(f(y),f(y')) \leq C d_Y(y,y') + R
\]
for all $y, y' \in Y$. $f$ is a Roe coarse equivalence if and only if it is a 
\emph{quasi-isometry} in that there are constants $c, C > 0$ and $r, R \geq 0$ 
such that
\[
    c d_Y(y,y') - r \leq d_X(f(y),f(y')) \leq C d_Y(y,y') + R
\]
for all $y, y' \in Y$ (evidently, one can always take $c = 1/C$ and $r = R$, as 
is conventional) and there is a constant $D \geq 0$ such that every point of 
$X$ is within distance $D$ of a point in the image of $f$.

One can replace the length space hypothesis with a weaker condition, but some 
hypothesis is necessary; for general metric spaces there are Roe coarse maps, 
and indeed Roe coarse equivalences, which are not large-scale Lipschitz. 
However, every proper large-scale Lipschitz map is evidently also Roe coarse, 
and every quasi-isometry is a coarse equivalence.
\end{remark}

\subsection{Continuous control}\label{subsect:cts-ctl}

Most of the following originates from \cites{MR1277522, MR1451755}, but see 
also, e.g., \cite{MR2007488}*{\S{}2.2}. In the following, all topological 
spaces will be assumed to be second countable and locally compact (and 
Hausdorff), whence paracompact. $X$ and $Y$ will always denote such spaces.

\begin{definition}
A \emph{compactified space} is a (second countable, locally compact) 
topological space $X$ equipped with a (second countable) compactification 
$\overline{X}$; its \emph{boundary} is the space $\die X \defeq \overline{X} 
\setminus X$.
\end{definition}

The \emph{continuously controlled Roe coarse structure} $\calR_{|X|_{\die 
X}^\TXTRoe}$ on $X$ (for the compactification $\overline{X}$, or for the 
boundary $\die X$) consists of the $E \subseteq X^{\cross 2}$ such that
\begin{equation}\label{subsect:cts-ctl:eq}
    \overline{E} \subseteq X^{\cross 2} \union 1_{\die X}
        \subseteq \overline{X}^{\cross 2},
\end{equation}
where $1_{\die X}$ is the diagonal subset of $(\die X)^{\cross 2}$ and the 
closure is taken in $\overline{X}$ (for the proof that this a Roe coarse 
structure, see, e.g., \cite{MR2007488}*{Thm.~2.27}). The associated coarse 
space (resulting from applying the discretization functor $\Disc$ to the above 
Roe coarse space $|X|_{\die X}^\TXTRoe$) is the \emph{continuously controlled 
coarse space} $|X|_{\die X}$ (for the compactification $\overline{X}$, or for 
the boundary $\die X$) whose entourages are the $E \in \calE_{|X|_\tau}$ (i.e., 
$E$ satisfying the topological properness axiom) which also satisfy 
\eqref{subsect:cts-ctl:eq}.

\begin{remark}
If $X$ is compact (so $\overline{X} = X$ and $\die X = \emptyset$), then 
$|X|_{\die X} = |X|_0^\TXTconn$ (i.e., $X$ equipped with the initial connected 
coarse structure).
\end{remark}

The following is standard.

%FIXME: where?
\begin{proposition}
Suppose $X$, $Y$ are compactified spaces. Any Roe coarse map $f \from |Y|_{\die 
Y}^\TXTRoe \to |X|_{\die X}^\TXTRoe$ determines a canonical \emph{continuous} 
map $\die Y \to \die X$ which we denote by $\die [f]$. Moreover, Roe coarse 
maps $f, f' \from |Y|_{\die Y}^\TXTRoe \to |X|_{\die X}^\TXTRoe$ are close if 
and only if $\die [f] = \die [f']$ (which justifies our notation).
\end{proposition}

The ``converse'' is also true: Any set map $Y \to X$ (not necessarily 
continuous, but necessarily topologically proper) which ``extends 
continuously'' to a continuous map $\die Y \to \die X$ is Roe coarse as a map 
$|Y|_{\die Y}^\TXTRoe \to |X|_{\die X}^\TXTRoe$. This is essentially 
tautological, since the definition of ``extends continuously'' is exactly the 
definition of ``is continuously controlled''.

\begin{proof}
Fix a Roe coarse map $f$. Given $y_\infty \in \die Y$, define $(\die 
[f])(y_\infty)$ as follows: By second countability, there is a sequence 
$\seq{y_n}_{n=1}^\infty$ in $Y$ which converges to $y_\infty$. Then the 
diagonal set $1_{\set{y_n \suchthat n \in \setN}}$ is in $\calR_{|Y|_{\die 
Y}^\TXTRoe}$, so $1_{\set{f(y_n) \suchthat n \in \setN}}$ must be in 
$\calR_{|X|_{\die X}^\TXTRoe}$. By topological properness, the limit points of 
$\seq{f(y_n)}_{n=1}^\infty$ in $\overline{X}$ (which exist by compactness) are 
all in $\die X \subseteq \overline{X}$; in fact there is only one limit point 
which we call $(\die [f])(y_\infty)$. Well-definedness follows from the 
observation that if $\seq{y'_n}_{n=1}^\infty \subseteq Y$ (possibly a 
subsequence of $\seq{y_n}_{n=1}^\infty$) also converges to $y_\infty$, then 
$1_{\set{(y_n,y'_n) \suchthat n \in \setN}} \in \calR_{|Y|_{\die Y}^\TXTRoe}$ 
hence $1_{\set{(f(y_n),f(y'_n)) \suchthat n \in \setN}} \in \calR_{|X|_{\die 
X}^\TXTRoe}$, so $\seq{f(y'_n)}_{n=1}^\infty$ and $\seq{f(y_n)}_{n=1}^\infty$ 
have the same limit points. To see that $\die [f]$ is continuous, one proves 
sequential continuity (which suffices) using the obvious diagonal argument.

$f$ (and similarly $f'$) ``extend continuously'' to maps $\overline{Y} \to 
\overline{X}$: e.g.,
\[
    \bar{f}(y) \defeq \begin{cases}
            f(y) & \text{if $y \in Y$, and} \\
            (\die [f])(y) & \text{if $y \in \die Y$.}
        \end{cases}
\]
The second assertion then follows using the observation that, for any $F \in 
\calR_{|Y|_{\die Y}^\TXTRoe}$,
\[
    \overline{(f \cross f')(F)}
        = (\bar{f} \cross \bar{f}')(\overline{F})
\]
(closures $\overline{X}^{\cross 2}$ and $\overline{Y}^{\cross 2}$).
\end{proof}

Temporarily let $\calC$ be the category of second countable, compact spaces 
(and continuous maps). If $M \in \Obj(\calC)$, $\setRplus \cross M$ 
compactified with boundary $M$ (so $\overline{\setRplus \cross M}$ is 
homeomorphic to $\ccitvl{0,1} \cross M$) is a compactified space. Then $M 
\mapsto |\setRplus \cross M|_M^\TXTRoe$ (on objects; $g \mapsto \id_{\setRplus} 
\cross g$ on functions) defines a (Roe) coarsely invariant functor 
$\calO_\TXTtop^\TXTRoe \from \calC \to \CATRoePCrs$. By the above Proposition,
\[
    [\calO_\TXTtop^\TXTRoe] \defeq \Quotient \circ \calO_\TXTtop^\TXTRoe
            \from \calC \to \CATRoeCrs
\]
is fully faithful. As $[\Disc] \from \CATRoeCrs \to \CATCrs$ is also fully 
faithful (Proposition~\ref{thm:RoeCrs-Disc-fullfaith}), the resulting 
composition
\[
    [\calO_\TXTtop] \defeq [\Disc] \circ [\calO_\TXTtop^\TXTRoe]
            \from \calC \to \CATCrs
\]
is again fully faithful. Note that
\[
    [\calO_\TXTtop] = \Quotient \circ \calO_\TXTtop,
\]
where $\calO_\TXTtop \defeq \calO_\TXTtop^\TXTRoe \from \calC \to \CATRoePCrs$ 
(a coarsely invariant functor).

\begin{definition}[see, e.g., \cite{MR1817560}*{\S{}6.2}]%
        \label{def:ctsctl-cone}
For any second countable, compact space $M$, the \emph{continuously controlled 
open cone} on $M$ is the coarse space
\[
    \calO_\TXTtop M \defeq |\setRplus \cross M|_M.
\]
\end{definition}

We saw above that $M \mapsto \calO_\TXTtop M$ is a coarsely invariant functor 
from the category of second countable, compact topological spaces to the 
precoarse category $\CATPCrs$.

\begin{remark}[compare \cite{MR1341817}*{Thm.~1.23 and Cor.~1.24}]%
        \label{rmk:ctsctl-cones}
All continuously controlled coarse spaces can be described as cones in a 
natural way. That is, for any compactified space $X$, there is a natural coarse 
equivalence
\[
    \calO_\TXTtop(\die X) \isoto |X|_{\die X}
\]
(indeed, there is a natural Roe coarse equivalence $|\setRplus \cross \die 
X|_{\die X}^\TXTRoe \isoto |X|_{\die X}^\TXTRoe$). Thus, up to coarse 
equivalence, $|X|_{\die X}$ only depends on the topology of the boundary $\die 
X$, and not of $X$ itself. We leave this to the reader.
\end{remark}

\begin{remark}\label{rmk:ctsctl-quot}
Suppose $M$ is a second countable, compact space, and $N \subseteq M$ is a 
closed subspace. There is a natural (coarse) inclusion $\iota \from 
\calO_\TXTtop N \injto \calO_\TXTtop M$ of continuously controlled open cones, 
hence a quotient coarse space
\[
    (\calO_\TXTtop M)/[\calO_\TXTtop N]
        \defeq (\calO_\TXTtop M)/[\iota](\calO_\TXTtop N)
\]
(see \S\ref{subsect:Crs-quot}). One can check that the quotient $(\calO_\TXTtop 
M)/[\calO_\TXTtop N]$ is naturally coarsely equivalent to the continuously 
controlled open cone $\calO_\TXTtop (M/N)$ on the topological quotient $M/N$.
\end{remark}

\begin{remark}\label{rmk:ctsctl-ray}
The continuously controlled ray
\[
    |\coitvl{0,1}|_{\set{1}} \cong |\setRplus|_{\ast} \cong |\setZplus|_{\ast} 
        \cong \calO_\TXTtop \ast
\]
(where $\ast$ is a one-point space) is coarsely equivalent to $|\setZplus|_1$, 
i.e., a countable set with the terminal coarse structure.
\end{remark}

\subsection{Metric coarse simplices}\label{subsect:met-simpl}

We index our simplices in the same way as Mac~Lane \cite{MR1712872}*{Ch.~VII 
\S{}5}, shifted by $1$ from most topologists' indexing. That is, our 
$n$-simplices are topologists' $(n-1)$-simplices (which have geometric 
dimension $n-1$) and we include the ``true'' $0$-simplex.

\begin{definition}
As sets, put $\Delta_0 \defeq \set{0}$, $\Delta_1 \defeq \setRplus \defeq 
\coitvl{0,\infty}$, \ldots, $\Delta_n \defeq (\setRplus)^n$, \ldots. For each 
$n = 0, 1, 2, \dotsc$, let $d \defeq d_n$ be the $l^1$ metric on $\Delta_n$, 
i.e.,
\[
    d_n((x_0,\dotsc,x_{n-1}),(x'_0,\dotsc,x'_{n-1}))
        \defeq |x_0 - x'_0| + \dotsb + |x_{n-1} - x'_{n-1}|,
\]
and denote the resulting coarse space, called the \emph{metric coarse 
$n$-simplex}, by
\[
    |\Delta_n| \defeq |\Delta_n|_\TXTmet \defeq |\Delta_n|_{d_n}
\]
(the metric coarse space defined in \S\ref{subsect:prop-met}). We may also 
substitute the coarsely equivalent unital subspaces $(\setZplus)^n \subseteq 
(\setRplus)^n$ for the $\Delta_n$ when convenient.
\end{definition}

Note that we may replace the $l^1$ metric with any $l^p$-metric ($1 \leq p \leq 
\infty$), since
\[
    \norm{x}_\infty \leq \norm{x}_p \leq \norm{x}_1 \leq n \norm{x}_\infty
\]
(for all $1 \leq p \leq \infty$, $x \in \Delta_n \subseteq \setR^n$); all these 
metrics yield the same Roe coarse structure and hence the same coarse structure 
on $\Delta_n$. See Proposition~\ref{prop:met-simp-univ} below for a bit more 
about the ``universality'' of metric coarse simplices.

For each $n = 0, 1, 2, \dotsc$, $j = 0, \dotsc, n$, define a coarse map 
$\delta_j \defeq \delta_j^n \from |\Delta_n| \to |\Delta_{n+1}|$ by
\begin{equation}\label{subsect:met-simpl:eq:delta}
    \delta_j(x_0, \dotsc, x_{n-1})
        \defeq (x_0, \dotsc, x_{j-1}, 0, x_j, \dotsc, x_{n-1})
\end{equation}
(for $n = 0$, let $\delta_0^0$ be the inclusion). For each $n = 1, 2, 3, 
\dotsc$, $j = 0, \dotsc, n-1$, define a coarse map $\sigma_j \defeq \sigma_j^n 
\from |\Delta_{n+1}| \to |\Delta_n|$ by
\begin{equation}\label{subsect:met-simpl:eq:sigma}
    \sigma_j(x_0, \dotsc, x_n)
        \defeq (x_0, \dotsc, x_{j-1}, x_j+x_{j+1}, x_{j+2}, \dotsc, x_n).
\end{equation}
It is easy to verify that the above maps are coarse and satisfy the 
\emph{cosimplicial identities} (see, e.g., \cite{MR1711612}*{I.1} or equations 
(11)--(13) in \cite{MR1712872}*{Ch.~VII \S{}5}). Consequently, we get a functor 
from the \emph{simplicial category} $\CATSimp$ to $\CATPCrs$. Composing with 
the quotient functor yields the \emph{metric coarse simplex functor}
\[
    |\Delta|_\TXTmet \from \CATSimp \to \CATCrs;
\]
for $n \in \Obj(\CATSimp) = \set{0, 1, 2, \dotsc}$, $|\Delta|_\TXTmet(n) = 
|\Delta_n|_\TXTmet$.

Proceeding as standard (see, e.g., \cite{MR1711612}), we may obtain 
\emph{metric coarse realizations} of any simplicial set (since $\CATCrs$ has 
all colimits), get a corresponding notion of (metric coarse) ``weak 
equivalence'', define \emph{metric coarse singular sets} and a resulting 
\emph{metric coarse singular homology}, and so on. We leave all of this to a 
future paper (or to the reader).

\begin{remark}
Mitchener has defined a related notion of \emph{coarse $n$-cells} and 
\emph{coarse $(n-1)$-spheres} (and resulting \emph{coarse $CW$-complexes}) 
\cites{MR1834777, MR2012966}. We will also defer the comparison of these with 
our coarse simplices (and resulting coarse simplicial complexes) to a future 
paper.
\end{remark}

The $l^1$ (or any $l^p$, $1 \leq p \leq \infty$) metric coarse structure on a 
$\Delta_n$ is the minimal ``good'' one, in the following sense. Fix $n \geq 0$, 
and consider the maps
\[
    \delta_{j_1}^m \circ \dotsb \circ \delta_{j_{n-m}}^{n-1}
        \from \Delta_m \to \Delta_n
\]
for all $0 \leq m < n$. (The $\delta_j$ all topologically embed their domains 
as closed subspaces of their codomains, and hence the same is true of 
compositions of the $\delta_j$.) Let us call the (set or topological) images of 
the each of the above maps a \emph{boundary simplex} of the topological space 
$\Delta_n$. We will not prove the following in full detail.

\begin{proposition}\label{prop:met-simp-univ}
Suppose $|\Delta_n|_\calR$ is a Roe coarse space with underlying topological 
space $\Delta_n \defeq (\setRplus)^n$ and Roe coarse structure $\calR$. Then 
there is a Roe coarse map $i \from |\Delta_n|_\TXTmet \to |\Delta_n|_\calR$ 
such that (as a set map) $i$ maps each boundary simplex of $\Delta_n$ to 
itself.
\end{proposition}

In fact, with a bit more trouble, one can even take $i$ to be a homeomorphism. 
The obvious discrete version of the above, with $(\setZplus)^n$ in place of 
$\Delta_n \defeq (\setRplus)^n$, is rather trivial. To get a nontrivial 
version, one should replace $|\Delta_n|_\calR$ with a ``sector'' which grows 
arbitrarily quickly away from the origin.

\begin{proof}[Sketch of proof]
It is trivial for $n = 0$, so suppose that $n \geq 1$. Fix an open 
neighbourhood $E_0 \in \calR$ of the diagonal $1_{\Delta_n}$. We will say that 
$B \subseteq \Delta_n$ is \emph{$E_0$-bounded} if $B^{\cross 2} \subseteq E_0$. 
In the following, \emph{disc} will mean ``closed $l^1$ metric disc in 
$\Delta_n$''; the \emph{diameter} of a disc will always be measured in the 
$l^1$ metric.

Tesselate $\Delta_n$ by discs diameter $1$ as in 
Figure~\ref{prop:met-simp-univ:fig-I} (we illustrate the case $n = 2$), and let
\[
    L_{2j} \defeq \set{x \in \Delta_n \suchthat j \leq \norm{x}_1 \leq j+1}
\]
for $j = 0, 1, 2, \dotsc$ be the ``layers'' of the tesselation. Then there is a 
refinement of this tesselation by discs as in 
Figure~\ref{prop:met-simp-univ:fig-II} such that each ``small'' disc of the 
refinement is $E_0$-bounded; label the layers of this tesselation $L'_{2j_0}, 
L'_{2j_1}, \dotsc$ as indicated in the Figure.

\begin{figure}
\resizebox{0.95\linewidth}{!}{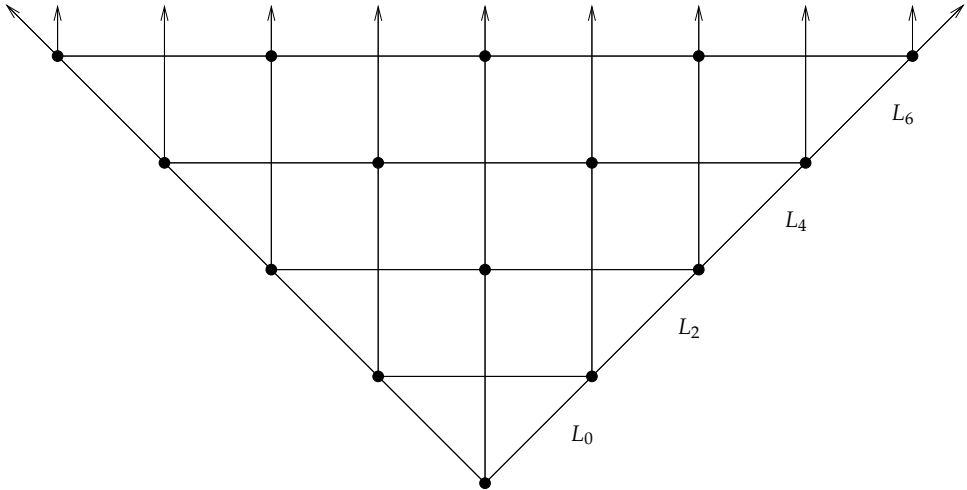}
\caption{\label{prop:met-simp-univ:fig-I}%
The tesselation of $\Delta_2$ by discs of $l^1$-diameter $1$.}
\end{figure}

\begin{figure}
\resizebox{0.95\linewidth}{!}{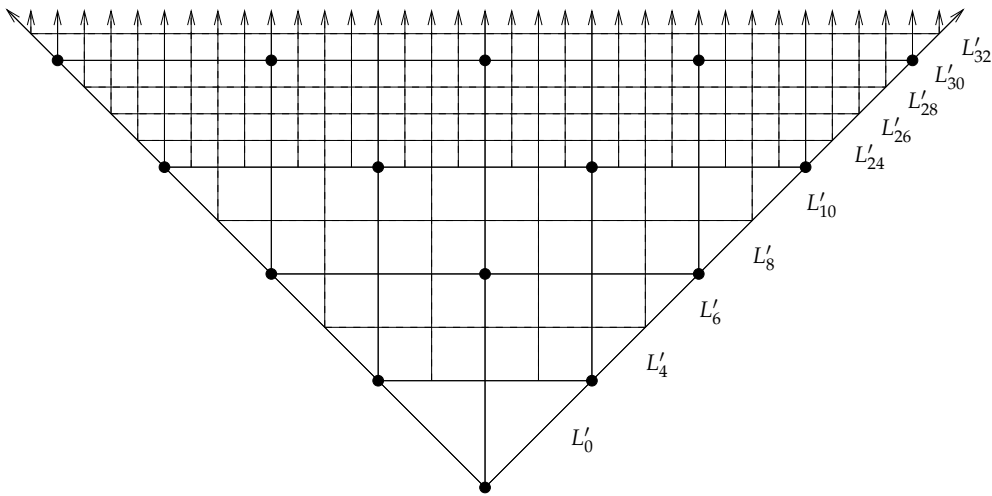}
\caption{\label{prop:met-simp-univ:fig-II}%
A refinement of the tesselation by $E_0$-controlled discs.}
\end{figure}

Define a continuous, ``tesselation preserving'' map $i \from \Delta_n \to 
\Delta_n$ which sends $L_{2j_0}$ to $L'_{2j_0}$, $L_{2j_1}$ to $L'_{2j_1}$, 
etc., collapsing the $L_{2j}$ which do not occur in the sequence $L_{2j_0}, 
L_{2j_1}, \dotsc$; in the example illustrated in the Figures, $L_2$ is 
collapsed to the ``level set'' $\set{x \in \Delta_2 \suchthat \norm{x}_1 = 1}$, 
$L_{12}$ through $L_{22}$ are collapsed to $\set{x \in \Delta_2 \suchthat 
\norm{x}_1 = 3}$, etc.

The map $i$ is (Roe) coarse. This map is proper and ``preserves'' the boundary 
simplices. Consider the cover $\set{B_1, B_2, \dotsc}$ of $\Delta_n$ by 
(overlapping) discs $B_k$ of diameter $2$, each a union of $2^n$ adjacent discs 
in the tesselation of Figure~\ref{prop:met-simp-univ:fig-I}. We have that
\[
    \bigunion_{k=1}^\infty (B_k)^{\cross 2}
\]
generates the Roe coarse structure of $|\Delta_n|_\TXTmet$. The collection 
$\set{i(B_1), i(B_2), \dotsc}$ is uniformly $(E_0 \circ E_0)$-bounded: the 
collection of unions of $2^n$, adjacent, \emph{equal-sized} discs in the 
tesselation of Figure~\ref{prop:met-simp-univ:fig-II} is uniformly $(E_0 \circ 
E_0)$-bounded, and each $i(B_k)$ is contained in such a disc (there are four 
cases to check in the latter assertion: (1) $i$ does not collapse $B_k$ at all, 
(2) $i$ completely collapses $B_k$, (3) $i$ collapses the ``top half'' of 
$B_k$, or (4) $i$ collapses the ``bottom half'' of $B_k$). This suffices to 
show that $i$ preserves all (Roe) entourages of $|\Delta_n|_\TXTmet$.
\end{proof}

%FIXME -- ref to coarse homotopy?
\begin{remark}
The above Proposition is not entirely satisfactory. $|\Delta_n|_\TXTmet$ should 
satisfy a stronger universal property (which I have not yet proven): 
$|\Delta_n|_\TXTmet$ should be \emph{coarse-homotopy}-universal with the above 
property. That is, if $\calS$ is any Roe coarse structure on $\Delta_n$ such 
that the above is true with $|\Delta_n|_\calS$ in place of 
$|\Delta_n|_\TXTmet$, then $|\Delta_n|_\calS$ should be coarse homotopy 
equivalent to $|\Delta_n|_\TXTmet$ in such a way that its boundary simplices 
are preserved (compare \cite{MR1243611}*{Thm.~7.3}).
\end{remark}

\subsection{Continuously controlled coarse simplices}

If the previously defined metric coarse structure on a simplex $\Delta_n$ is 
the minimal ``good'' one, the continuously controlled coarse structure on 
$\Delta_n$ defined below is the maximal ``good'' one (again, we will not make 
this precise in this paper).

For $n = 1, 2, 3, \dotsc$, let $\overline{\Delta}_n$ be the obvious 
compactification of the topological space $\Delta_n \defeq (\setRplus)^n$ by 
the standard topological simplex of geometric dimension $n-1$. Alternatively 
(and equivalently, for our purposes), put
\begin{align*}
    \Delta_n & \defeq \Bigset{(x_0, \dotsc, x_{n-1}) \in (\setRplus)^n
            \suchthat \sum_{j=0}^{n-1} x_j < 1}
\quad\text{and} \\
    \overline{\Delta}_n & \defeq \Bigset{(x_0, \dotsc, x_{n-1})
            \in (\setRplus)^n \suchthat \sum_{j=0}^{n-1} x_j \leq 1},
\end{align*}
so that $\die \Delta_n \defeq \overline{\Delta}_n \setminus \Delta_n$ really is 
the standard topological $(n-1)$-simplex. Put $\Delta_0 \defeq \set{0}$ which 
is compact, so $\overline{\Delta}_0 = \Delta_0$ and $\die \Delta_0 = 
\emptyset$.

\begin{definition}
For $n = 0, 1, 2, \dotsc$, the \emph{continuously controlled coarse 
$n$-simplex} is the continuously controlled coarse space
\[
    |\Delta_n| \defeq |\Delta_n|_\TXTtop \defeq |\Delta_n|_{\die \Delta_n}.
\]
\end{definition}

Equivalently (see Remark~\ref{rmk:ctsctl-cones}), we can define 
$|\Delta_n|_\TXTtop$ to be the continuously controlled open cone 
$\calO_\TXTtop(\die \Delta_n)$ (with underlying set $\setRplus \cross (\die 
\Delta_n)$).

Again, as in \S\ref{subsect:met-simpl}, we can define various coarse maps 
$\delta_j$ and $\sigma_j$ between the continuously controlled coarse simplices. 
Indeed (using either of the above descriptions of the $\Delta_n$), we may 
define them using the same formul\ae\ \eqref{subsect:met-simpl:eq:delta} and 
\eqref{subsect:met-simpl:eq:sigma}, and hence they also satisfy the 
cosimplicial identities. Consequently, we get a \emph{continuously controlled 
coarse simplex functor}
\[
    |\Delta|_\TXTtop \from \CATSimp \to \CATCrs,
\]
and everything that comes along with it: \emph{continuously controlled coarse 
realizations} of simplicial sets, a notion of (continuously controlled coarse) 
``weak equivalence'', \emph{continuously controlled coarse singular sets} and 
\emph{homology}, etc.

\begin{remark}%\label{rmk:ctsctl-singhom-iso}
If $X = \calO_\TXTtop M$ for a second countable compact topological space $M$ 
(where $\calO_\TXTtop M$ is the continuously controlled open cone on $M$ from 
Def.~\ref{def:ctsctl-cone}), then it is easy to see that the continuously 
controlled coarse singular homology of $\calO_\TXTtop M$ is exactly the 
singular homology of $M$ (in this case, we would want to discard our 
$0$-simplices and shift our indexing to match the topologists').
Continuously controlled coarse simplices have another nice feature: 
$|\Delta_1|_\TXTtop$ is the continuously controlled ray, which is coarsely 
equivalent to $|\setZplus|_1$, so $|\Delta_1|_\TXTtop$ is a product identity 
for most coarse spaces which arise in practice (those in $\CATCrs_{\preceq 
|\Delta_1|_\TXTtop}$, which includes all those which are coarsely equivalent to 
countable coarse spaces). However, continuously controlled simplices have a 
fundamental problem: they are too coarse, and so many coarse spaces $X$ of 
interest (e.g., metric coarse spaces) do not even admit a coarse map 
$|\Delta_1|_\TXTtop \to X$.
\end{remark}

\subsection{\pdfalt{\maybeboldmath $\sigma$-coarse spaces and $\sigma$-unital
        coarse spaces}{sigma-coarse spaces and sigma-unital coarse spaces}}

In \cite{MR2225040}*{\S{}2}, Emerson--Meyer consider increasing sequences of 
coarse spaces. Their coarse spaces are equipped with topologies and are 
connected and unital (i.e., are \emph{Roe coarse spaces} in the terminology of 
\S{}\ref{subsect:Roe-crs-sp}). We will simply handle the discrete case. (This 
is perhaps at significant loss of generality, since in a sense Emerson--Meyer 
are largely interested in ``non-locally-compact coarse spaces'' which we do not 
really examine in this paper; see Remark~\ref{rmk:top-crs-sp}.) For our 
purposes, we may safely discard the connectedness assumption, though we still 
need unitality.

\begin{definition}[\cite{MR2225040}*{\S{}2}]
A (discrete) \emph{$\sigma$-coarse space} $(X_m)$ is a nondecreasing sequence
\[
    X_0 \subseteq X_1 \subseteq X_2 \subseteq \dotsb
\]
of unital coarse spaces such that, for all $m \geq 0$, $X_m$ is a coarse 
subspace of $X_{m+1}$ (i.e., is a subset and has the subspace coarse 
structure).
\end{definition}

\begin{remark}
Given a sequence $(X_m)$ which is a $\sigma$-coarse space in the sense of 
Emerson--Meyer (i.e., each $X_m$ is a Roe coarse space thus may have nontrivial 
topology), one can obtain a nondecreasing sequence of coarse spaces by applying 
our discretization functor $\Disc$ to each $X_m$. However, $\Disc(X_m)$ is 
typically not unital. It may be interesting to remove the unitality assumption 
from the above Definition, and thus be able to consider $(\Disc(X_m))$ as a 
``nonunital $\sigma$-coarse space''.
\end{remark}

Until otherwise stated (near the end of this section), $(X_m)$ and $(Y_n)$ will 
always denote $\sigma$-coarse spaces.

\begin{definition}[\cite{MR2225040}*{\S{}4}]
A \emph{coarse map} $(f_n) \from (Y_n) \to (X_m)$ of $\sigma$-coarse spaces is 
a map of directed systems in $\CATPCrs$ (taken modulo cofinality).
\end{definition}

That is, a coarse map $(f_n) \from (Y_n) \to (X_m)$ is represented by a 
sequence of coarse maps $f_n \from Y_n \to X_{m(n)}$, $n = 0, 1, \dotsc$, 
(where $0 \leq m(0) \leq m(1) \leq \dotsb$ is a nondecreasing sequence) such 
that the obvious diagram commutes (in $\CATPCrs$, not modulo closeness); two 
representative sequences $(f_n)$, $(f'_n)$ are considered to be equivalent if, 
for all $n$, the compositions
\begin{equation}\label{sect:sigma-crs:eq:maps}
    Y_n \nameto{\smash{f_n}} X_{m(n)} \injto X_{\max\set{m(n),m'(n)}}
\quad\text{and}\quad
    Y_n \nameto{\smash{f'_n}} X_{m(n)} \injto X_{\max\set{m(n),m'(n)}}
\end{equation}
are equal.

Actually, Emerson--Meyer consider maps $\bigunion_n Y_n \to \bigunion_m X_m$, 
i.e., maps between set colimits which restrict to give sequences of coarse 
maps. This is equivalent to our definition (which avoids set colimits).

\begin{definition}[\cite{MR2225040}*{\S{}4}]\label{def:sigma-crs}
Coarse maps $(f_n), (f'_n) \from (Y_n) \to (X_m)$ are \emph{close} if, for all 
$n$ (and any, hence all, representative sequences $(f_n)$, $(f'_n)$, 
respectively), the compositions \eqref{sect:sigma-crs:eq:maps} are close. We 
denote the \emph{closeness} (equivalence) \emph{class} of $(f_n)$ by $[f_n]$.
\end{definition}

Equivalently, coarse maps $(f_n)$, $(f'_n)$ are close if they yield maps of 
directed systems in $\CATCrs$ which are equivalent modulo cofinality.

Since the system $X_0 \to X_1 \to \dotsb$ consists of inclusion maps, the 
precoarse colimit $\pfx{\CATPCrs}\OBJcolim X_m$ exists; one may take it to be
\[
    X \defeq \pfx{\CATPCrs}\OBJcolim X_m \defeq \bigunion_m X_m
\]
as a set, with coarse structure
\[
    \calE_X \defeq \langle \calE_{X_m} \suchthat m = 0, 1, \dotsc \rangle_X
\]
generated by the coarse structures of all the $X_m$. In fact, since $X_m$ is a 
coarse subspace of $X_{m+1}$ for all $m$,
\[
    \calE_X = \bigunion_m \calE_{X_m}
\]
(and $X_m$ is a subspace of $X$); conversely, we get, for each $m$, that
$\calE_{X_m} = \calE_X |_{X_m}$.

Until otherwise stated, let $X$ be as above and similarly $Y \defeq 
\pfx{\CATPCrs}\OBJcolim Y_m \defeq \bigunion_n Y_n$.

The coarse colimit $\pfx{\CATCrs}\OBJcolim X_m$ also exists (since all colimits 
in $\CATCrs$ exist), and maps canonically to $X$ in $\CATCrs$. The following is 
easy to show.

\begin{proposition}
$\pfx{\CATCrs}\OBJcolim X_m = \pfx{\CATPCrs}\OBJcolim X_m \eqdef X$. More 
precisely, the canonical arrow
\[
    \pfx{\CATCrs}\OBJcolim X_m \to \pfx{\CATPCrs}\OBJcolim X_m \eqdef X
\]
is an isomorphism (in $\CATCrs$).
\end{proposition}

By definition, any coarse map $(f_n) \from (Y_n) \to (X_m)$ of $\sigma$-coarse 
spaces yields a well-defined coarse map $f \from Y \to X$. (Of course, $f$ is 
just, as a set map, given by $f(y_n) \defeq f_n(y_n)$ for all $n$ and $y_n \in 
Y_n$.) Likewise, its closeness class $[f_n]$ yields a well-defined closeness 
class $[f] \from Y \to X$.

Let $\calP\calS$ be the category of $\sigma$-coarse spaces and coarse maps, and 
$\calS$ be the category of $\sigma$-coarse spaces and closeness classes of 
coarse maps. We have defined functors
\[
    \calL \defeq \pfx{\CATPCrs}\OBJcolim \from \calP\calS \to \CATPCrs
\quad\text{and}\quad
    [\calL] \defeq \pfx{\CATCrs}\OBJcolim \from \calS \to \CATCrs.
\]

\begin{proposition}\label{prop:sigma-crs:L-full-faith}
The functor $\calL \from \calP\calS \to \CATPCrs$ is fully faithful.
\end{proposition}

(Recall that ``faithful'' does not require injectivity on object sets!)

\begin{proof}
Faithfulness: Clear, since representative sequences $(f_n), (f'_n) \from (Y_n) 
\to (X_m)$ are cofinally equivalent if and only if they are equal on colimits 
(i.e., $f = f'$).

Fullness: To show that $\calL$ maps $\Hom_{\calP\calS}((Y_n),(X_m))$ to 
$\Hom_{\CATPCrs}(Y,X)$ surjectively, we must use the unitality of the $Y_n$. 
Suppose $f \from Y \to X$ is a coarse map (not a priori in the image of 
$\calL$). For each $n$, $Y_n$ is a unital subspace of $Y$, and hence $f(Y_n)$ 
is a unital subspace of $X$. Then $1_{f(Y_n)}$ must be an entourage of some 
$X_m$; let $m(n)$ be the least such $m$. Since $\calE_{X_{m(n)}} = \calE_X 
|_{X_{m(n)}}$, $f_n \defeq f |_{Y_n}^{X_{m(n)}} \from Y_n \to X_{m(n)}$ is a 
coarse map. It follows that $(f_n)$ is a coarse map of $\sigma$-coarse spaces, 
and that $\calL((f_n)) = f$.
\end{proof}

The following shows that unitality of the $Y_n$ really is needed for fullness.

\begin{example}\label{ex:sigma-crs:L-nonunital-not-full}
Put, for each $m$, $X_m \defeq |\set{0, \dotsc, m-1}|_1$, so that $X \defeq 
\calL((X_m))$ is just $\setZplus$ as a set, with entourages the finite subsets 
of $(\setZplus)^{\cross 2}$. Put, for all $n$, $Y_n \defeq X$, so that colimit 
$Y \defeq X$ is nonunital ($(Y_n)$ is not a $\sigma$-coarse space in our 
terminology). The identity map $Y \to X$ is coarse, but its image is not 
contained in any single $X_m$ so is no ``coarse map'' $(f_n) \from (Y_n) \to 
(X_m)$ which yields $f$.
\end{example}

A $\sigma$-coarse space $(X_m)$ includes as a part of its structure the 
``filtration'' $X_0 \subseteq X_1 \subseteq \dotsb$. However, the particular 
choice of ``filtration'' is not important, since maps of $\sigma$-coarse spaces 
are taken modulo cofinality.

\begin{corollary}
If $X \defeq \calL((X_m))$ is isomorphic in $\CATPCrs$ to $Y \defeq 
\calL((Y_n))$ (i.e., there is a \emph{bijection} of sets $f \from Y \to X$ such 
that $f$ and $f^{-1}$ are both coarse maps), then $(X_m)$ is isomorphic to 
$(Y_n)$ in $\calP\calS$ (in particular, this is the case when $X = Y$ as coarse 
spaces).
\end{corollary}

The situation modulo closeness parallels the above.

\begin{proposition}\label{prop:sigma-crs:QL-full-faith}
The functor $[\calL] \from \calS \to \CATCrs$ is fully faithful.
\end{proposition}

\begin{proof}
Faithfulness: Since each $Y_n$ is a subspace of $Y$ and each $X_m$ a subspace 
of $X$, closeness of $f = \calL((f_n))$ to $f' = \calL((f'_n))$ implies 
closeness of the compositions \eqref{sect:sigma-crs:eq:maps} (noting that $f_n 
= f |_{Y_n}^{X_{m(n)}}$ and similarly for $f'_n$).

Fullness: Here, we implicitly use the unitality condition. We have a 
commutative diagram
\[\begin{CD}
    \calP\calS @>{\calL}>> \CATPCrs \\
    @V{\Quotient}VV @V{\Quotient}VV \\
    \calS @>{[\calL]}>> \CATCrs
\end{CD}\,.\]
Since $\calL$ is full and evidently the quotient functors are also full and map 
surjectively onto object sets, $[\calL]$ is full.
\end{proof}

It is not clear to me whether fullness of $[\calL]$ fails if the unitality 
condition is removed from Definition~\ref{def:sigma-crs}; the counterexample of 
Example~\ref{ex:sigma-crs:L-nonunital-not-full} fails.

\begin{corollary}
If $X \defeq \calL((X_m))$ is coarsely equivalent (i.e., isomorphic in 
$\CATCrs$) to $Y \defeq \calL((Y_n))$, then $(X_m)$ is isomorphic to $(Y_n)$ in 
$\calS$.
\end{corollary}

It follows from Propositions \ref{prop:sigma-crs:QL-full-faith} 
and~\ref{prop:sigma-crs:L-full-faith} that $\calL$ and $[\calL]$ are 
equivalences (of categories) onto their images. We now consider what the images 
of these functors are (and how one constructs ``inverse'' functors).

Let us ``reset'' our notation: $X$, $Y$ are just coarse spaces, not necessarily 
coming from $\sigma$-coarse spaces, and $(X_m)$, $(Y_n)$ are not assumed to 
have any meaning.

\begin{definition}
A coarse space $X$ is \emph{$\sigma$-unital} if there is a nondecreasing 
sequence
\[
    X_0 \subseteq X_1 \subseteq \dotsb \subseteq X
\]
of unital subspaces of $X$ such that each unital subspace $X' \subseteq X$ is 
contained in some $X_m$ ($m$ depending on $X'$).
\end{definition}

It is implied that $X = \bigunion_m X_m$, though this equality certainly does 
not imply that each unital subspace of $X$ is contained in some $X_m$.

Let $\CATPCrs_{\bfsigma} \subseteq \CATPCrs$ and $\CATCrs_{\bfsigma} \subseteq 
\CATCrs$ denote the full subcategories of $\sigma$-unital coarse spaces. 
Clearly, $\calL$ and $[\calL]$ map both into and onto $\CATPCrs_{\bfsigma}$ and 
$\CATCrs_{\bfsigma}$, respectively. We get the following.

\begin{theorem}
The functors $\calL \from \calP\calS \to \CATPCrs_{\bfsigma}$ and $[\calL] 
\from \calS \to \CATCrs_{\bfsigma}$ are equivalences of categories.
\end{theorem}

It is also easy to construct ``inverse'' functors. Choose, for each 
$\sigma$-unital $X$, a ``filtration'' $(X_m)$. Then $X \mapsto (X_m)$ (and, for 
$f \from Y \to X$, $f \mapsto (f_n)$, where $f_n$ is an appropriate range 
restriction of $f |_{Y_n}$) gives a functor ``inverse'' to $\calL \from 
\calP\calS \to \CATPCrs_{\bfsigma}$. Choosing representative coarse maps, one 
does the same to obtain an ``inverse'' to $[\calL] \from \calS \to 
\CATCrs_{\bfsigma}$.

\subsection{Quotients and Roe algebras}\label{subsect:quot-Cstar}

We shall assume that the reader is familiar with the definition and 
construction of the Roe algebras $C^*(X)$ for $X$ a (Roe) coarse space (see, 
e.g., \cite{MR1451755}); the generalization to our nonunital situation is 
straightforward. We will follow the standard, abusive practice of pretending 
that $X \mapsto C^*(X)$ is a functor. (The situation is slightly complicated by 
our nonunital situation. However, there are a number of ways of obtaining an 
actual functor, just not to the category of $C^*$-algebras. One could, for 
example, construct a coarsely invariant functor from $\CATCrs$ to the category 
of $C^*$-categories \cite{MR1881396}.) The important fact is that, applying 
$K$-theory, one gets a coarsely invariant functor $X \mapsto 
K_\grstar(C^*(X))$. The following should be regarded as a sketch, with more 
details to follow in a future paper.

Fix a coarse space $X$ and a subspace $Y \subseteq X$, and denote the inclusion 
$Y \injto X$ by $\iota$. We note that the following does not depend on our 
generalizations, and even works in the ``classical'' unital context; if $X$ is 
a Roe coarse space in the sense of \S\ref{sect:top-crs}, $Y$ should be closed 
in $X$. Recall that we simply denote the quotient $X/[\iota](Y)$ (defined in 
\S\ref{subsect:Crs-quot}) by $X/[Y]$. The coarse space $X/[Y]$ is easy to 
describe: It is just $X$ as a set, with coarse structure generated by the 
entourages of $X$ and those of $\Terminate(Y)$ (if $X$ is unital, the latter 
are just those of the terminal coarse structure on $Y$).

The quotient $Y/[Y] \Terminate(Y)$ is a subspace of $X/[Y]$, with
$\utilde{\iota} \from Y/[Y] \injto X/[Y]$ an inclusion. We get a commutative 
square
\[\begin{CD}
    Y @>{\iota}>> X \\
    @V{\tilde{q}}VV @V{q}VV \\
    Y/[Y] @>{\utilde{\iota}}>> X/[Y]
\end{CD}\quad,\]
where $\tilde{q}$ and $q$ represent the quotient maps (which one can take to be 
identity set maps). This square gives rise to a commutative diagram
\[\begin{CD}
    0 @>>> C^*_X(Y) @>{\iota_*}>> C^*(X) @>>> Q_{X,Y} @>>> 0 \\
    @. @V{\tilde{q}_*}VV @V{q_*}VV @V{\utilde{q}_*}VV \\
    0 @>>> C^*_{X/[Y]}(Y/[Y]) @>{\utilde{\iota}_*}>> C^*(X/[Y])
        @>>> Q_{X/[Y],Y/[Y]} @>>> 0
\end{CD}\]
of $C^*$-algebras whose rows are exact; $C^*_X(Y)$ denotes the ideal of 
$C^*(X)$ of operators supported near $Y$ (which can be identified with 
$C^*(\OBJcoim [\iota])$, where $\OBJcoim [\iota] = X$ as a set with the 
nonunital coarse structure of entourages of $X$ supported near $Y$; see 
Def.~\ref{def:Crs-coimage}) and $Q_{X,Y}$ is the quotient $C^*$-algebra (and 
similarly for the second row).

Next, one observes that $\utilde{q}_*$ is an isomorphism \emph{of 
$C^*$-algebras} hence induces an isomorphism on $K$-theory. Let us specialize 
to the case when $X$ is unital (from which it follows that $Y$ and the quotient 
coarse spaces are also unital), and examine the consequences. If $Y$ is finite 
(or compact, in the Roe coarse space version) then $X = X/[Y]$ and $Y = Y/[Y]$, 
so $\tilde{q}_*$ and $q_*$ are identity maps on the level of $C^*$-algebras and 
hence the diagram is trivial.

On the other hand, if $Y$ is infinite, then $Y/[Y] = |Y|_1$ has the terminal 
coarse structure and one can show by a standard ``Eilenberg swindle'' (see, 
e.g., \cite{MR1817560}*{Lem.~6.4.2}) that $K_\grstar(C^*_{X/[Y]}(Y/[Y])) = 0$. 
Thus we get a canonical isomorphism
\[
    K_\grstar(C^*(X/[Y])) \isoto K_\grstar(Q_{X/[Y],Y/[Y]})
        \isoto K_\grstar(Q_{X,Y})
\]
on $K$-theory. Consequently, using the isomorphism $K_\grstar(C^*(Y)) \isoto
K_\grstar(C^*_X(Y))$ (which is easy to prove under most circumstances), we get 
a long (or six-term) exact sequence
\begin{equation}\label{subsect:quot-Cstar:eq:lx}
    \dotsb \nameto{\smash{\die}} K_\grstar(C^*(Y))
        \to K_\grstar(C^*(X))
        \to K_\grstar(C^*(X/[Y]))
        \nameto{\smash{\die}} K_{\grstar-1}(C^*(Y))
        \to \dotsb \,.
\end{equation}

\begin{remark}[continuous control]\label{rmk:ctsctl-quot-Cstar}
In the above situation, suppose that $X = \calO M$ and $Y = \calO N$ are 
continuously controlled open cones, where $N$ is a nonempty closed subspace of 
a second countable, compact space $M$ and we abbreviate $\calO \defeq 
\calO_\TXTtop$. Then there are natural isomorphisms
\begin{equation}\label{rmk:ctsctl-quot-Cstar:eq}
    K_\grstar(C^*(\calO M)) \cong \tilde{K}^{1-\grstar}(C(M))
        = \tilde{K}_{\grstar-1}(M)
\quad\text{and}\quad
    K_\grstar(C^*(\calO N)) \cong \tilde{K}_{\grstar-1}(N),
\end{equation}
where $\tilde{K}$ is reduced $K$-homology (see, e.g., 
\cite{MR1817560}*{Cor.~6.5.2}). One can check that there is also a natural 
isomorphism
\[
    K_\grstar(C^*(\calO M/[\calO N])) \cong K_{\grstar-1}(M,N)
\]
(to relative $K$-homology), so that the above long exact sequence 
\eqref{subsect:quot-Cstar:eq:lx} naturally maps isomorphically to the reduced 
$K$-homology sequence
\[
    \dotsb \nameto{\smash{\die}} \tilde{K}_{\grstar-1}(N)
        \to \tilde{K}_{\grstar-1}(M)
        \to K_{\grstar-1}(M,N)
        \nameto{\smash{\die}} K_{\grstar-2}(N)
        \to \dotsb \,.
\]
We have three natural isomorphisms
\begin{align*}
    K_\grstar(C^*(\calO M/[\calO N])) & \cong K_\grstar(C^*(\calO (M/N))), \\
    K_\grstar(C^*(\calO (M/N))) & \cong \tilde{K}_{\grstar-1}(M/N),
\qquad\qquad\quad\text{and} \\
    K_{\grstar-1}(M,N) & \cong \tilde{K}_{\grstar-1}(M/N),
\end{align*}
from Remark~\ref{rmk:ctsctl-quot}, as in \eqref{rmk:ctsctl-quot-Cstar:eq} 
above, and by excision for $K$-homology, respectively; these are mutually 
compatible in the obvious sense.
\end{remark}

\begin{example}[\maybeboldmath $K$-theory of $\calO_\TXTtop S^n$]
We give yet another version of a standard calculation (see, e.g., 
\cite{MR1817560}*{Thm.~6.4.10}). For $n \geq 0$, denote the topological 
$n$-sphere by $S^n$ and, for $n \geq 1$, the closed $n$-disc by $D^n$; recall 
that $D^n$ has ``boundary'' $S^{n-1}$ and that $D^n/S^{n-1} \cong S^n$. Again 
we abbreviate $\calO \defeq \calO_\TXTtop$.

First, we compute the $K$-theory of $X \defeq \calO S^0$. Put $Y \defeq \calO 
\set{-1} \subseteq X$, $X' \defeq \calO \set{1} \subseteq X$, and $Y' \defeq 
\set{0} \subseteq Y \intersect X'$. It is well known that $K_\grstar(C^*(X')) = 
0$ and $K_\grstar(C^*_X(Y)) = 0$ (by the aforementioned ``Eilenberg swindle''), 
and that
\[
    K_\grstar(C^*_{X'}(Y')) = \begin{cases}
                \setZ & \text{if $\grstar \equiv 0 \AMSdisplayoff\pmod{2}$,
                        and} \\
                0 & \text{otherwise}
            \end{cases}
\]
(since $C^*_{X'}(Y')$ is just the compact operators). We have a map of short 
exact sequences
\[\begin{CD}
    0 @>>> C^*_{X'}(Y') @>>> C^*(X') @>>> Q' @>>> 0 \\
    @. @VVV @VVV @VVV \\
    0 @>>> C^*_{X}(Y) @>>> C^*(X) @>>> Q @>>> 0
\end{CD}\quad.\]
But one checks easily that the map $Q' \to Q$ is an isomorphism of 
$C^*$-algebras, hence from the $K$-theory long exact sequences we get
\begin{multline*}
    K_\grstar(C^*(\calO S^0)) = K_\grstar(Q) = K_\grstar(Q') \\
        = K_{\grstar-1}(C^*_{X'}(Y')) = \begin{cases}
                \setZ & \text{if $\grstar \equiv 1 \AMSdisplayoff\pmod{2}$,
                        and} \\
                0 & \text{otherwise.}
            \end{cases}
\end{multline*}

We proceed to calculate the $K$-theory of $\calO S^n$, $n \geq 1$, by 
induction. Put $X \defeq \calO D^n$ and $Y \defeq \calO S^{n-1} \subseteq X$, 
and recall that $X/[Y] = \calO (D^n/S^{n-1}) = \calO S^n$. Then, by 
Remark~\ref{rmk:ctsctl-quot-Cstar} above, we have a long exact sequence
\[
    \dotsb \nameto{\smash{\die}} K_\grstar(C^*(Y))
        \to K_\grstar(C^*(X))
        \to K_\grstar(C^*(\calO S^n))
        \nameto{\smash{\die}} K_{\grstar-1}(C^*(Y))
        \to \dotsb \,.
\]
By another ``Eilenberg swindle'', one shows that $K_\grstar(C^*(X)) = 0$ and 
hence
\[
    K_\grstar(C^*(\calO S^n)) = K_{\grstar-1}(C^*(Y))
        = \begin{cases}
                \setZ & \text{if $\grstar \equiv n-1 \AMSdisplayoff\pmod{2}$,
                        and} \\
                0 & \text{otherwise.}
            \end{cases}
\]
\end{example}

\begin{example}[\maybeboldmath suspensions in $K$-homology]%
        \label{ex:ctsctl-quot-Khom}
Suppose that $A$ is a separable $C^*$-algebra. It is known that, for $n \geq 
0$, elements of the Kasparov $K$-homology group $K^{n+1}(A)$ can be represented 
by (equivalence classes of) $C^*$-algebra morphisms
\begin{equation}\label{ex:ctsctl-quot-Khom:eq:A}
    \phi \from A \to C^*(\calO S^n)
\end{equation}
(see \cites{MR1627621, my-thesis}; I caution that, in my opinion, this is 
probably not the ``best'' coarse geometric description of $K$-homology, but 
work remains ongoing). The pairing of $K_m(A)$ with a $K$-homology class 
represented by such $\phi$ is given simply by applying $K$-theory to $\phi$ 
(and using the computation as in the previous Example).

Fix $n \geq 1$ and suppose that we are given an element of $K^n(\Sigma A)$, 
where $\Sigma A \defeq C_0(\ooitvl{0,1}) \tensor A$ is the $C^*$-algebraic 
suspension of $A$, represented by a morphism
\begin{equation}\label{ex:ctsctl-quot-Khom:eq:SA}
    \tilde{\psi} \from \Sigma A \to C^*(\calO S^{n-1}).
\end{equation}
Actually, let us assume something stronger, that we are given a morphism $\psi$ 
which fits into the following commutative diagram whose \emph{rows} are with 
exact:
\[\begin{CD}
    0 @>>> \Sigma A @>>> CA @>>> A @>>> 0 \\
    @. @V{\psi |_{\Sigma A}}VV @V{\psi}VV @VVV \\
    0 @>>> C^*_X(Y)
        @>>> C^*(X)
        @>>> Q @>>> 0 \\
    @. @VVV @VVV @VVV \\
    0 @>>> C^*_{X/[Y]}(Y/[Y])
        @>>> C^*(X/[Y])
        @>>> Q' @>>> 0
\end{CD}\quad,\]
where $CA \defeq C_0(\coitvl{0,1}) \tensor A$ is the cone on $A$, $X \defeq 
\calO D^n$, and $Y \defeq \calO S^{n-1} \subseteq X$. (In fact, given a 
$\tilde{\psi}$, one can find a $\psi$ such that
\[\begin{CD}
    K_\grstar(\Sigma A) @>{\tilde{\psi}}>> K_\grstar(C^*(Y)) \\
    @V{=}VV @V{\sim}VV \\
    K_\grstar(\Sigma A) @>{\psi |_{\Sigma A}}>>
        K_\grstar(C^*_X(Y))
\end{CD}\]
commutes. This is not easy to prove, and seems to require that $A$ be 
separable.)

Denote the composition $A \to Q \to Q'$ by $\utilde{\phi}$. From the previous 
Example, we have natural isomorphisms
\[
    K_\grstar(\calO S^{n-1}) = K_\grstar(C^*_X(Y)) = K_{\grstar+1}(Q)
        = K_{\grstar+1}(Q') = K_{\grstar+1}(X/[Y]) = K_{\grstar+1}(\calO S^n).
\]
Moreover, since $K_\grstar(CA) = 0$, we have $K_\grstar(\Sigma A) = 
K_{\grstar+1}(A)$. These isomorphisms are all compatible, in the sense that 
$\psi$ and $\utilde{\phi}$ are naturally equivalent on $K$-theory (with a 
dimension shift).

In fact, one can ``lift'' the morphism $\utilde{\phi}$ to a morphism $\phi 
\from A \to C^*(X/[Y]) = C^*(\calO S^n)$ in the weak sense that the composition 
$A \nameto{\smash{\phi}} C^*(X/[Y]) \to Q'$ is equal to $\utilde{\phi}$ on the 
level of $K$-theory. (This is not too difficult, but again seems to require 
that $A$ be separable.) This provides a map from the $K$-homology group 
$K^n(\Sigma A)$ (described as classes of morphisms as in 
\eqref{ex:ctsctl-quot-Khom:eq:SA}) to the group $K^{n+1}(A)$ (described as in 
\eqref{ex:ctsctl-quot-Khom:eq:A}).
\end{example}

%%%%%%%%%%%%%%%%%%%%%%%%%%%%%%%%%%%%%%%%%%%%%%%%%%%%%%%%%%%%%%%%%%%%%%%%%%%%%%%%

%%%%%%%%%%%%%%%%%%%%%%%%%%%%%%%%%%%%%%%%%%%%%%%%%%%%%%%%%%%%%%%%%%%%%%%%%%%%%%%%

\end{document}